\documentclass[12pt,a4paper]{article}
\usepackage[left=2.0cm, right=2.0cm, top=2cm, bottom=2.5cm]{geometry}

\usepackage[utf8]{inputenc}
\usepackage[T1]{fontenc}
\usepackage{amsmath, amssymb, amsfonts}
\usepackage{amsthm}
\newtheorem{assumption}{Assumption}
\newtheorem{condition}{Condition}
\newtheorem{remark}[assumption]{Remark}
\usepackage{authblk}
\usepackage[numbers]{natbib}
\usepackage{graphicx}
\usepackage{hyperref}
\usepackage{tabularx}
\newcolumntype{C}{>{\centering\arraybackslash}X}
\usepackage{algorithm, algorithmicx, algpseudocode}

\usepackage{caption}
\usepackage[inline, shortlabels]{enumitem}
\usepackage{amsmath,amsthm}
\usepackage{xspace}
\usepackage{color} 
\usepackage{url} 
\usepackage{graphicx} 
\usepackage{thmtools}
\usepackage{booktabs}
\usepackage{cleveref}
\usepackage{amssymb}
\usepackage{pifont} 
\usepackage{makecell}

\newcommand{\calP}{\mathcal{P}}
\newcommand{\calY}{\mathcal{Y}}

\newcommand{\Reals}{\ensuremath{\mathbb{R}}}
\newcommand*{\NNReals}{\ensuremath{\mathbb{R}_{\geq 0}}}
\newcommand*{\ExtReals}{\ensuremath{\bar{\mathbb{R}}}}
\newcommand*{\PReals}{\ensuremath{\mathbb{R}_{\geq 0}}}

\newcommand{\EE}{\ensuremath{\mathbb{E}}}
\newcommand{\II}{\ensuremath{\mathbb{I}}}

\newcommand{\hatS}{\hat{S}}
\newcommand{\hatF}{\hat{F}}

\newcommand{\etal}{\textit{et al.}\xspace}

\declaretheorem[within=section,name=Definition]{Def}
\declaretheorem[sibling=Def,name=Proposition]{Prop}
\declaretheorem[sibling=Def,name=Lemma]{Lem}

\crefname{Lem}{lemma}{lemmas}
\Crefname{Lem}{Lemma}{Lemmas}
\crefname{Def}{definition}{definitions}
\Crefname{Def}{Definition}{Definitions}
\crefname{Cor}{corollary}{corollaries}
\Crefname{Cor}{Corollary}{Corollaries}
\crefname{Prop}{proposition}{propositions}
\Crefname{Prop}{Proposition}{Propositions}	

\newcommand{\indep}{\perp \!\!\! \perp}

\newcommand{\Exp}{\operatorname{Exp}}

\title{When Are Scoring Rules Proper? Bridging Theory and Practice in Survival Model Evaluation}

\author[1,2]{John Zobolas*}
\author[3]{Raphael Sonabend*}
\author[4]{Riccardo De Bin}
\author[5,6,7]{Johannes Piller}
\author[5,6]{Philipp Kopper}
\author[5,6,8,9]{Lukas Burk}
\author[5,6]{Andreas Bender}

\affil[1]{Department of Cancer Genetics, Institute for Cancer Research, Oslo University Hospital, Norway}
\affil[2]{Department of Biostatistics, Faculty of Medicine, Oslo Centre for Biostatistics and Epidemiology (OCBE), University of Oslo, Norway}
\affil[3]{OSPO Now, London, UK}
\affil[4]{Department of Mathematics, Faculty of Mathematics and Natural Sciences, University of Oslo, Norway}
\affil[5]{Department of Statistics, LMU Munich, Germany}
\affil[6]{Munich Center for Machine Learning (MCML), LMU Munich, Germany}
\affil[7]{Statistical Consulting Unit (StaBLab), LMU Munich, Germany}
\affil[8]{Leibniz Institute for Prevention Research and Epidemiology - BIPS, Bremen, Germany}
\affil[9]{Faculty of Mathematics and Computer Science, University of Bremen, Germany}

\date{} 

\begin{document}
\footnotetext{*Contributed equally}
\maketitle

\begin{abstract}
Proper scoring rules encourage probabilistic predictions that match the true underlying distribution and are central to model evaluation, with increasing relevance in automated workflows such as AutoML.
In survival analysis, however, their behavior under censoring is not fully understood.
We study commonly used squared and logarithmic scoring rules for right-censored survival data under independent censoring, introducing a notion of marginal properness based on observable outcomes.
Within this framework, we show that the Survival Brier Score (SBS), evaluated at a fixed time point, along with its integrated version (ISBS) and the Right-Censored Log-Loss (RCLL) are strictly proper when all individuals eventually experience the event, but can become improper under finite follow-up or in the presence of cure fractions.
For the SBS, we derive a closed-form expression that reveals the true mechanism: residual mass, corresponding to individuals who remain event-free at study end, systematically biases the score toward underestimating survival, with the effect increasing at later evaluation times and under heavier censoring.

Through simulation experiments, we examine how these issues manifest in finite samples and under misspecification.
The SBS exhibits pronounced improperness at late evaluation times and poor discrimination between models.
The ISBS is more robust due to temporal integration but remains sensitive to tail regularity violations, exhibiting detectable improperness and reduced discriminatory power.
The RCLL behaves consistently with strict properness and effectively separates misspecified models.
Although its absolute scale is affected by both the censoring rate and the prediction grid used for density approximation, relative model rankings remain stable across grid resolutions.
Our findings suggest that reliable evaluation requires adequate test set sizes ($n \ge 50-100$); for ISBS, truncating integration at, e.g., the 80th percentile of follow-up; and for RCLL, using a dense, consistent prediction grid while emphasizing relative model rankings rather than absolute score values in multi-dataset benchmarks.
Under these conditions, prior ISBS‑based benchmarks remain largely trustworthy.

Overall, our results demonstrate how theoretical improperness can translate into misleading model comparisons, underscoring the need for further methodological development in survival model evaluation under censoring and realistic data conditions.
\end{abstract}

\section{Introduction}
\label{sec:intro}

In survival analysis, model evaluation has traditionally focused on concordance indices, which assess the ranking ability of predictions \cite{Harrell1984, Antolini2005, Gonen2005, Uno2007, Gerds2013, Blanche2019, Lillelund2025StopModels}.
While widely used, these measures solely evaluate discrimination, i.e., a model’s ability to distinguish between high- and low-risk individuals.
However, survival models can also be used to make predictions about the distribution of an individual's risk over time, necessitating evaluation metrics beyond concordance indices.
Scoring rules provide an alternative by assessing the overall predictive ability of a model, incorporating both calibration and discrimination \cite{Murphy1973, Gneiting2007}. 
Over the past decades, these measures have gained popularity as probabilistic forecasts are recognized to be superior to deterministic predictions in quantifying uncertainty \cite{Stigler1975, Dawid1984, Dawid1986}.
Formalisation and development of scoring rules have primarily been due to Dawid \cite{Dawid1984, Dawid1986, Dawid2014} and Gneiting and Raftery \cite{Gneiting2007}; although the earliest measures promoting ``rational'' and ``honest'' decision-making date back to the 1950s \cite{Brier1950, Good1952}.
In classification and probabilistic regression settings there are established definitions for scoring rules \cite{Gressmann2018}, most popular of which are the Brier score (or `squared loss')\cite{Brier1950} and the Log loss (or `logarithmic loss')\cite{Good1952}.
Scoring rules have also been extensively studied in Bayesian decision theory \cite{Bernardo2008BayesianTheory}, where they serve as tools for evaluating and comparing probabilistic predictions.
In the context of survival analysis, early discussions on loss functions for censored data can be found in Korn \etal \cite{Korn1990MeasuresData}.
A major development was introduced by Graf \etal \cite{Graf1999}, who proposed the use of inverse probability of censoring weighting (IPCW) to adapt proper scoring rules—most notably the Brier score—to right-censored survival data.
Since then, IPCW-based losses have become a standard approach for extending classical scoring rules to the survival setting and form the basis of several widely used loss functions.

In this paper, we investigate the properness of a subset of scoring rules proposed in survival analysis.
Properness is examined both \textit{theoretically} and \textit{empirically}; the latter assesses whether losses exhibit practical improperness in finite-sample simulation studies and evaluates how sensitively they differentiate models with varying degrees of misspecification.
Motivated by their prominence in classification and regression, we focus on squared and logarithmic losses, which constitute two of the most widely used classes of evaluation measures for survival models, as evidenced by recent benchmarking studies \cite{Jaeger2024, Wissel2023, Li2024, Bommert2022BenchmarkData, Herrmann2021}.
For comprehensive reviews incorporating absolute losses and other variations of scoring rules, including $R^2$ forms, we point readers to the earlier studies of Hielscher \etal \cite{Hielscher2010} and Bøvelstad \etal \cite{Bvelstad2011}, as well as to the subsequent surveys by Choodari-Oskooei \etal \cite{Choodari2012a, Choodari2012b} and Rahman \etal \cite{Rahman2017}.
We limit our scope to measures used for external validation of survival models and therefore exclude losses that require additional hyperparameter estimation, as is common in the optimization of deep learning models \cite{Lee2018, Ren2019, Tjandra2021, han2021, Yanagisawa2023}.
\\\\
In \Cref{sec:notation} we introduce mathematical notation and terminology used throughout; \Cref{sec:properdef} revisits the definition of properness for losses in the classification setting and extends it to the survival setting, emphasizing marginal definitions over observed data; in \Cref{sec:losses} we review several survival scoring losses and analyze their properness under both idealized and practically relevant settings;
in \Cref{sec:experiments} we examine properness in simulation experiments; in \Cref{sec:discussion} we provide a summary of our work and propose guidelines for future practitioners.
All technical proofs are given in Appendices \ref{sec:appendixA} and \ref{sec:appendixB}.

\section{Notation}\label{sec:notation}

Throughout this paper, we assume a single-event survival setting with continuous time and the presence of right-censoring.
All time-to-event variables are non-negative and take values in $\PReals := [0,\infty)$, with $\mathcal{T} \subseteq \PReals$ a set of evaluation times used in certain loss definitions.
Let $Y$ and $C$ be such random variables representing the latent survival time (time to the event of interest) and the latent censoring time, respectively.
We define the observable time $T := \min(Y, C)$ and the event indicator $\Delta := \II(Y \le C)$, where $\II(\cdot)$ denotes the indicator function; 
$\Delta = 1$ corresponds to an observed event, while $\Delta = 0$ indicates right-censoring.
Lower-case letters denote observed realizations of random variables; in particular $(t,\delta)$ denotes an observed survival outcome (time and censoring status) corresponding to $(T,\Delta)$.

For a non-negative continuous random variable $Z$, we denote its probability density function by $f_Z$, its cumulative distribution function by $F_Z(z) = \int_0^z f_Z(u)\,du$ and its survival function by $S_Z(z) = 1 - F_Z(z)$.
All densities are defined almost everywhere with respect to Lebesgue measure on $[0,\infty]$.
Subscripts indicate the random variable inducing the corresponding distribution; for example $S_Y$ and $S_C$ denote the survival functions of the distributions of $Y$ (survival time) and $C$ (censoring time), respectively (analogously for $F_Z$ and $f_Z$).

We use $\hat S$ to denote an estimated survival function, typically obtained from a fitted survival model.
This may represent an unconditional estimate $\hat S(t)$ or a covariate-dependent estimate $\hat S(t \mid X)$, where $X \in \Reals^p$ denotes a vector of covariates.
Correspondingly, $\hat F$ and $\hat f$ denote the estimated cumulative distribution and density functions.
When appearing within the same expression, these are assumed to satisfy
\[
\hat S(t) = 1 - \hat F(t) = 1 - \int_0^t \hat f(u)\,du.
\]

For convenience, we also introduce a latent predicted survival variable $\hat Y$ with survival function $S_{\hat Y}(t) = \hat S(t)$.
Thus $\hat Y$ follows the distribution implied by the prediction model.

Finally, let $\ExtReals := \Reals \cup \{-\infty, +\infty\}$ denote the extended real line.
Throughout the paper, expectations are denoted by $\EE[\cdot]$, with the underlying random variables indicated using underset notation when needed (e.g.\ $\underset{(T,\Delta)}{\EE}[\cdot]$ or $\underset{Y}{\EE}[\cdot]$).

\section{Properness definitions}\label{sec:properdef}

Before turning to survival outcomes, we briefly revisit the notion of properness in probabilistic classification (with probabilistic regression being analogous), which provides the conceptual foundation for proper scoring rules \cite{Gneiting2007, Dawid2014}.
This setting also clarifies why \textit{strict} properness is the relevant notion when evaluating predictive models.

\begin{Def}[Proper Scoring Rule - Binary Classification]
Let $\calP$ be a family of distributions over $[0,1]$, let $\calY = \{0,1\}$, and let $p_Y \in \calP$ denote the true distribution of the binary outcome $Y \in \calY$.
Here, $Y$ is the random variable representing a binary classification outcome of interest, $p_Y$ is the true probability of $Y = 1$, and $p \in \calP$ is a candidate predictive probability provided by a model.
A loss function $L_C : \calP \times \calY \rightarrow \ExtReals$ assigns a numerical score to each probabilistic prediction $p$ given an observed outcome $Y$.  
The loss $L_C$ is called

\begin{enumerate}
    \item \textbf{proper} if 
    $\underset{Y \sim p_Y}{\EE}[L_C(p_Y, Y)] \leq \underset{Y \sim p_Y}{\EE}[L_C(p, Y)] \quad \text{for all } p \in \calP$;
    
    \item \textbf{strictly proper} if the above inequality holds with equality if and only if $p = p_Y$.
\end{enumerate}
\end{Def}

A loss is proper if the expected loss is minimized at the true data-generating distribution.
However, properness alone is insufficient for meaningful model comparison.
As a trivial counterexample, the constant loss $L_{42}(p, Y) := 42$, is technically proper, since all predictions achieve the same (and hence minimal) expected loss, including the true distribution. Nevertheless, such a loss provides no incentive for truthful or accurate probabilistic predictions and cannot discriminate between competing models.
This illustrates why strict properness is essential in practice.
Strictly proper scoring rules uniquely reward the true distribution, thereby encouraging honest probabilistic forecasts and enabling principled model comparison \cite{Gneiting2007, Good1952}.
For this reason, strict properness is typically regarded as a defining property of scoring rules in both theory and applications.
\\\\
Early discussions of loss functions for censored survival data can be found in Korn \etal \cite{Korn1990MeasuresData} and Graf \etal \cite{Graf1999}. 
While these works introduced losses that are now widely used, they did not formally define the notion of a scoring rule or properness in the survival setting.
Subsequent contributions similarly proposed survival scoring rules without explicitly formalizing properness, often relying on analogies to probabilistic regression and classification (e.g., Avati \etal \cite{Avati2020}).
To our knowledge, the first explicit definitions of properness for survival scoring rules appear in Sonabend \cite{Sonabend2021}, Rindt \etal \cite{Rindt2022}, and Yanagisawa \cite{Yanagisawa2023}.
These works differ primarily in the probability distribution with respect to which the expected loss is defined.

Properness can be defined with expectations taken either over the full data-generating process $(Y,C)$ \citep{Yanagisawa2023, Rindt2022} or directly over the observed outcomes $(T,\Delta)$ \cite{Sonabend2021, Hothorn2005}.
Because $(T,\Delta)$ is a measurable function of $(Y,C)$, these formulations are mathematically equivalent in the sense that expectations over one can be expressed as expectations over the other.
However, the conceptual interpretation differs.
Defining properness in terms of $(T,\Delta)$ more closely aligns with the original motivation of proper scoring rules as articulated by Gneiting and Raftery \cite{Gneiting2007}, namely to reward predictions that are close to the actually observed outcome.
In survival analysis, $(T,\Delta)$ constitute the fully observed data, while the true event time $Y$ and the censoring time $C$ are only partially observed: the data reveal which time occurred first and at what time, but the other quantity remains unknown.
Formulating properness directly in terms of $(T,\Delta)$ therefore emphasizes that evaluation must be based solely on observable quantities.

A separate distinction concerns whether properness is defined \textit{marginally} or \textit{conditionally} on covariates $X \in \Reals^p$.
Marginal properness defines the expected loss with respect to the joint distribution of outcomes (e.g., $Y$ in the binary classification setting, $(Y,C)$ or $(T,\Delta)$ in the survival setting) without conditioning on $X$.
Conditional properness instead requires that, for each covariate vector $x$, the true conditional survival distribution minimizes the expected loss under the conditional distribution $(Y,C)\mid X=x$, as for example in Rindt \etal \citep{Rindt2022}.
This distinction is closely tied to assumptions about the censoring mechanism.
In particular, many works assume conditional independence of survival and censoring times given covariates, $Y \indep C \mid X$ \citep{han2021, Yanagisawa2023, Rindt2022, Gerds2006, Kvamme2023, Prince2025}.
Under this assumption, defining properness conditionally is natural: the target of prediction is the conditional survival distribution $P(Y\mid X=x)$, and independence ensures that censoring does not introduce bias after conditioning on $X$.

Another widely used censoring assumption in survival analysis is that of (unconditional) independent censoring, i.e. $Y \indep C$, also known as \textit{random} censoring \citep{Klein2003SurvivalData}.
It is important to note that neither censoring assumption implies the other in general.
In particular, even if $Y \indep C \mid X$, the variables $Y$ and $C$ may still be dependent marginally through their shared dependence on $X$, i.e. $Y \not\!\perp\!\!\!\perp C$.
For example, suppose a binary covariate $X\in\{0,1\}$ determines risk level, and conditional on $X$ (i.e., within the stratum $X = 0$ or $X = 1$), the survival time $Y$ and censoring time $C$ are independent but both tend to be smaller in high-risk individuals.
Then, although $Y$ and $C$ are independent within each risk stratum, pooling over risk groups induces association: individuals with shorter event times also tend to have shorter censoring times because both are driven by the same underlying covariate.

Marginal properness, which operates under independent censoring $Y \indep C$, offers a conceptually simpler framework for studying scoring rules: it avoids modeling the censoring distribution and provides a transparent baseline for theoretical analysis.
However, it is less flexible, as it ignores covariate-dependent censoring and is therefore misspecified for data-generating processes where censoring depends on observed covariates.
Conditional properness, on the other hand, accommodates such settings by requiring that the expected loss be minimized conditional on covariates $X$.
When the conditional censoring distribution can be consistently estimated, this formulation can mitigate bias induced by informative censoring and yield more faithful model comparisons \citep{Gerds2006, Prince2025}.
The cost, however, is stronger modeling assumptions and additional estimation uncertainty introduced through the censoring model.
\\\\
In this work, we focus on marginal properness under the assumption of independent censoring, $Y \indep C$, defined through expectations taken over the population distribution of the observed outcomes $(T,\Delta)$.
Importantly, this concerns the data-generating mechanism and does not restrict predictions to be unconditional; the predicted survival distribution $\hat S$ may still depend on covariates through $\hat S(\cdot \mid X)$.
This framework permits a simpler characterization of properness and isolates the properties of the scoring rules themselves from additional modeling assumptions about the censoring mechanism.
We view conditional notions of properness as complementary and important, and refer the reader to Rindt \etal \cite{Rindt2022} and Yanagisawa \etal \cite{Yanagisawa2023} for further details in that setting.
With this foundation, we now state the censoring assumption and formalize marginal properness for survival scoring rules.

\begin{assumption}[Independent Censoring]
\label{assum:indep_cens}
We assume that the survival and censoring times are independent, i.e. $Y \indep C$.
\end{assumption}

\begin{Def}[Marginally Proper Scoring Rule]
\label{def:surv_proper}
Let $\calP$ be a family of distributions over $\PReals$, let $\mathcal{T} \subseteq \PReals$, and let $S_Y \in \calP$ denote the true survival distribution.
For any $S_Y, \hat S \in \calP$ and for $Y \sim S_Y$ and a censoring variable $C$ taking values in $\PReals$ (independent of $Y$), define $T := \min\{Y,C\}$ and $\Delta := \II(Y \le C)$.
Let $\hat Y \sim \hat S$ denote the predicted survival time.

A survival loss $L: \calP \times \mathcal{T} \times \{0, 1\}\mapsto \ExtReals$ is called

\begin{itemize}
    \item \textbf{proper} if 
    $\underset{(T,\Delta)}{\EE}[L(S_Y, T, \Delta)] \leq \underset{(T,\Delta)}{\EE}[L(\hat S, T, \Delta)] \quad \text{for all } \hat S \in \calP$;
    
    \item \textbf{strictly proper} if the above inequality holds with equality if and only if $\hat S = S_Y$ (equivalently $\hat Y \stackrel{d}{=} Y$).
\end{itemize}
\end{Def}

To clarify the notation, $\mathcal{P}$ denotes a class of time-to-event distributions, $S_Y \in \mathcal{P}$ represents the true survival distribution generating the event time $Y$, and $\hat S$, with corresponding predicted variable $\hat Y \sim \hat S$, will typically denote a model-based approximation to $S_Y$.
Such predictions are often covariate-dependent, e.g., $x \mapsto \hat S(\cdot \mid X = x)$ from an accelerated failure time model, a Cox proportional hazards model, or a machine learning method.
Throughout this work, the prediction $\hat S$ may therefore depend on covariates; for notational simplicity, however, this dependence is suppressed in Definition~\ref{def:surv_proper} and elsewhere.

To accommodate both theoretical analysis and practical applications, we study properness under two settings: one in which the following regularity conditions are satisfied, and one in which they are violated.

\begin{condition}[Regularity Conditions]
\label{cond:regularity}
Let $Y$ denote the event time, $C$ the censoring time, and $\hat S \in \calP$ a candidate survival distribution.
The regularity conditions are:
\begin{enumerate}
\item \textbf{Absolute continuity:} The survival functions $S_Y$, $S_C$, and $\hat S$ are absolutely continuous in $[0,\infty]$, with densities $f_Y$, $f_C$, and $\hat f$, respectively (defined a.e.).
\item \textbf{Boundary behavior:} $S_Y(0)=S_C(0)=\hat S(0)=1$ and $\lim_{t\to\infty} S_Y(t)=\lim_{t\to\infty} S_C(t)=\lim_{t\to\infty} \hat S(t)=0$.
\end{enumerate}
\end{condition}

The second condition is of particular practical interest because real-world survival studies almost always have finite follow-up.
When events are not fully observed by the end of the study $t_{\max}$, the true survival function retains positive probability mass beyond the observation window, i.e., $S_Y(t_{\max}) > 0$.
This situation commonly arises under administrative censoring, where follow-up is truncated at a fixed calendar date \cite{Klein2003SurvivalData}.
Since survival models are fitted to such data, their predicted survival functions will likewise typically retain positive tail mass beyond $t_{\max}$.
Consequently, the regularity conditions are frequently violated in practice, motivating an explicit analysis of both settings.

Finally, we state the following proposition which will be used repeatedly throughout the remainder of the paper when analyzing the properness of survival losses.

\begin{Prop}[\textbf{Expectation of Losses under Independent Censoring}]
\label{prop:exp}
Let $\phi$ be a real-valued, measurable and integrable function, such as an evaluation loss, and let $Y,C,T,S_Y,f_Y,\hat S,\hat f$ be as defined in Section \ref{sec:notation} with $Y \indep C$. Then,

\begin{equation}
\label{eq:expprop}
\underset{(T,\Delta)}{\EE}[\phi(\hatS, T, \Delta)] = \int^\infty_0 \left[ f_Y(t)\,S_C(t)\,\phi(\hatS, t\mid\delta=1) + f_C(t)\,S_Y(t)\,\phi(\hatS, t\mid\delta=0) \right] dt \\
\end{equation}
\end{Prop}

\begin{proof}
See Appendix \ref{sec:appendixA}.
\end{proof}

As noted in Remark~\ref{rmk:cov_dep_exp}, Proposition~\ref{prop:exp} applies equally to covariate-dependent predictions $\hat S(\cdot\mid X)$ without requiring conditional independence between $Y$ and $C$.
Lastly, if there is no censoring in the data generating process then $S_C(t) = 1, f_C(t) = 0 \ \forall t$ and the above reduces to $\int^\infty_0 f_Y(t)\,\phi(\hatS, t) \ dt= \underset{Y}{\EE}[\phi(\hatS, t)]$, recovering the uncensored (regression) setting \cite{Hothorn2005}.

\section{Properness of Squared and Logarithmic Survival Losses}\label{sec:losses}

In this section we focus on squared and logarithmic scoring rules, two broad classes of losses in survival analysis.
We use the terms \textit{scoring rule} and \textit{loss} interchangeably and adopt the convention that losses are \textit{negatively oriented}, i.e. smaller values indicate better predictive performance.
This aligns with standard practice for proper scoring rules, where the negative log-likelihood, for example, is minimized when predictions match the true distribution.
We examine several commonly used survival losses within each class.
For each loss, we review existing claims regarding properness in the literature and formulate propositions that either establish or refute these claims under the framework introduced in \Cref{sec:properdef}.
Our discussion distinguishes between \textit{exact} (\Cref{sec:exact}) and \textit{IPCW} losses (\Cref{sec:ipcw}), where the latter explicitly incorporate the censoring distribution in the loss calculation.

Throughout, losses are presented at the \emph{observation level}.
In practice, model evaluation is performed by averaging these losses over a sample, yielding the empirical risk
\[
\frac{1}{n}\sum_{i=1}^n L(\hat S_i, t_i, \delta_i).
\]

where $\hat S_i := \hat S(\cdot \mid X_i)$ denotes the predicted survival distribution for the $i$-th observation.
We note that properness, is not a property of a single observed loss value but of the \emph{population risk} induced by the loss function, i.e. its expectation over the observable outcomes under the data-generating process.
The empirical risk above is a uniform convex combination of observation-level losses, with weights $w_i = 1/n$ satisfying $\sum_{i=1}^n w_i = 1$.
By \Cref{prop:proper-sum}, the risk induced by such an average is proper (respectively, strictly proper) if and only if the risk induced by the underlying observation-level loss is proper (respectively, strictly proper).
Accordingly, it suffices to analyze observation-level loss functions, with properness understood through the expected loss they induce.

A related distinction is between a loss function, the population risk it defines, and the estimation of that risk in finite samples.
Throughout this section, properness is studied as a property of the population risk under the assumptions specified in \Cref{sec:properdef}.
Questions concerning the estimation of that risk, including the estimation of nuisance quantities such as the censoring distribution (as required by IPCW-based losses), are conceptually separate and do not affect the definition of properness itself.
When discussing IPCW losses, we therefore distinguish between the theoretical loss formulated using the true censoring distribution and practical estimators that replace it with an estimate.

\subsection{Exact survival losses}\label{sec:exact}

We begin by reviewing a class of \emph{exact} survival losses, by which we mean scoring rules that do not explicitly incorporate the censoring distribution into the loss function.
Censoring may still enter implicitly through the event indicator $\delta$, but no inverse probability weighting or explicit modeling of $S_C$ is required.

We consider the following losses:

\begin{itemize}
\item Survival Continuous Ranked Probability Score (SCRPS) \cite{Avati2020}
\begin{equation}
\label{eq:SCRPS}
L_{SCRPS}(\hatS, t, \delta) = \int^{t}_0 \hatF^2(\tau) \ d\tau + \delta \int^\infty_t \hatS^2(\tau) \ d\tau
\end{equation}
\item Negative Log-Likelihood (NLL)
\begin{equation}
\label{eq:NLL}
L_{NLL}(\hatS, t) = - \log\hat f(t)
\end{equation}
\item Right-Censored Log-Loss (RCLL) \cite{Avati2020}
\begin{equation}
\label{eq:RCLL}
L_{RCLL}(\hatS, t, \delta) = - \left(\delta \log\hat f(t) + (1 - \delta)\log\hat S(t)\right)
\end{equation}
\end{itemize}

Avati \etal \cite{Avati2020} introduced the SCRPS as a survival analogue of the continuous ranked probability score and argued for its properness by interpreting it as a weighted CRPS.
However, no formal proof is provided, and the argument does not explicitly account for censoring.
Rindt \etal \cite{Rindt2022} subsequently demonstrated that this weighting argument is insufficient to guarantee properness and provided both theoretical arguments and simulation evidence indicating that SCRPS is improper.

The negative log-likelihood (NLL) is a standard strictly proper scoring rule in regression and density estimation.
Its use in survival analysis has remained limited, for example to augment maximum likelihood objectives in neural network optimization, as in Goldstein et al. \cite{Goldstein2020}.
This is because the formulation in \Cref{eq:NLL} is only defined for fully observed outcomes and implicitly treats all observations as events.
As such, it does not constitute a valid survival scoring rule under censoring, and, to the best of our knowledge, no formal properness result has been established in the survival setting.

Finally, Avati \etal \cite{Avati2020} proposed the right-censored log-loss (RCLL), a special case of the censored likelihood-based scoring rules studied by Dawid \etal \cite{Dawid2014}
Rindt \etal \cite{Rindt2022} and Yanagisawa \cite{Yanagisawa2023} established strict properness of the RCLL under their respective conditional properness definitions, assuming $Y \indep C \mid X$.
Yanagisawa \cite{Yanagisawa2023} also noted that this result assumes a sufficiently dense time grid, recommending at least 16–32 intervals in practice.

Note that in the absence of censoring, all three losses reduce to their uncensored counterparts: the SCRPS reduces to the continuous ranked probability score (CRPS), and both the NLL and RCLL reduce to the logarithmic loss.
In the uncensored world, all three scoring rules are strictly proper \cite{Gneiting2007}.
\\\\
Under marginal properness (\Cref{def:surv_proper}), a survival loss is proper if and only if its expected risk is uniquely minimized by the true survival distribution $\hat S = S_Y$.
Using Proposition~\ref{prop:exp}, we construct explicit counterexamples showing that both the SCRPS (\Cref{eq:SCRPS}) and the NLL (\Cref{eq:NLL}) are not marginally proper under independent censoring: there exist misspecified predictions $\hat S \neq S_Y$ that attain a strictly lower expected loss than the true survival distribution (Appendices~\ref{proof:scrps_not_proper} and \ref{proof:nll_not_proper}).
For the right-censored log-loss we obtain the following result:

\begin{Prop}[\textbf{RCLL Strict Properness}]
\label{prop:rcll_strict_proper}
Assume that the event time $Y$ and censoring time $C$ are independent and continuous on $[0,\infty)$, with survival functions $S_Y$ and $S_C$, respectively.
Let $\hat S$ be a candidate survival prediction, and suppose the regularity conditions hold (Condition~\ref{cond:regularity}).
We also require that the associated densities $f_Y$, $f_C$ and $\hat f$ are strictly positive for all finite $t \ge 0$.
Then the right-censored log-loss (RCLL, \Cref{eq:RCLL})
is strictly proper under the marginal properness definition (\Cref{def:surv_proper}).
\end{Prop}

\begin{proof}
See Appendix~\ref{proof:rcll-proof}.
\end{proof}

This result complements earlier findings by Rindt \etal \cite{Rindt2022} and Yanagisawa \cite{Yanagisawa2023}, who establish strict properness of the RCLL under conditional notions of properness.
Our contribution shows that the RCLL remains strictly proper under a marginal definition based solely on the observable outcomes $(T,\Delta)$.

The proof highlights that strict properness relies on several nontrivial assumptions, which are important to interpret both theoretically and in practice.
First, the RCLL requires all involved densities—the true event density $f_Y$, the censoring density $f_C$, and the predicted event density $\hat f$—to be strictly positive on $[0,\infty)$.
This requirement is essential for the proof and, in particular, for the RCLL to be well defined, as zero density values lead to undefined terms through $\log(0)$.
As discussed by Rindt \etal \cite{Rindt2022} in the context of the \textit{intractable likelihood} problem, this poses practical challenges: many survival models do not directly estimate densities and may thus assign zero density to parts of the time axis (either implicitly or by construction), necessitating careful smoothing or interpolation of the predicted survival or hazard functions.

Second, strict properness depends on regularity conditions of the survival distributions, including the predicted survival function.
When these conditions are violated, the expected RCLL can be potentially minimized by a misspecified prediction $\hat S$ that has a smaller residual mass than the true survival distribution $S_Y$.
This observation is formalized in the following proposition.

\begin{Prop}[\textbf{RCLL Improperness under Regularity Violations}]
\label{prop:rcll_improper}
Assume the setting of Proposition~\ref{prop:rcll_strict_proper}, except that the boundary behavior specified in Condition~\ref{cond:regularity} allows non-vanishing tails $\varepsilon_Y := \lim_{t\to\infty} S_Y(t) > 0$ and $\varepsilon_C := \lim_{t\to\infty} S_C(t) > 0$, with $\varepsilon := \varepsilon_Y \varepsilon_C > 0$.
Let $\hat S$ be a predicted survival function with residual mass $\varepsilon_{\hat Y} := \lim_{t\to\infty} \hat S(t)$.
If $\varepsilon_{\hat Y} < \varepsilon_Y$, then the RCLL scoring rule is improper: the true model no longer minimizes the expected RCLL.
\end{Prop}

\begin{proof}
See Appendix~\ref{proof:rcll-proof}.
\end{proof}

Taken together, these results show that among the exact survival losses considered, only the RCLL preserves marginal properness under independent censoring when the regularity assumptions hold (i.e., at least one of the event time or censoring distributions has a vanishing tail).
When these regularity conditions are violated—specifically, when both $Y$ and $C$ have positive residual mass ($\varepsilon_Y > 0$ and $\varepsilon_C > 0$)—and the predicted survival function has a strictly smaller residual mass than the true model ($\varepsilon_{\hat Y} < \varepsilon_Y$), the RCLL becomes improper.
The extent to which such violations of the regularity assumptions affect the practical relevance of the RCLL for real-world survival data will be investigated in \Cref{sec:experiments}.

\subsection{IPCW survival losses}
\label{sec:ipcw}
\vspace{1em}

A second major class of survival losses consists of scoring rules based on \emph{inverse probability of censoring weighting} (IPCW) \citep{Robins1992}.
The central idea of IPCW is to reweight each observation's loss contribution according to the inverse probability of remaining uncensored, thereby replacing the unavailable complete-data loss with an observable, unbiased analogue \citep{vanDerLaan2003}.
Most IPCW-based survival losses trace back to the seminal work of Graf \etal \cite{Graf1999}, who proposed censoring-adjusted versions of squared and logarithmic error measures for survival prediction.

A key feature of IPCW losses is that they depend explicitly on the censoring distribution.
Throughout this section we formulate these losses in terms of the true censoring survival function $S_C$, rather than an estimator (typically denoted $\hat G$ in the literature).
This allows us to study the intrinsic properness
properties of the loss functions themselves, as defined at the population level.
Issues related to estimating $S_C$—including consistency assumptions ($\hat G \to S_C$ as $n \to \infty$) and the choice between marginal and covariate-dependent censoring models \citep{Gerds2006}—are practically important but pertain to risk estimation rather than population-level properness; we discuss these issues in \Cref{sec:discussion}.
For a comprehensive treatment of IPCW estimation in semiparametric models, we refer to van der Laan and Robins \citep{vanDerLaan2003}.

Let $\tau^* \in \mathcal{T}$ denote a fixed positive evaluation time (or upper integration limit), chosen such that $S_C(\tau^*) > 0$.
Let $(t,\delta)$ denote the observed survival outcome (time and event status), and let all remaining notation follow Section~\ref{sec:notation}.
We consider the following IPCW survival losses:

\begin{itemize}
\item Survival Brier Score (SBS) \cite{Graf1999}
\begin{equation} \label{eq:SBS}
L_{SBS}(\hatS, t, \delta, \tau^*) = \frac{\hatS^2(\tau^*)\,\II(t \leq \tau^*, \delta=1)}{S_C(t)} + \frac{\hatF^2(\tau^*)\,\II(t > \tau^*)}{S_C(\tau^*)}
\end{equation}
\item Integrated Survival Brier Score (ISBS) \cite{Graf1999}
\begin{equation} \label{eq:ISBS}
L_{ISBS}(\hat S, t, \delta, \tau^*) = \int^{\tau^*}_0 L_{SBS}(\hat S, t, \delta, \tau) \ d\tau = \int^{\tau^*}_0 \Big[ \frac{\hat S^2(\tau)\,\II(t \leq \tau, \delta=1)}{S_C(t)} + \frac{\hat F^2(\tau)\, \II(t > \tau)}{S_C(\tau)} \ \Big]d\tau
\end{equation}
\item Integrated Binomial Log-Loss (IBLL) \cite{Graf1999}
\begin{equation}
L_{IBLL}(\hat S, t, \delta, \tau^*) = -\int^{\tau^*}_0 \Big[ \frac{\log\hat F(\tau)\,\II(t \leq \tau, \delta=1)}{S_C(t)} + \frac{\log\hat S(\tau)\,\II(t > \tau)}{S_C(\tau)} \ \Big] d\tau
\end{equation}
\item Integrated Survival Absolute Score (ISAS) \cite{Schmid2011}
\begin{equation}
L_{ISAS}(\hat S, t, \delta, \tau^*) = \int^{\tau^*}_0 \Big[ \frac{\hat S(\tau)\,\II(t \leq \tau, \delta=1)}{S_C(t)} + \frac{\hat F(\tau)\,\II(t > \tau)}{S_C(\tau)} \ \Big] d\tau
\end{equation}
\end{itemize}

All scoring rules listed above share a common decomposition into two IPCW components.
The first term applies when the event has occurred ($t \le \tau, \delta = 1$) and penalizes predictions with high survival probability $\hat S(\tau)$, encouraging the predicted cumulative incidence $\hat F(\tau)$ to be large around the observed event time.
Conversely, the second term applies when the observation has not yet experienced either the event or censoring ($t > \tau$) and penalizes predictions with low survival probability, thereby rewarding high $\hat S(\tau)$.
Thus, despite using different link functions (squared, absolute, or negated logarithm), all scores above share the same underlying structure of weighting contributions across time.

Graf \etal claimed properness for the \emph{uncensored} integrated Brier Score (IBS; see the formula on p. 2537~\cite{Graf1999}), which does not involve IPCW, although neither a formal definition of properness nor a proof was provided.
This result is well established in probabilistic regression, where the IBS reduces to an integrated squared-error score and is known to be strictly proper \citep{Gneiting2007, Gressmann2018}.
For the right-censored case, Rindt \etal \cite{Rindt2022} showed that both the SBS and the ISBS are proper only if both $Y \indep C$ and $C \indep X$.
When censoring depends on covariates, these scoring rules are generally improper, and the authors demonstrate by simulation that a misspecified survival distribution can achieve a lower Brier score than the truth.
We also note the work of Kvamme \etal \cite{Kvamme2023}, who proposed an ISBS variant for \emph{pure administrative censoring}.
Their approach exploits known administrative censoring times to avoid estimating the censoring distribution, but assumes that administrative censoring is the only censoring mechanism, making it inapplicable to general right-censoring settings.
Although the authors derive the expected score under these assumptions and show its approximate agreement with the expected uncensored and IPCW Brier scores when censoring is independent of covariates, they do not formulate or prove a properness result.

The IBLL was first introduced by Graf \etal \cite{Graf1999}, without discussion of properness.
Han \etal\ \cite{han2021} subsequently claimed that the IBLL is a proper scoring rule under conditional expectations given covariates and independent censoring.
Under the same conditional framework, Rindt \etal\ \cite{Rindt2022} established properness of the IBLL under covariate-independent censoring ($C \indep X$), by analogy with their proof for the ISBS, and demonstrated via simulation that the IBLL becomes improper when this assumption is violated.
Despite its conceptual connection to the binomial log-likelihood, the IBLL has seen far less adoption than its squared-error counterpart, the ISBS—both in methodological work \cite{han2021, Rindt2022, Kvamme2019, Gandy2025, Burk2026AData} and in software—with implementations are currently restricted to \texttt{mlr3proba} \cite{pkgmlr3proba} and \texttt{pycox} \cite{pkgpycox}.
Given its limited practical use and properness properties that closely parallel those of the ISBS, we do not consider the IBLL further.

The ISAS was introduced by Schmid \etal \cite{Schmid2011}, as an IPCW-based alternative to the ISBS, building on the earlier absolute-distance measure between true and predicted survival functions \cite{Schemper2000}.
Unlike squared-error–based measures, ISAS offers a more interpretable, outlier-robust assessment of prediction error, reflecting the mean absolute deviation between survival distributions.
Under independent censoring, ISAS serves as a consistent estimator of the mean absolute error; however, neither Schmid \etal\ nor subsequent works provide formal proofs or claims regarding its properness as a scoring rule.
Since our focus is on squared and logarithmic losses, we include the ISAS here primarily for completeness.
\\\\
We now turn to the Survival Brier Score (SBS) and its integrated version (ISBS), which are by far the most widely used IPCW losses in survival analysis.
We begin with the SBS, which evaluates predictive accuracy at a fixed evaluation time $\tau^*>0$.
Our first result establishes that, under standard regularity conditions, the SBS is strictly proper.

\begin{Prop}[\textbf{SBS Strict Properness}]
\label{prop:sbs_proper}
Assume that the event time $Y$ and censoring time $C$ are independent, with survival functions $S_Y$ and $S_C$, respectively.
Let $\tau^* \in \mathcal T$ be a fixed evaluation time with $S_C(\tau^*)>0$, and let $\hat S$ denote a candidate survival prediction.
Also assume that the standard regularity conditions hold (Condition~\ref{cond:regularity}).
Then, under the marginal properness definition (\Cref{def:surv_proper}), the Survival Brier Score (SBS, \Cref{eq:SBS}) is strictly proper at $\tau^*$.
In particular,
\begin{equation}
\label{eq:sbs_expectation_simple}
\underset{(T,\Delta)}{\EE}\big[L_{SBS}(\hat S,T,\Delta)\big]
= \hat S^2(\tau^*)\,F_Y(\tau^*) + \hat F^2(\tau^*)\,S_Y(\tau^*),
\end{equation}
which is strictly convex as a function of $\hat S(\tau^*)$ and uniquely minimized at $\hat S(\tau^*) = S_Y(\tau^*)$.
\end{Prop}

\begin{proof}
See Appendix~\ref{proof:sbs-proof}.
\end{proof}

As a scoring rule, the SBS is well defined under minimal assumptions: unlike the RCLL, it does not involve logarithmic terms and therefore does not require strictly positive densities, relying only on the positivity assumption $S_C(\tau^*) > 0$.
Properness, however, largely depends on the regularity of the underlying survival distributions.
When these conditions are violated—most notably when survival distributions retain residual probability mass beyond the effective observation horizon—the SBS is no longer strictly proper.
In this case, the expected loss remains convex but is minimized at a misspecified survival probability.
This behavior is characterized by the following proposition.

\begin{Prop}[\textbf{SBS Improperness under Regularity Violations}]
\label{prop:sbs_improper}
Assume the setting of Proposition~\ref{prop:sbs_proper}, except that the boundary behavior specified in Condition~\ref{cond:regularity} allows non-vanishing tails $\varepsilon_Y := \lim_{t\to\infty} S_Y(t) > 0$ and $\varepsilon_C := \lim_{t\to\infty} S_C(t) > 0$, with $\varepsilon := \varepsilon_Y \varepsilon_C > 0$.
Then the expected SBS becomes

\begin{equation}
\label{eq:sbs_expectation_eps}
\underset{(T,\Delta)}{\EE}\big[L_{SBS}(\hat S,T,\Delta)\big]
= \hat S^2(\tau^*)F_Y(\tau^*) + \hat F^2(\tau^*)\Big(S_Y(\tau^*) - \frac{\varepsilon}{S_C(\tau^*)}\Big)
\end{equation}
The above expression is strictly convex in $\hat S(\tau^*)$ and is uniquely minimized at
\begin{equation}
\label{eq:sbs-minimizer-prop}
x^\star = \frac{S_Y(\tau^*) - \varepsilon / S_C(\tau^*)}{1 - \varepsilon / S_C(\tau^*)},
\end{equation}
which satisfies $x^\star < S_Y(\tau^*)$ whenever $S_Y(\tau^*) \in (0,1)$ and $S_C(\tau^*) > \varepsilon$.
Consequently, the SBS is improper under marginal properness in this setting, as the expected loss is minimized by a systematically underestimated survival probability at the evaluation time $\tau^*$.
\end{Prop}

\begin{proof}
See Appendix~\ref{proof:sbs-proof}.
\end{proof}

To make this improperness result explicit, we express the deviation between the unique minimizer of the expected SBS and the true survival probability at $\tau^*$ as

\begin{equation}
\label{eq:sbs-shift}
x^\star - S_Y(\tau^*)
= -\frac{\varepsilon\bigl(1 - S_Y(\tau^*)\bigr)}{S_C(\tau^*) - \varepsilon}
\le 0 .
\end{equation}

This term can be interpreted as a systematic downward \textit{bias} induced by residual tail mass.

If $S_Y(\tau^*) < 1$ and $S_C(\tau^*) > \varepsilon$—that is, if the evaluation time lies within the effective support of the data (i.e., before censoring survival approaches its tail limit)—then for any $\varepsilon>0$ the bias is strictly negative.
Also, since the bias is monotone in $\varepsilon$ (the derivative with respect to $\varepsilon$ is strictly negative), larger residual tail mass leads to greater underestimation.
As such, the expected SBS is minimized by a survival prediction that underestimates the true survival probability $S_Y(\tau^*)$.
On the other hand, if at least one of the survival distributions $S_Y$ or $S_C$ has no residual mass at infinity (i.e. $\varepsilon = 0$), the bias disappears and the minimizer reduces to $x^\star = S_Y(\tau^*)$, thereby restoring strict properness as stated in Proposition~\ref{prop:sbs_proper}.

\Cref{eq:sbs-shift} also reveals how the choice of evaluation time affects the magnitude of the improperness.
For small residual mass $0 < \varepsilon \ll S_C(\tau^*)$, a first-order approximation yields (Remark~\ref{rem:sbs_bias_approx}):

\begin{equation}
\label{eq:sbs-shift-approx}
x^\star - S_Y(\tau^*) \;\approx\; -\; \varepsilon \;\cdot\; \frac{1 - S_Y(\tau^*)}{S_C(\tau^*)}.
\end{equation}

Thus, for fixed  $\varepsilon$, the magnitude of the bias is determined by the ratio term in \Cref{eq:sbs-shift-approx}.
At early evaluation times, $S_Y(\tau^*) \approx 1$, so the numerator is close to zero and the bias is negligible, largely independent of the censoring distribution.
As $\tau^*$ increases, the cumulative event probability $F_Y(\tau^*) = 1 - S_Y(\tau^*)$ grows while the censoring survival typically decreases.
Both effects increase the ratio and thereby induce larger bias.
In particular, heavy censoring, i.e. small $S_C(\tau^*)$, acts multiplicatively, making the improperness more pronounced even at intermediate time points.
Consequently, even small residual tail mass may lead to practically relevant bias when the SBS is evaluated late in follow-up or under substantial censoring.
This behavior is examined empirically in Section~\ref{sec:robustness}.
\\\\
While the SBS is well suited to assess model predictive accuracy at a specific task-relevant time point of interest, in practice, it is computed across a grid of time points to form a \emph{prediction error curve} \cite{Mogensen2012}.
While informative, such curves do not yield a single summary score for model comparison, motivating the use of integrated scores such as the ISBS, which we analyze next.

\begin{Prop}[\textbf{ISBS Properness Characterization}]
\label{prop:isbs_proper}
Assume the setting of Proposition~\ref{prop:sbs_proper}, but with the boundary conditions relaxed, by allowing for possibly non-vanishing tails $\varepsilon_Y := \lim_{t\to\infty} S_Y(t) \ge 0$ and $\varepsilon_C := \lim_{t\to\infty} S_C(t) \ge 0$, with $\varepsilon := \varepsilon_Y \varepsilon_C \ge 0$.
Then, the expected ISBS (\Cref{eq:ISBS}) equals the time integral of pointwise expected SBS losses:

\begin{equation}
\label{eq:isbs-expectation}
\begin{aligned}
\underset{(T,\Delta)}{\EE}\Big[L_{ISBS}(\hat S, T, \Delta)\Big] 
&= \int_0^{\tau^*} \underset{(T,\Delta)}{\EE}\Big[L_{SBS}(\hat S, T, \Delta)\Big]\,d\tau\\
&= \int_0^{\tau^*} \Big[\hat S^2(\tau)(1-S_Y(\tau)) + \hat F^2(\tau)\Big(S_Y(\tau)-\frac{\varepsilon}{S_C(\tau)}\Big)\Big]\,d\tau.
\end{aligned}
\end{equation}

If $\varepsilon = 0$, the integrand is uniquely minimized at $\hat S(\tau)=S_Y(\tau)$ for all $\tau \in [0,\tau^*]$, and the ISBS is strictly proper.
When both tails persist ($\varepsilon>0$), then by Proposition~\ref{prop:sbs_improper} the pointwise expected SBS at each time $\tau$ is strictly convex and uniquely minimized at a value strictly below the corresponding true survival probability $S_Y(\tau)$, rendering the ISBS improper.
\end{Prop}

\begin{proof}
See Appendix~\ref{proof:isbs-proof}.
\end{proof}

In summary, the SBS and ISBS are strictly proper under independent censoring and the vanishing-tail conditions of Propositions~\ref{prop:sbs_proper}--\ref{prop:isbs_proper}.
When both the event-time and censoring distributions retain positive mass beyond the observation window, as may occur under administrative censoring or when the event-time distribution exhibits persistent tail mass (e.g., due to cure fractions or insufficient follow-up) \cite{Nie2011InferenceCensoring,Srivastava2021ImpactStudies}, both losses are minimized by systematically underestimated survival predictions.

Critically, this improperness is intrinsic to the Brier score itself: it arises at the population-risk level under the true data-generating mechanism and persists even if the censoring distribution $S_C$ is known exactly.
Thus, the pathology is not merely an estimation issue, as has previously been suggested \citep{Rindt2022}, but a fundamental property of the score under non-vanishing tails.
Moreover, the quantity $\varepsilon$ governing the degree of improperness depends on features of the latent event and censoring distributions that are not recoverable from the observed outcomes $(T,\Delta)$, reflecting the broader non-identifiability of latent failure times under censoring \citep{Tsiatis1975}.
Consequently, practitioners cannot generally verify whether the conditions required for strict properness hold.
Even under independent censoring and with a perfectly known censoring distribution, these widely used IPCW Brier scores may therefore be improper in realistic settings, with potentially unknown implications for model evaluation and comparison.
The magnitude and practical consequences of this improperness are investigated empirically in Section~\ref{sec:experiments}.

\section{Experiments}\label{sec:experiments}

In this section, we examine whether survival scoring rules satisfy marginal properness \emph{in practice}, as defined in \Cref{def:surv_proper}.
Although the previous sections establish theoretical guarantees and characterize conditions under which properness holds or fails, applied model evaluation is conducted with finite samples, estimated censoring distributions, and potentially misspecified models.
It is therefore essential to investigate whether finite-sample effects and violations of the regularity conditions identified in the theoretical analysis lead to practically relevant differences in empirical model rankings, potentially affecting the conclusions of previously reported benchmarking studies.

We begin by introducing a simulation framework for assessing empirical properness and detecting violations (\Cref{sec:sim_design}).
Within this framework, we evaluate the right-censored log-likelihood (RCLL; \Cref{eq:RCLL}), the Survival Brier Score (SBS; \Cref{eq:SBS}), and the Integrated Survival Brier Score (ISBS; \Cref{eq:ISBS}).
Using this setup, we first study properness using a very large sample size to approximate the population risk in \Cref{def:surv_proper} via the law of large numbers, thereby examining alignment between theoretical results and empirical behavior (\Cref{sec:large-n-properness}).
We then analyze finite-sample robustness, where limited follow-up and non-negligible residual tail mass can induce empirical violations despite asymptotic guarantees (\Cref{sec:robustness}).
Lastly, we investigate sensitivity to model misspecification by quantifying how strongly each scoring rule penalizes increasingly misspecified predictions relative to the true survival distribution, and assess the impact of sample size and density approximation to the empirical behavior of the RCLL (\Cref{sec:sensitivity}).

\subsection{Simulation Design}\label{sec:sim_design}

To assess the empirical properness of survival scoring rules, we generate synthetic right-censored survival data from parametric Weibull distributions.
The Weibull family provides sufficient flexibility to induce a wide range of hazard shapes and tail behaviors, making it suitable for testing properness across heterogeneous settings.
We perform a total of $K = 10,000$ independent simulations, indexed by $k \in \{1, \dots, K\}$.
To simplify notation, we suppress the simulation index $k$ in the distributional quantities and sampled variables whenever no ambiguity arises.
For each simulation $k$, we draw $m = 1,000$ distinct Weibull \emph{distribution triplets} $(S_Y^d,S_C^d,{\hat S}^d)$, indexed by $d$.
Each triplet consists of the true survival distribution $S_Y$, the censoring distribution $S_C$ and a candidate prediction $\hat S$.
As the draws are independent, the candidate prediction is almost always misspecified relative to the true data-generating process.
From each triplet, we sample $n$ i.i.d. observations:

\begin{equation}
\label{eq:draw_times}
\begin{aligned}
y_i^d &\overset{\text{i.i.d.}}{\sim} \text{Weibull}(\alpha_Y^d, \sigma_Y^d), \quad 
c_i^d \overset{\text{i.i.d.}}{\sim} \text{Weibull}(\alpha_C^d, \sigma_C^d),\\    
t_i^d &= \min\{y_i^d, c_i^d\}, \quad \delta_i^d = \mathbb{I}(y_i^d \le c_i^d), \quad i = 1,...,n.
\end{aligned}
\end{equation}

Here, $y_i^d$ denotes the true event time, $c_i^d$ the censoring time, and $(t_i^d,\delta_i^d)$ the observed outcome (times and censoring indicator).
The maximum observed time is $t_{\max}^d=\max\limits_{1 \le i \le n} t_i^d$.
As event and censoring times are generated independently, $Y \indep C$ holds by construction, satisfying Assumption~\ref{assum:indep_cens}.
The sample size $n$ represents the number of test-set observations for which predictions are evaluated, with index $i$ enumerating individual samples within a given distribution triplet.
The candidate prediction $\hat S^d$ is defined parametrically and used directly in scoring.

The Weibull distributions are parameterized by the shape and scale parameters ($\alpha>0, \sigma>0$), independently drawn for event, censoring, and candidate distributions:

\begin{equation}
\label{eq:weibull_params}
(\alpha_Y^d,\sigma_Y^d),(\alpha_C^d,\sigma_C^d),(\hat \alpha^d,\hat \sigma^d) \overset{i.i.d.}{\sim} \text{Uniform}(0.5, 5)
\end{equation}

\begin{figure}[h]
\centering
\includegraphics[width=0.8\textwidth]{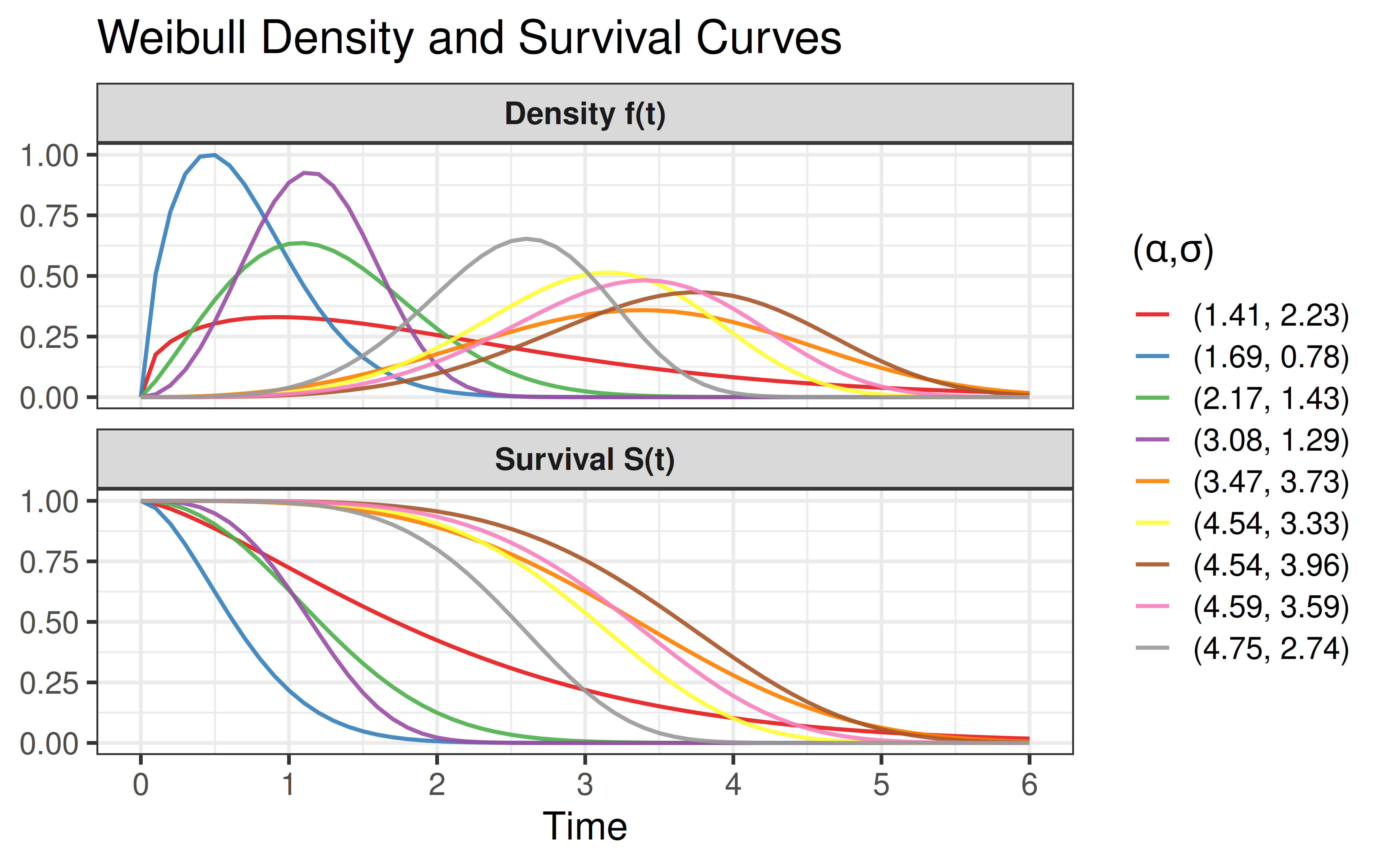}
\caption{Weibull density and survival curves for randomly generated shape ($\alpha$) and scale ($\sigma$) parameters.
Each curve corresponds to a different parameter pair, illustrating the broad variability across simulated distributions.}
\label{fig:weibull}
\end{figure}

This wide parameter range induces high variability across event and censoring distributions (\Cref{fig:weibull}), thereby increasing the likelihood of detecting violations of properness, while avoiding extreme values (e.g., near-zero or very large shape parameters) that could cause numerical instability.
Also, since Weibull survival functions satisfy $S(t) = \exp(-(t/\sigma)^{\alpha}) \rightarrow 0$ as $t \rightarrow \infty$, the survival mass vanishes in the tail.
In our setup, as $n$ increases, the empirical maximum $t_{\max}^d$ grows, and Proposition~\ref{prop:s_tmax} implies $S^d(t_{\max}^d)\rightarrow0$ in expectation.

For a scoring rule $L$, we compute the empirical loss difference between the true and predicted survival functions:

\begin{equation}
\label{eq:delta_j}
D_{L}^{k,d} = \frac{1}{n} \sum_{i=1}^{n} \left[ L(S_Y^d, t_i^d, \delta_i^d) - L(\hat{S}^d, t_i^d, \delta_i^d) \right]
\end{equation}

which averages the loss difference across the $n$ i.i.d.\ observations generated from distribution triplet $d$ in simulation $k$.

We then aggregate over the $m$ independently sampled distribution triplets to obtain the simulation-level mean loss difference and its empirical standard deviation:

\begin{equation}
\label{eq:mean_score}
\bar{D}_{L}^k = \frac{1}{m} \sum_{d=1}^{m} D_{L}^{k,d} \quad,\quad
s^k_L = \sqrt{ \frac{1}{m - 1} \sum_{d=1}^{m} \left( D^{k,d}_{L} - \bar{D}_{L}^{k} \right)^2}
\end{equation}

The first averaging step (over $n$) yields the empirical loss difference under a fixed data-generating mechanism $(S_Y^d, S_C^d,\hat S^d)$.
The second averaging step (over $m$) aggregates these quantities across a diverse collection of randomly sampled Weibull triplets.
For a given sample size $n$, $\bar{D}_{L}^k$ therefore serves as a Monte Carlo estimate of the empirical loss difference under our simulated design.
For a proper scoring rule $L$, we expect $\bar{D}_{L}^{k} \le 0$; values $\bar{D}_{L}^{k} > 0$ indicate empirical violations, meaning that, on average, a misspecified prediction $\hat S$ attains a lower expected loss than the true survival distribution $S_Y$.

To assess whether the observed mean difference $\bar{D}_{L}^k$ indicates a \textbf{statistically significant violation of properness}—accounting for Monte Carlo variability across simulations—we construct a 95\% confidence interval (CI) based on the $t$-distribution, since the population variance of $D^{k,d}_{L}$ is unknown:

\begin{equation}
\label{eq:ci}
\bar{D}_{L}^k \pm t_{0.975, m - 1} \cdot \frac{s^k_L}{\sqrt{m}}
\end{equation}

The simulation $k$ is flagged as a violation if $\bar{D}_L^k > \epsilon$ (with $\epsilon>0$) and the CI lower bound exceeds zero.
We adopt a conservative threshold $\epsilon = 10^{-4}$ to focus on practically meaningful violations.
This procedure allows us to assess whether each scoring rule reliably distinguishes models with genuinely different predictive performance.
Minor score differences between similarly performing models are often unavoidable but unlikely to affect practical decisions such as model selection or tuning.
In practice, such small inconsistencies are typically mitigated via standard resampling procedures (e.g., repeated K-fold cross-validation).

Overall, this design evaluates empirical properness across a large and heterogeneous collection of survival distributions and sample sizes.
A pseudo-code summary is provided in Algorithm~\ref{alg:properness}.

\begin{algorithm}[ht]
\caption{Simulation procedure to assess empirical marginal properness of survival scoring rules}
\label{alg:properness}
\begin{algorithmic}[1]
\State \textbf{Input:} Survival loss $L$, simulations $K$, distributions per simulation $m$, sample size $n$, positive threshold $\epsilon$

\For{$k = 1$ to $K$}
    \For{$d = 1$ to $m$}
    
        \State Sample Weibull parameters 
        $(\alpha_Y, \sigma_Y)$, 
        $(\alpha_C, \sigma_C)$, 
        $(\hat{\alpha}, \hat{\sigma}) \sim \text{Uniform}(0.5,5)$
        
        \State Define survival functions 
        $S_Y^d$, $S_C^d$, and candidate prediction $\hat S^d$
        
        \For{$i = 1$ to $n$}
            \State Draw $y_i^d \sim \text{Weibull}(\alpha_Y, \sigma_Y)$
            \State Draw $c_i^d \sim \text{Weibull}(\alpha_C, \sigma_C)$
            \State $t_i^d = \min(y_i^d, c_i^d)$
            \State $\delta_i^d = \mathbb{I}(y_i^d \le c_i^d)$            
        \EndFor

        \State Compute empirical mean loss under truth:
        $L_Y^{k,d} = \frac{1}{n}\sum_{i=1}^n L(S_Y^d, t_i^d, \delta_i^d)$
        
        \State Compute empirical mean loss under prediction:
        $L_{\hat S}^{k,d} = \frac{1}{n}\sum_{i=1}^n L(\hat S^d, t_i^d, \delta_i^d)$
        
        \State Set $D_L^{k,d} = L_Y^{k,d} - L_{\hat S}^{k,d}$ using \eqref{eq:delta_j}
    \EndFor
    
    \State Compute simulation-level mean and standard deviation $(\bar{D}_L^k,s_L^k)$ using \eqref{eq:mean_score}
    \State Construct 95\% confidence interval $(CI_{\text{lower}}, CI_{\text{upper}})$ via \eqref{eq:ci}
    
    \If{$\bar{D}_L^k > \epsilon$ \textbf{and} $CI_{\text{lower}} > 0$}
        \State Record violation
    \EndIf
    
\EndFor

\State \textbf{Output:} $\{\bar{D}_L^k,s_L^k\}_{k=1}^K$, confidence intervals, and number of violations
\end{algorithmic}
\end{algorithm}

\subsection{Large-Sample Empirical Properness}\label{sec:large-n-properness}

We examine the empirical properness of RCLL, SBS, and ISBS in a large-sample setting under the Weibull simulation design (\Cref{sec:sim_design}; $K=10{,}000$ simulations), using $n = 10{,}000$ observations and $m = 1{,}000$ independently sampled Weibull distribution triplets per simulation.
At this sample size, the empirical average loss within each sampled distribution triplet closely approximates its corresponding population risk.
Aggregating these risks across the $m$ randomly sampled Weibull triplets further approximates the expectation over the distribution family $\calP$ appearing in \Cref{def:surv_proper}.
Consequently, the empirical comparison closely mirrors the population-level marginal properness criterion studied in our theoretical results.
Moreover, the observed follow-up period is sufficiently long that the probabilities of remaining event-free or uncensored by the end of follow-up are negligible.
As a result, this effectively satisfies the vanishing-tail conditions underlying the strict properness results for RCLL, SBS, and ISBS (Propositions~\ref{prop:rcll_strict_proper}, \ref{prop:sbs_proper}, and \ref{prop:isbs_proper}).

The RCLL (\Cref{eq:RCLL}) depends on both the survival function $S(t)$ and the corresponding density $f(t)$.
Since both the data-generating and predicted distributions are Weibull, these quantities are obtained directly from the corresponding parametric distributions.
For the SBS (\Cref{eq:SBS}), evaluation is performed at fixed time points $\tau^*$ corresponding to the 10th ($q_{0.1}$), 50th (median, $q_{0.5}$), and 90th ($q_{0.9}$) percentiles of the observed times.
The ISBS (\Cref{eq:ISBS}) extends the SBS by integrating over time.
We approximate this integral numerically using the trapezoidal rule on a grid of 50 equidistant time points spanning the interval between the 5th and 80th percentiles of the observed outcome times.
To facilitate comparison across simulation settings with different time horizons, the resulting integral is normalized by the length of the corresponding evaluation interval.

Both SBS and ISBS require the censoring survival function $S_C(t)$ for inverse probability of censoring weighting.
We first consider an oracle setting in which the censoring distribution is known and compute $S_C$ directly from the generating Weibull distribution.
To reflect practical applications, we also estimate $S_C$ via the Kaplan–Meier estimator $\hat G(t)$ under independent censoring, applied to $(t_i^d, 1-\delta_i^d)$ within each distribution triplet $d$.
Because SBS and ISBS may require $\hat G(t)$ beyond the observed support, we use constant interpolation (and extrapolation) to ensure the estimator is defined over the full time range, as is standard practice \cite{survival-package}.

The results are summarized in the last rows of \Cref{tab:sbs_violations} and \Cref{tab:sbs_violations_est} for $n=10{,}000$.
Neither SBS nor ISBS exhibits statistically significant violations of properness.
RCLL, which does not depend on the censoring distribution, likewise remains strictly proper.
This behavior persists across all sample sizes considered (\Cref{sec:robustness}); for brevity, the corresponding RCLL results are omitted from the tables.
Overall, both under the idealized setting with known $S_C$ and under Kaplan--Meier estimation of $S_C$ in the simulated independent censoring setting, the empirical findings are consistent with the theoretical strict properness results for RCLL, SBS, and ISBS.

\begin{table}[ht]
\centering
\scriptsize
\caption{Empirical violations of properness for the SBS and ISBS across $K = 10{,}000$ simulations with known survival and censoring Weibull distributions. 
SBS was evaluated at the 10th percentile ($q_{0.1}$), median ($q_{0.5}$), and 90th percentile ($q_{0.9}$) of observed times.
ISBS was computed using 50 equidistant time points between the 5th and 80th percentiles of observed times.
Columns report the number of violations (\#Viol.), the violation rate (Rate), and the mean score difference (true minus predicted) among simulations flagged as violations ($\bar{D}_{L}^{\text{Viol}}$).}
\label{tab:sbs_violations}
\renewcommand{\arraystretch}{1.2}
\begin{tabularx}{\textwidth}{c|CCC|CCC|CCC|CCC}
\toprule
$n$ & \multicolumn{3}{c|}{\textbf{SBS}($q_{0.1}$)} & \multicolumn{3}{c|}{\textbf{SBS}($q_{0.5}$)} & \multicolumn{3}{c|}{\textbf{SBS}($q_{0.9}$)} & \multicolumn{3}{c}{\textbf{ISBS}} \\
 & \#Viol. & Rate & $\bar{D}_{\text{SBS}}^{\text{Viol}}$ 
 & \#Viol. & Rate & $\bar{D}_{\text{SBS}}^{\text{Viol}}$ 
 & \#Viol. & Rate & $\bar{D}_{\text{SBS}}^{\text{Viol}}$ 
 & \#Viol. & Rate & $\bar{D}_{\text{ISBS}}^{\text{Viol}}$ \\
\midrule
10    & 4458 & 0.446  & 0.00562 & 743  & 0.0743 & 0.00721 & 1207 & 0.121  & 0.0147  & 433 & 0.0433 & 0.00231 \\
25    & 476  & 0.0476 & 0.00124 & 406  & 0.0406 & 0.00380 & 637  & 0.0637 & 0.00642 & 121 & 0.0121 & 0.000532 \\
50    & 438  & 0.0438 & 0.000741 & 321 & 0.0321 & 0.00181 & 469 & 0.0469 & 0.00369 & 34 & 0.0034 & 0.000235 \\
100   & 169  & 0.0169 & 0.000397 & 210 & 0.021  & 0.000977 & 318 & 0.0318 & 0.00206 & 3  & 0.0003 & 0.000139 \\
250   & 56   & 0.0056 & 0.000192 & 102 & 0.0102 & 0.000433 & 182 & 0.0182 & 0.000869 & 0  & 0      & --       \\
500   & 15   & 0.0015 & 0.000128 & 73  & 0.0073 & 0.000255 & 129 & 0.0129 & 0.000473 & 0  & 0      & --       \\
1000  & 0    & 0      & --       & 33  & 0.0033 & 0.000156 & 67  & 0.0067 & 0.000282 & 0  & 0      & --       \\
2500  & 0    & 0      & --       & 0   & 0      & --       & 20  & 0.002  & 0.000152 & 0  & 0      & --       \\
5000  & 0    & 0      & --       & 0   & 0      & --       & 2   & 0.0002 & 0.000102 & 0  & 0      & --       \\
10000 & 0    & 0      & --       & 0   & 0      & --       & 0   & 0      & --       & 0  & 0      & --       \\
\bottomrule
\end{tabularx}
\end{table}

\subsection{Robustness at Small Sample Sizes}\label{sec:robustness}

We next investigate empirical properness in finite-sample regimes where the tail regularity conditions required for strict properness may not hold.
Although the Weibull event and censoring distributions used in our simulation design have vanishing tails at infinity ($S_Y(t) \rightarrow 0$ and $S_C(t) \rightarrow 0$ as $t \rightarrow \infty$; \Cref{sec:sim_design}), the realized maximum observation time $t_{\max}$ in a finite sample may not be sufficiently large for the corresponding survival probabilities, $S_Y(t_{\max})$ and $S_C(t_{\max})$, to be negligible.
Consequently, the residual tail product $S_Y(t_{\max})\cdot S_C(t_{\max})$ can remain strictly positive, 
serving as a finite-sample proxy for the theoretical residual mass $\varepsilon$ that drives improperness in the expected SBS and ISBS losses (Propositions~\ref{prop:sbs_improper} and \ref{prop:isbs_proper}).
We therefore expect empirical violations of properness to occur in small-sample regimes, and we investigate whether the observed violations follow the predicted theoretical mechanism.

For each sample size $n \in \{10, 25, 50, 100, 250, 500, 1000, 2500, 5000, 10000\}$, we perform $K = 10{,}000$ simulation replicates, each based on $m = 1{,}000$ independently generated Weibull distribution triplets (\Cref{sec:sim_design}).
All quantities—including $S_Y$, $S_C$, and candidate predictions $\hat S$—are generated under independent censoring (Assumption~\ref{assum:indep_cens}).
Results using the true censoring distribution are reported in \Cref{tab:sbs_violations}; replacing $S_C$ with the Kaplan–Meier estimator $\hat G(t)$ yields largely similar findings (\Cref{tab:sbs_violations_est}), indicating that estimating the censoring distribution has little impact on empirical properness in this setting.

To assess the practical relevance of improperness, we report the average loss difference only among simulations exhibiting statistically significant violations.
Let

\[
\mathcal{V} = \{k : \bar{D}_L^k > \epsilon \text{ and } CI_{\text{lower}} > 0\}
\]

denote the set of simulations flagged as violations (Algorithm~\ref{alg:properness}).
We then compute

\[
\bar{D}_{L}^{\text{Viol}} = \frac{1}{|\mathcal{V}|} \sum_{k \in \mathcal{V}} \bar{D}_L^k,
\]

where $\bar{D}_L^k$ is the simulation-level mean loss difference (\Cref{eq:mean_score}).
This positive conditional average quantifies how strongly a misspecified prediction outperforms the true survival distribution when properness is violated, thereby reflecting the practical relevance of such deviations.

Across all sample sizes, the RCLL exhibited no violations, indicating strong finite-sample robustness; it is therefore omitted from \Cref{tab:sbs_violations}.
In contrast, the SBS and ISBS showed more nuanced behavior.
SBS violations were frequent at small sample sizes, with rates at $n=10$ ranging from 7.4\% at the median ($q_{0.5}$) to 44.6\% at the 10th percentile ($q_{0.1}$) and 12.1\% at the 90th percentile ($q_{0.9}$).
As sample size increased, violations decreased rapidly: for $n \ge 250$, rates fell below 2\% across all evaluation times, converging to zero for large $n$.
For all sample sizes $n \ge 25$, violations were more frequent at later evaluation times, with $\text{SBS}(q_{0.9})$ consistently exceeding $\text{SBS}(q_{0.5})$ and $\text{SBS}(q_{0.1})$.
Despite the relatively high frequency of violations at very small sample sizes, their practical magnitude remained modest, with $\bar{D}_{\text{SBS}}^{\text{Viol}} < 0.015$ throughout.
The ISBS demonstrated stronger robustness: violations occurred only for $n \le 50$, with rates of 4.3\% at $n=10$, 1.2\% at $n=25$, and 0.3\% at $n=50$; for $n > 50$, violations were rare.
Moreover, the magnitude of violations was substantially smaller than for SBS, with $\bar{D}_{\text{ISBS}}^{\text{Viol}} < 0.001$ for $n \ge 25$, suggesting that ISBS converges rapidly to empirical properness and remains reliable at moderate sample sizes.

The observed SBS behavior is explained by the theoretical mechanism driving improperness (\Cref{prop:sbs_improper}).
For each sampled distribution triplet $d$ and sample size $n$, we compute the residual tail product at the maximum observed time $t_{\max}^d$:

\begin{equation}
\label{eq:empirical_epsilon_d}
\varepsilon^d = S_Y^d(t_{\max}^d) \cdot S_C^d(t_{\max}^d).
\end{equation}

This quantity approximates the remaining joint survival mass beyond the observed follow-up.
We then evaluate $S_Y^d$ and $S_C^d$ at the evaluation times $\tau^*$ (10th, 50th, and 90th percentiles), allowing us to compute the implied minimizer $x^\star$ of the expected SBS (\Cref{eq:sbs-minimizer-prop}) and its deviation from $S_Y(\tau^*)$ via \Cref{eq:sbs-shift}.
Averaging over the $m = 1{,}000$ distribution triplets yields simulation-level estimates of $S_Y(\tau^*)$, $S_C(\tau^*)$, the bias $x^\star - S_Y(\tau^*)$, and the residual tail mass:

\begin{figure}[ht]
\centering
\includegraphics[width=1.0\textwidth]{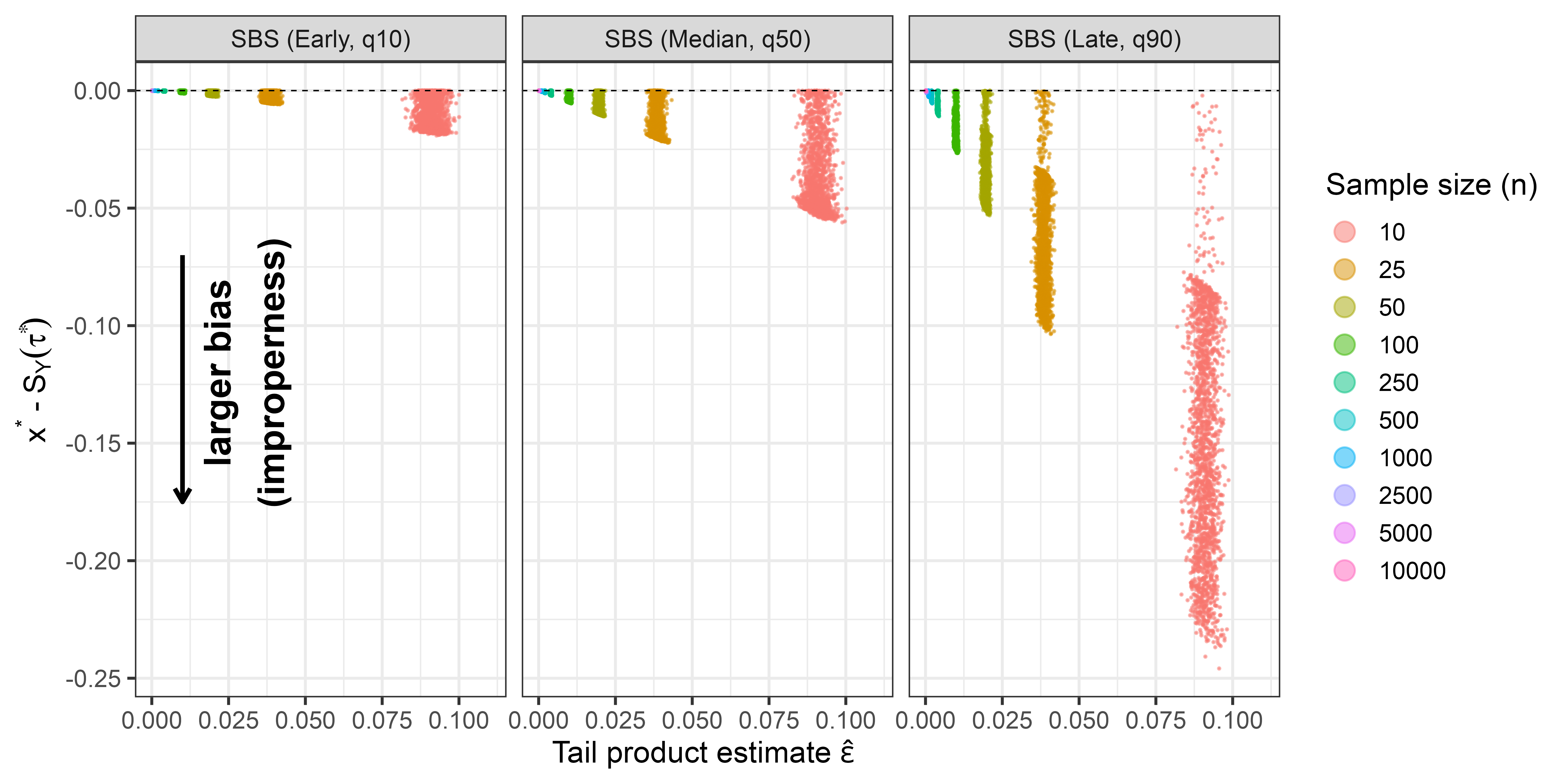}
\caption{
Bias of the SBS minimizer as a function of the residual tail mass.
The figure shows the empirical bias of the minimizer of the expected SBS loss relative to the true survival probability $x^\star - S_Y(\tau^*)$ (\Cref{eq:sbs-shift}) against the estimated residual tail mass $\hat\varepsilon$ (\Cref{eq:empirical_epsilon}), based on $K = 10{,}000$ simulations.
A random 20\% subsample of the simulations is shown for visualization.
Colors denote the sample size $n$.
Panels correspond to early, median, and late evaluation times $\tau^*$.
Negative values indicate underestimation of the true survival probability.}
\label{fig:epsilon_vs_bias}
\end{figure}

\begin{equation}
\label{eq:empirical_epsilon}
\hat\varepsilon=\frac{\sum_{d=1}^{1000} \varepsilon^d}{1000}.
\end{equation}

\Cref{fig:epsilon_vs_bias} confirms the monotonic relationship implied by \Cref{eq:sbs-shift}: larger values of $\hat\varepsilon$ produce increasingly negative values of $x^\star - S_Y(\tau^*)$, indicating greater underestimation of the true survival probability by the minimizer of the expected SBS. 
As sample size increases, the observed maximum $t_{\max}^d$ increases, causing both $S_Y(t_{\max}^d)$ and $S_C(t_{\max}^d)$ to decrease (\Cref{prop:s_tmax}).
Consequently, $\hat\varepsilon$ approaches zero via Equations~(\ref{eq:empirical_epsilon_d})-(\ref{eq:empirical_epsilon}), recovering the regularity conditions required for strict properness.
Accordingly, the minimizer converges to the true survival probability, as reflected by the vanishing bias in \Cref{fig:epsilon_vs_bias}.
For the smallest sample sizes ($n = 10,25,50$), the residual tail mass remains substantial ($\hat\varepsilon\approx 0.1 - 0.02$), leading to noticeable downward bias, particularly at the median and late evaluation times, whereas the bias rapidly diminishes as the sample size increases and becomes negligible for sufficiently large $n$.
These trends closely mirror the violation frequencies reported in \Cref{tab:sbs_violations}: smaller sample sizes are associated with larger residual tail mass, greater bias, and more frequent improperness violations, all of which disappear as the regularity conditions are progressively recovered.
The visible clustering of points into vertical bands simply reflects the discrete set of simulated sample sizes, with each band corresponding to a characteristic range of residual tail masses induced by a fixed $n$; under a finer grid of sample sizes, these bands would merge into a continuous point cloud.

\begin{figure}[ht]
\centering
\includegraphics[width=1.0\textwidth]{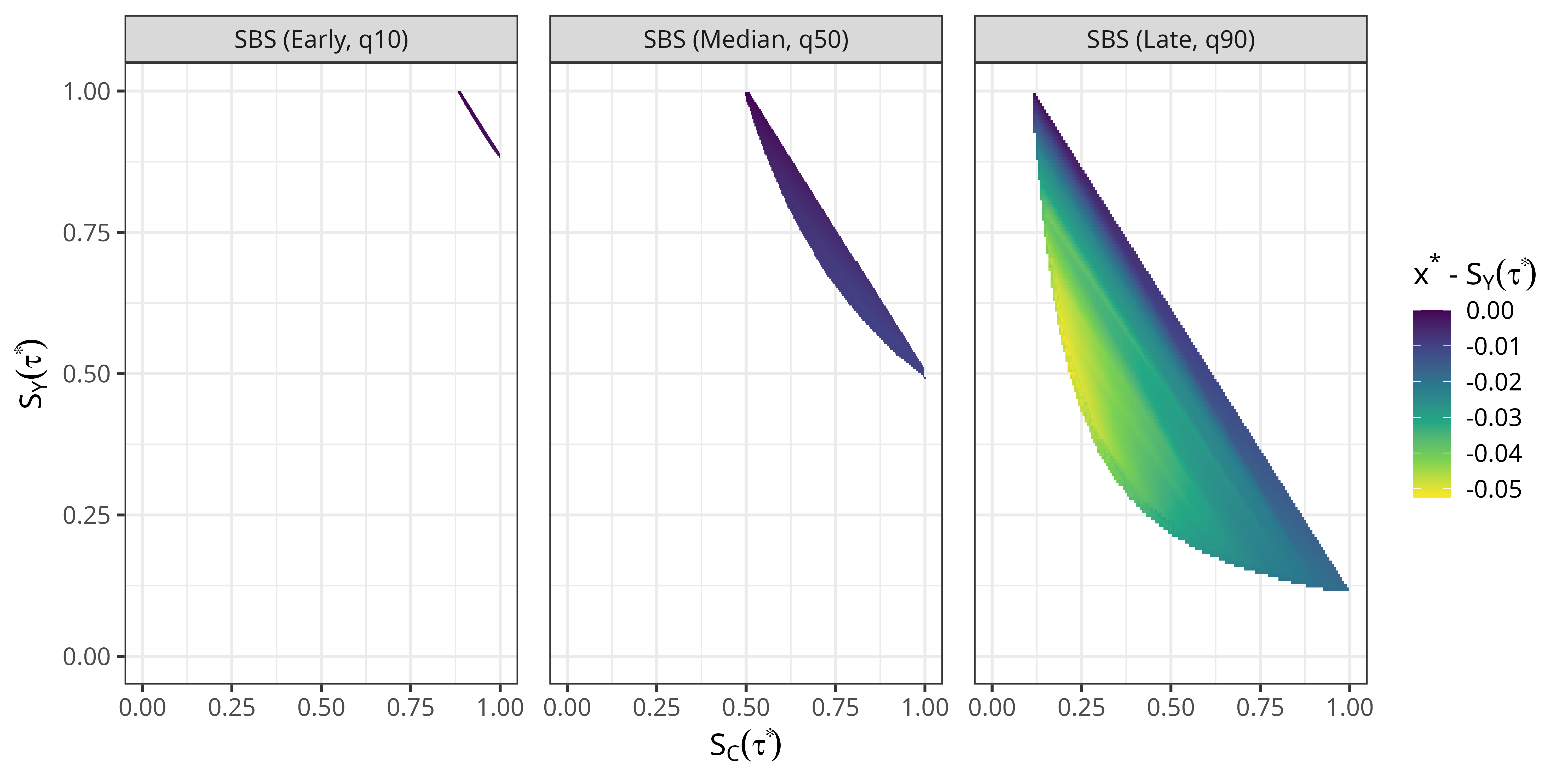}
\caption{
Bias of the SBS minimizer as a function of the event and censoring survival probabilities.
The figure shows the interpolated surface of the empirical difference between the minimizer of the expected SBS loss and the true survival probability $x^\star - S_Y(\tau^*)$ (\Cref{eq:sbs-shift}) against the true event survival probability $S_Y(\tau^*)$ and censoring survival probability $S_C(\tau^*)$, for $n = 50$ ($\hat\varepsilon \approx 0.02$), based on $K = 10{,}000$ simulations.
Panels correspond to early, median, and late evaluation times $\tau^*$.
Negative values indicate underestimation of the true survival probability.}
\label{fig:bias_surface}
\end{figure}

To further understand the role of $S_Y(\tau^*)$ and $S_C(\tau^*)$ in determining the magnitude of the SBS bias, \Cref{fig:bias_surface} visualizes the deviation $x^\star - S_Y(\tau^*)$ as a function of these quantities for $n = 50$ ($\hat\varepsilon \approx 0.02$).
The results closely follow the theoretical first-order approximation (\Cref{eq:sbs-shift-approx}).
For fixed $\varepsilon > 0$, the bias is proportional to the ratio:

\[
\frac{1 - S_Y(\tau^*)}{S_C(\tau^*)}.
\]

Consequently, the bias increases as the cumulative event probability $1 - S_Y(\tau^*)$ increases and the censoring survival probability $S_C(\tau^*)$ decreases.
At early evaluation times ($\text{SBS}(q_{0.1})$), both $S_Y(\tau^*)$ and $S_C(\tau^*)$ remain close to one, keeping the ratio small and the bias negligible.
At median evaluation times ($\text{SBS}(q_{0.5})$), the range of both quantities broadens, but the ratio remains moderate, resulting in only small deviations from properness.
At late evaluation times ($\text{SBS}(q_{0.9})$), lower event survival $S_Y(\tau^*)$ together with lower censoring survival $S_C(\tau^*)$ substantially increases the ratio, producing the pronounced negative bias visible in the yellow region of \Cref{fig:bias_surface}.
Accordingly, the strongest underestimation occurs when the cumulative event probability is high and the censoring survival probability is low.

In our simulation design, evaluation times are chosen as quantiles of the observed time distribution, restricting $S_Y(\tau^*)$ and $S_C(\tau^*)$ to moderate values ($\gtrsim 0.1$).
Consequently, the bias remains bounded even at late evaluation times.
More extreme settings, such as $S_C(\tau^*) \rightarrow 0$, would produce substantially larger bias because the amplification factor $(1 - S_Y(\tau^*)) / S_C(\tau^*)$ diverges.
\\\\
Taken together, the empirical results support the theoretical characterization of SBS improperness under independent censoring.
The residual tail mass $\varepsilon$ determines whether bias can occur, while the interaction between $S_Y(\tau^*)$ and $S_C(\tau^*)$—and hence the choice of evaluation time $\tau^*$—determines the magnitude of the underestimation.
As $\varepsilon \to 0$, both the bias and the frequency of violations vanish, consistent with the asymptotic properness established in \Cref{sec:large-n-properness}.
In our simulation design, this convergence is driven by the increasing support of the observed follow-up as the sample size grows.

This conclusion, however, does not extend to settings with administrative censoring, where the observation horizon remains fixed.
In such settings, $\varepsilon$ can remain strictly positive even as $n \rightarrow \infty$, so SBS may remain asymptotically improper.
To illustrate this, we conducted an additional simulation study using the same data-generating mechanism as in \Cref{sec:sim_design}, but imposing administrative censoring at the 80th percentile of the true event-time distribution.
Even at $n = 10{,}000$, SBS exhibited persistent improperness: violations occurred in 6.35\% of simulations at the late evaluation time, with $\hat\varepsilon \approx 0.06$ and $\bar{D}_{\text{SBS}}^{\text{Viol}} \approx 0.027$.
Additional details are provided in the Supplementary Material (\Cref{sec:sup}).

These findings also help explain the empirical behavior of ISBS.
Because ISBS averages SBS over the evaluation interval, local deviations from properness are attenuated, resulting in substantially greater empirical robustness than SBS.
Nevertheless, the theoretical analysis indicates that contributions from sufficiently late evaluation times—where $S_C(\tau^*)$ usually becomes small—can still amplify improperness.
In practice, restricting integration to a sufficiently early horizon (e.g., up to the 80-90th percentile of observed times) mitigates these issues, consistent with prior recommendations \cite{Bender2021AAnalysis, Stuber2023ACRC, Wissel2025}.

\subsection{Sensitivity to Model Misspecification}
\label{sec:sensitivity}

While properness guarantees that a scoring rule is minimized in expectation by the true survival distribution, it does not characterize how strongly the score penalizes misspecified predictions in practice.
To study this, we compare a correctly specified model to increasingly misspecified alternatives and evaluate empirical score differences between the true survival distribution and misspecified predictions.
This analysis provides a practical assessment of how effectively each scoring rule discriminates between correct and incorrect models.

\subsubsection{Simulation Design}

We generated training data from a non-proportional hazards data-generating process (DGP) in which the true event time distribution $Y$ follows a log-normal accelerated failure time (AFT) model with covariate-dependent location and group-specific scale.
Specifically, omitting the observation index $i$ for notational simplicity, the true event-time model is

\begin{equation}
\label{eq:DGP}
\log(Y)=\eta + \sigma Z, \qquad Z \sim N(0,1),    
\end{equation}

with linear predictor
\[
\eta
=
\beta_0
+ \beta_1 x_{1}
+ \beta_2 x_{2}
+ \beta_{12} x_{1} x_{2},
\]
where the regression coefficients were set to
\[
\beta_0 = 1.15,\quad
\beta_1 = 0.15,\quad
\beta_2 = -0.55,\quad
\beta_{12} = -0.75,
\]
and covariates
\[
x_1 \sim N(0,1),
\qquad
x_2 \sim \mathrm{Bernoulli}(0.5),
\]

where $x_2$ represents a binary treatment indicator (e.g., drug or surgical intervention).
To induce strong non-proportional hazards with crossing survival curves between the two treatment groups, the scale parameter $\sigma$ was allowed to depend on $x_2$ through
\[
\sigma =
\begin{cases}
0.5, & x_2 = 0, \\
1.5, & x_2 = 1.
\end{cases}
\]

We chose this configuration as it is both challenging for misspecified models (e.g., under proportional hazards or restrictive parametric assumptions) and representative of real-world clinical studies where the treatment shows time-varying effects~\cite{Rothwell1999PredictionStudy,Li2015StatisticalMethods,Dormuth2022WhichGuideline}.

Independent right-censoring ($Y \indep C$) was imposed through exponential censoring times, i.e. $C \sim \mathrm{Exp}(\lambda_C)$, together with administrative censoring at $t_{\max} = 10$ to reflect finite study follow-up, as commonly encountered in practice.
Although the true distribution of event times $Y$ depends on covariates $X$, censoring $C$ is generated independently of $Y$, placing this simulation within the marginal properness framework of Section~\ref{sec:properdef}.
The observed survival outcomes were given by
\[
T = \min(Y, C, t_{\max}),
\qquad
\Delta = \II\!\left(Y \leq \min(C,t_{\max})\right).
\]

\begin{figure}[ht]
\centering
\includegraphics[width=1.0\textwidth]{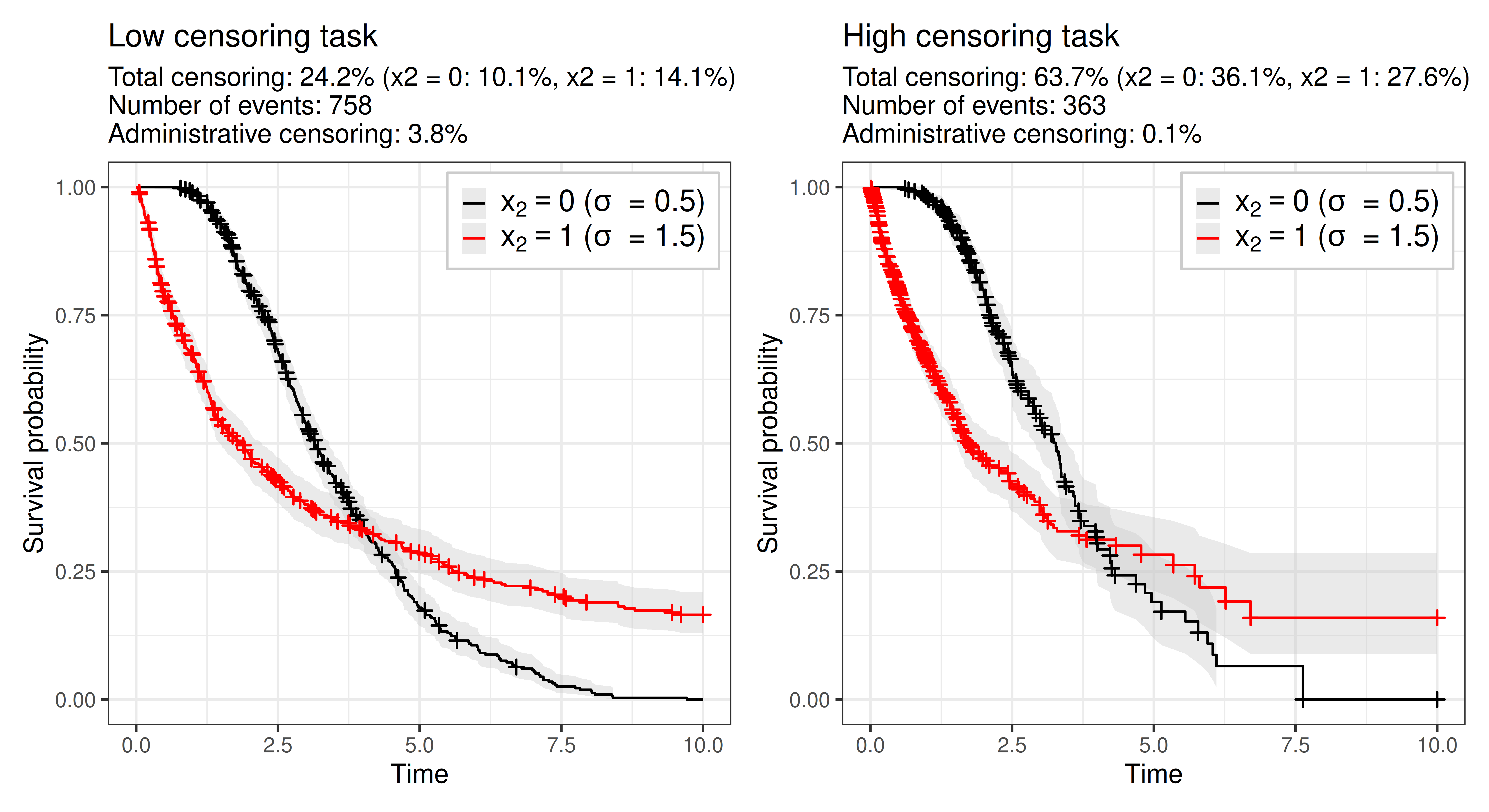}
\caption{Kaplan--Meier curves for the two simulation tasks, stratified by the binary covariate $x_2$, showing the induced crossing survival structure from the non-proportional hazards data-generating process. Left: low censoring task ($\simeq 24$\% total censoring).
Right: high censoring task ($\simeq 64$\% total censoring).
Both settings include administrative censoring at $t_{\max}=10$.}
\label{fig:KM_DGP}
\end{figure}

We consider two training datasets, hereafter referred to as \textit{tasks}, corresponding to low and high censoring, generated by varying the exponential censoring rate while keeping the event-time mechanism fixed.
For the low-censoring task, we used $\lambda_C = 0.075$, yielding approximately 24\% total censoring, while for the high-censoring task we used $\lambda_C = 0.45$, yielding approximately 64\% total censoring.
In both settings, censoring arose from the combination of independent exponential censoring and administrative censoring at $t_{\max}=10$, with some variation across the two $x_2$-groups.
Kaplan--Meier curves for both tasks are shown in \Cref{fig:KM_DGP}, illustrating the induced crossing survival structure and differing censoring regimes.

For each task, models were trained on \(n_{\mathrm{train}} = 1000\) observations, which was sufficient for reliable estimation of the oracle model while preserving meaningful misspecification among competing approaches.
We compared the correctly specified log-normal model (Oracle) to a sequence of increasingly misspecified alternatives, including models with omitted interaction terms, misspecified scale structure, alternative parametric families, as well as semiparametric and nonparametric learners.
An overview of all model configurations is provided in \Cref{tab:sensitivity_models}.
Parametric models were fitted using \texttt{flexsurv}~\citep{Jackson2016Flexsurv:R}, the Cox and Kaplan--Meir models via \texttt{survival}~\citep{survival-package}, and random survival forests via \texttt{ranger}~\citep{Wright2017}.
All models were interfaced through \texttt{mlr3}~\citep{Lang2019} and \texttt{mlr3proba}~\citep{pkgmlr3proba}, with additional components from \texttt{mlr3extralearners}~\citep{mlr3extralearners} and \texttt{mlr3pipelines}~\citep{mlr3pipelines2021} used to ensure consistent model specification (e.g., inclusion of interaction terms).
Further details on the simulation setup are provided in the Supplementary Material (\Cref{sec:sup}).

\begin{table}[ht]
\footnotesize
\centering
\caption{Models considered in the misspecification benchmark. The oracle model matches the data-generating process, while subsequent models introduce progressively different forms of misspecification.}
\label{tab:sensitivity_models}
\renewcommand{\arraystretch}{1.1}
\begin{tabularx}{\textwidth}{>{\raggedright\arraybackslash}l l l c c}
\toprule
Model & ID & Type & Varies with $x_2$$^{*}$ & Interaction \\
\midrule
Log-normal (interaction, $x_2$-scale) & \textbf{Oracle} & Parametric & \checkmark & \checkmark \\
Log-logistic (interaction, $x_2$-shape) &  LLogis & Parametric & \checkmark & \checkmark \\
Weibull (interaction, $x_2$-shape) & Weibull & Parametric & \checkmark & \checkmark \\
Log-normal ($x_2$-scale only) & LogNorm\_scale & Parametric & \checkmark & \ding{55} \\
Log-normal (interaction only) & LogNorm\_int & Parametric & \ding{55} & \checkmark \\
Log-normal (both omitted) & LogNorm & Parametric & \ding{55} & \ding{55} \\
Cox Proportional Hazards (interaction) & Cox\_int & Semi-parametric & -- & \checkmark \\
Random Survival Forest & RSF & Nonparametric & -- & implicit$^{\dagger}$ \\
Kaplan--Meier & KM & Baseline & -- & -- \\
\bottomrule
\end{tabularx}

\vspace{2mm}
$^{*}$ For the log-normal model this refers to the scale parameter $\sigma$; for the Weibull and log-logistic models, it refers to the shape parameter $\alpha$ \citep{Jackson2016Flexsurv:R}.
In all cases, the parameter is allowed to depend on the binary covariate $x_2$.
$^{\dagger}$ Random Survival Forest can capture nonlinear covariate effects and interactions through recursive partitioning and ensemble learning, but do not explicitly model covariate-dependent scale parameters.
\end{table}

We generated 100 Monte Carlo test sets from each data-generating process, each with $n_{\mathrm{test}} = 1000$ observations.
Using a fixed, relatively large test size reduces Monte Carlo variability in the estimated losses and enables reliable comparison of the scoring rules.
For each test set, we evaluated the true survival function $S_Y(t)$ across observations, obtained directly from the DGP (\Cref{eq:DGP}), and the predicted survival functions $\hat S(t)$ from all models, including the Oracle and misspecified alternatives (\Cref{tab:sensitivity_models}).
All survival functions were evaluated on a common prediction grid defined by the unique event times of the corresponding training task (758 and 363 time points for the low- and high-censoring tasks, respectively).
Parametric models were evaluated directly on this grid, random survival forests produce predictions on these time points by construction, and Cox and Kaplan--Meier predictions were restricted to the common grid via piecewise-constant survival interpolation implemented in \texttt{survdistr}~\citep{survdistr2026}, ensuring that all models were evaluated at identical prediction time points.

To quantify the discrepancy between $S_Y(t)$ and $\hat S(t)$ across each test set, we computed the mean integrated absolute error (MIAE) across observations, up to the last observed event time $\tau$, as
\begin{equation}
\label{eq:MIAE}
\mathrm{MIAE}(S_Y, \hat S)
=
\frac{1}{n_{\mathrm{test}}}
\sum_{i=1}^{n_{\mathrm{test}}}
\int_0^{\tau}
\bigl| S_Y^{(i)}(t) - \hat S^{(i)}(t) \bigr| \, dt.
\end{equation}

The integral was approximated using the trapezoidal rule over the common evaluation grid defined by the unique event times of the corresponding training task.
For each test set, we further computed the mean empirical losses $L_{\mathrm{true}}(S_Y)$ and $L_{\mathrm{pred}}(\hat S)$ over the test observations, and report their difference:

\begin{equation}
\label{eq:excess_loss}
\Delta L = L(\hat S) - L(S_Y),
\end{equation}

which represents the excess loss of a predicted survival function relative to the true model.
Negative values indicate empirical properness violations, where a misspecified model attains a lower loss than the true survival distribution.

We considered three scoring rules $L$: the RCLL (\Cref{eq:RCLL}), the Survival Brier Score (SBS; \Cref{eq:SBS}) evaluated at $t=5$, chosen to lie beyond the crossing point (\Cref{fig:KM_DGP}), and the Integrated Survival Brier Score (ISBS; \Cref{eq:ISBS}) computed using 100 equidistant evaluation points up to the 90th percentile of observed times (approximately 6.19 for the low-censoring task and 3.35 for the high-censoring task).
The censoring survival function $S_C(t)$ for IPCW was estimated from the training data using the Kaplan--Meier estimator under independent censoring—with piecewise-constant interpolation where required—consistent with the marginal properness framework adopted throughout this work.
Corresponding Kaplan--Meier estimates are shown in the Supplementary Material (\Cref{sec:sup}) and closely match the true censoring distributions, as expected for the sample size considered.
Since all models provide survival predictions but the RCLL additionally requires the density $f(t)$, we estimated densities uniformly across models via linear interpolation of $S(t)$ over the ordered grid of unique event times from the training data, following Rindt \etal\cite{Rindt2022}.
This approach was also applied to the true model, i.e. estimating $f_Y(t)$ via linear interpolation of $S_Y(t)$, to ensure comparability.
When extrapolation beyond the observed support was required, we extended the final linear segment, preserving monotonicity of the survival function and ensuring non-negative density estimates.
The dense evaluation grid is expected to yield sufficiently accurate density approximations; this assumption is further examined in a subsequent experiment (\Cref{sec:rcll_grid}).

\subsubsection{Model Discrimination under Misspecification}

\begin{figure}[ht]
\centering
\includegraphics[width=1.0\textwidth]{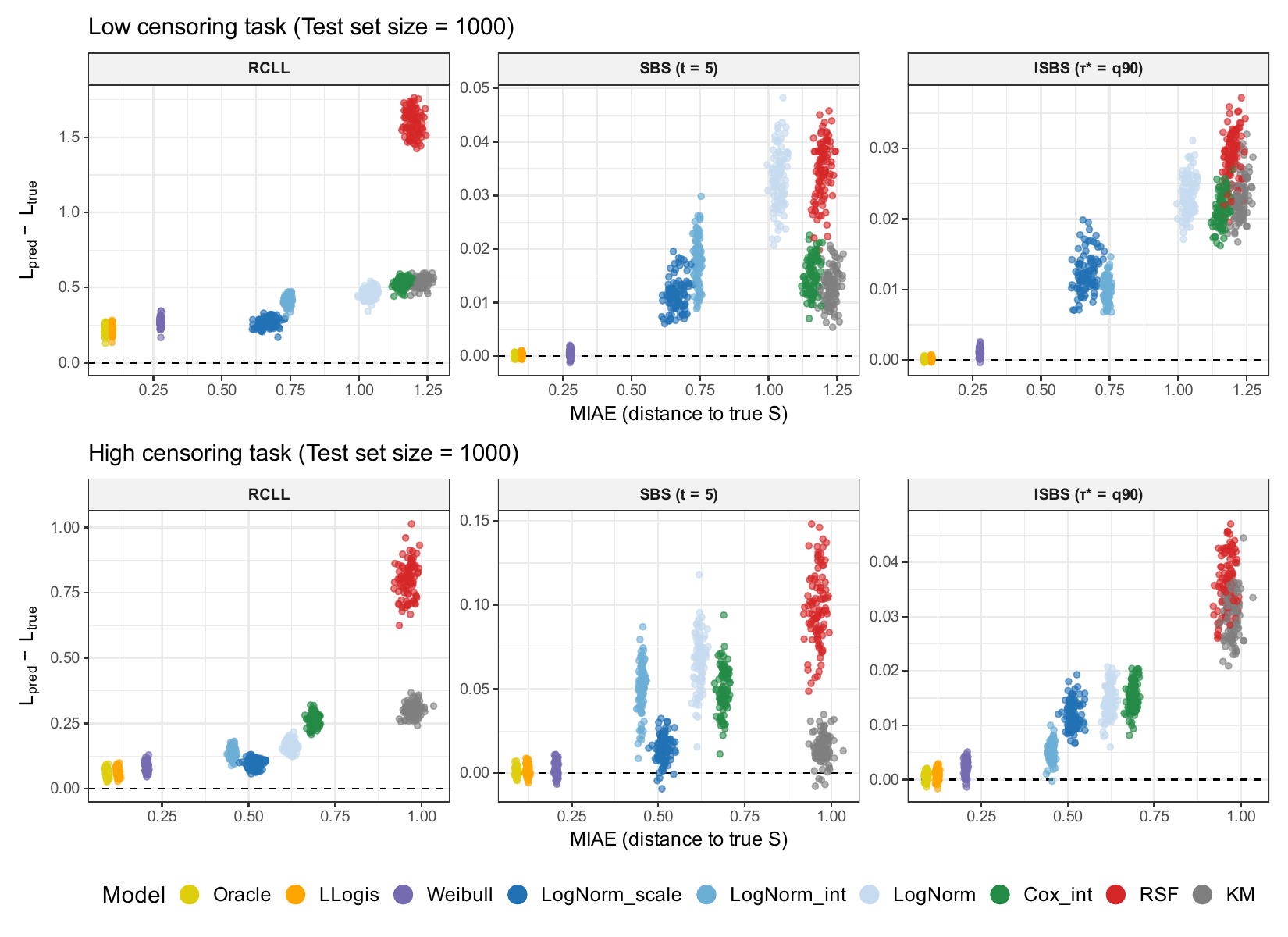}
\caption{Sensitivity to model misspecification under low (top) and high (bottom) censoring.
Each point corresponds to one of 100 Monte Carlo test replicates ($n_{\mathrm{test}} = 1000$ per replicate) and shows the excess loss $\Delta L = L_{\mathrm{pred}}(\hat S) - L_{\mathrm{true}}(S_Y)$ versus the mean integrated absolute error (MIAE) between $\hat S$ and $S_Y$.
Models are arranged from left to right by increasing discrepancy from the true survival function (larger MIAE).
Facets correspond to RCLL, SBS evaluated at $t = 5$, and ISBS integrated up to the 90th percentile of observed times.}
\label{fig:sens_large_n}
\end{figure}

In \Cref{fig:sens_large_n}, we plot the excess loss ($\Delta L$; \Cref{eq:excess_loss}) against the mean integrated absolute error (MIAE; \Cref{eq:MIAE}).
Across both censoring regimes, a consistent pattern emerges.
The SBS exhibits substantial overlap between models, weak alignment between $\Delta L$ and MIAE, and high variability across Monte Carlo replicates—particularly under high censoring—making it unreliable for discriminating between models with different degrees of misspecification.

In contrast, both RCLL and ISBS show increasing excess loss with increasing MIAE, yielding a similar overall ranking of models according to their deviation from the true survival distribution. However, RCLL provides clearer separation between neighboring models and exhibits lower variability across replicates, whereas ISBS operates on a more compressed scale, producing broader overlap between models with similar predictive accuracy.
This difference is particularly apparent for models with comparable MIAE.
For example, in the low-censoring setting, the misspecified log-normal variants (`LogNorm\_scale' and `LogNorm\_int') have similar MIAE values but are well separated by RCLL, whereas their ISBS values overlap substantially.
Likewise, RCLL clearly distinguishes `RSF' from `Cox\_int' and `KM', while ISBS provides considerably less separation.
Models whose survival predictions are genuinely very similar, such as the Oracle and the log-logistic model, remain difficult to distinguish under either score.
An interesting exception occurs in the high-censoring setting, where the MIAE ordering of `LogNorm\_scale' and `LogNorm\_int' reverses.
ISBS follows this ordering in terms of excess loss, whereas RCLL assigns a slightly larger mean excess loss to `LogNorm\_int', indicating that the two scores emphasize different aspects of predictive error under heavier censoring.
Nevertheless, RCLL continues to distinguish neighboring models more consistently than ISBS across the remaining comparisons.

\begin{figure}[ht]
\centering
\includegraphics[width=1.0\textwidth]{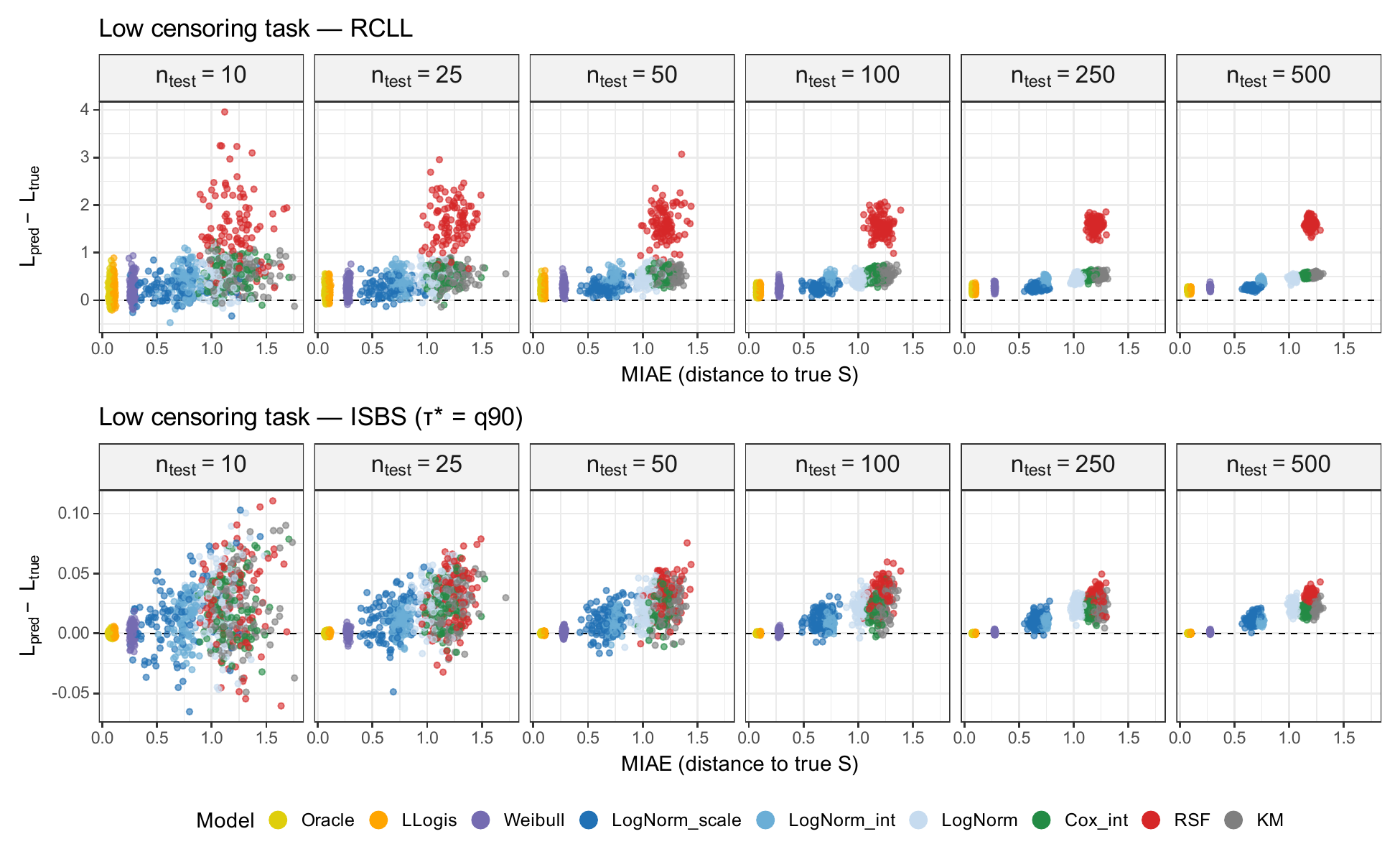}
\caption{Effect of test set size on sensitivity to model misspecification under low censoring for RCLL (top) and ISBS (bottom).
Each panel corresponds to a test set size $n_{\mathrm{test}} \in \{10, 25, 50, 100, 250, 500\}$.
Points represent 100 Monte Carlo replicates and show the excess loss $\Delta L = L_{\mathrm{pred}}(\hat S) - L_{\mathrm{true}}(S_Y)$ versus the mean integrated absolute error (MIAE) between $\hat S$ and $S_Y$.
Models are ordered along the horizontal axis by increasing discrepancy from the true survival function (larger MIAE).
ISBS is integrated up to the 90th percentile of observed times.}
\label{fig:sens_vary_n}
\end{figure}

Importantly, for the large test sets considered ($n_{\mathrm{test}} = 1000$), RCLL yields uniformly positive excess loss across all models, consistent with empirical properness.
In contrast, ISBS exhibits occasional negative excess loss values ($\Delta L < 0$), primarily for models close to the truth, indicating finite-sample violations of properness.
This behavior is consistent with the earlier analysis showing that administrative censoring can prevent the regularity conditions required for strict properness of the Brier score from being fully satisfied, even at large sample sizes, particularly for the $x_2 = 1$ group (\Cref{fig:KM_DGP}).
To assess the robustness of these findings, we repeat the experiment for varying test set sizes $n_{\mathrm{test}} \in \{10, 25, 50, 100, 250, 500\}$.
Results for the low censoring task are shown in \Cref{fig:sens_vary_n} (analogous results for the high censoring task are provided in \Cref{fig:sens_vary_n_high_cens}).

As the test set size decreases, the discriminatory ability of both RCLL and ISBS deteriorates markedly.
Excess loss becomes increasingly variable, leading to greater overlap between models and weaker correspondence with model misspecification.
Nevertheless, RCLL retains better separation between strongly misspecified models—most notably `RSF'—whereas ISBS exhibits substantially greater overlap across the entire range of test-set sizes.
Smaller test sets also increase the frequency of empirical properness violations.
For ISBS, negative excess losses become more pronounced as $n_{\mathrm{test}}$ decreases, whereas for $n_{\mathrm{test}} \ge 100$ such violations are uncommon, consistent with the earlier simulation results (\Cref{tab:sbs_violations}).
In contrast, RCLL remains largely centered above zero, with only occasional negative excess losses for very small sample sizes ($n_{\mathrm{test}} \leq 25$), attributable to finite-sample variability in the empirical log-likelihood.
At the same time, reducing $n_{\mathrm{test}}$ while keeping the integration grid fixed increases the variability of the estimated MIAE (\Cref{eq:MIAE}), producing the wider horizontal spread and reduced separation between models observed in \Cref{fig:sens_vary_n}.

\subsubsection{Sensitivity to Prediction Grid Resolution}
\label{sec:rcll_grid}

The above results indicate that RCLL can provide a more refined discrimination between models with different levels of misspecification, provided the test sample size is sufficiently large.
A natural question is how sensitive this behavior is to the numerical approximation of the density $f(t)$, which is required for RCLL but not for SBS or ISBS.
In previous experiments, densities were obtained by linearly interpolating the predicted survival function over a dense grid of event times, yielding a close approximation to the implied density $\hat f(t) = - \frac{d}{dt}\hat S(t)$.
To assess the impact of this approximation, we progressively reduce the prediction grid by retaining only 80\%, 50\%, 20\%, 10\%, 5\% and 2\% of the unique event times for each task.
Survival prediction and density estimation are both carried out on the reduced grids, mimicking settings with limited temporal resolution (e.g., smaller training samples or discretized time measurements), while MIAE is still computed on the original event-time grid (\Cref{eq:MIAE}).
We generate 100 Monte Carlo test sets with fixed size $n_{\mathrm{test}} = 250$, for which RCLL already exhibits stable empirical behavior under both censoring settings (\Cref{fig:sens_vary_n}), thereby isolating the effect of prediction grid resolution.
Corresponding experiments for $n_{\mathrm{test}}\in\{50,100,500\}$ together with tables of absolute RCLL scores are provided in the Supplementary Material (\Cref{sec:sup}).

\begin{figure}[ht]
\centering
\includegraphics[width=1.0\textwidth]{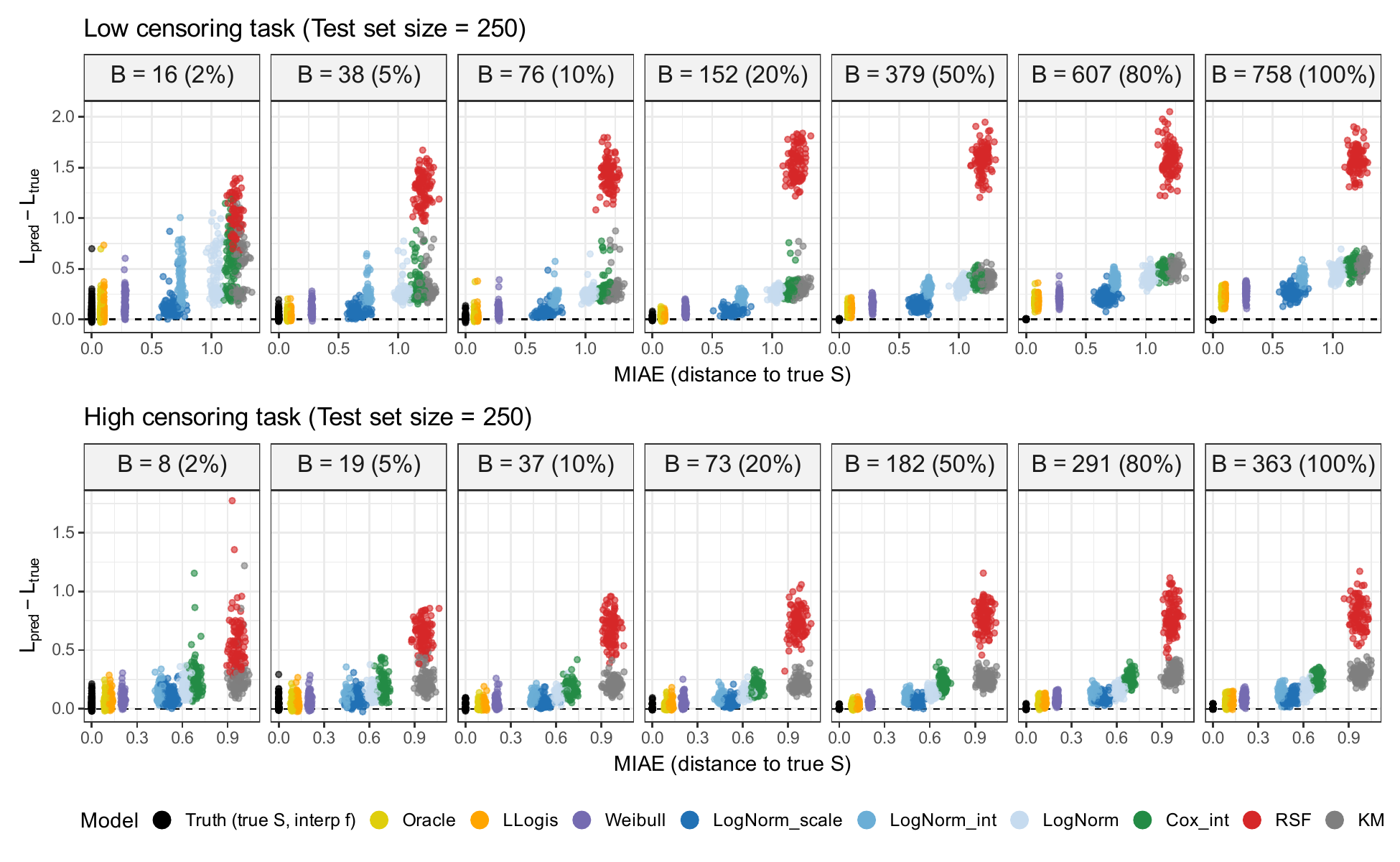}
\caption{
Effect of prediction grid resolution on sensitivity to model misspecification for RCLL under low (top) and high (bottom) censoring.
Each panel corresponds to a grid of size $B$, defined by subsampling a proportion of the unique event times from the corresponding training task (increasing from left to right).
Points represent 100 Monte Carlo replicates ($n_{\mathrm{test}} = 250$ per replicate) and show the excess loss $\Delta L = L_{\mathrm{pred}}(\hat S) - L_{\mathrm{true}}(S_Y)$ versus the mean integrated absolute error (MIAE) between $\hat S$ and $S_Y$.
Models are ordered along the horizontal axis by increasing discrepancy from the true survival function (larger MIAE).}
\label{fig:sens_rcll_grid}
\end{figure}

To isolate the effect of prediction grid resolution on the density approximation underlying the RCLL, excess loss $\Delta L$ (\Cref{eq:excess_loss}) is computed relative to the exact right-censored likelihood obtained from the true survival function $S_Y$ and true density $f_Y$.
The interpolated reference used in the previous experiments is therefore included as an additional model, denoted `Truth (true $S$, interp $f$)'.
\Cref{fig:sens_rcll_grid} shows that, provided all models are evaluated on the same prediction grid, grid resolution has little influence on their relative ranking.
Across all grid sizes, the ordering of the fitted models remains largely unchanged under both censoring settings, although very coarse grids ($\le 20\%$ of the original event times) increase the variability of $\Delta L$ and slightly reduce separation between neighboring models, particularly under high censoring.
Moreover, $\Delta L$ remains almost uniformly positive across all grid resolutions, indicating that empirical properness is largely preserved despite reduced temporal resolution.

While model ranking remains largely unaffected, the absolute RCLL values vary systematically with the prediction grid resolution, as shown in \Cref{tab:rcll_grid} for the low-censoring task and \Cref{tab:rcll_grid_high_cens} for the high-censoring task.
As expected, the exact likelihood (`Truth (true $S$, true $f$)') is essentially invariant to grid resolution, with only negligible Monte Carlo variation.
In contrast, the interpolated reference (`Truth (true $S$, interp $f$)') exhibits a systematic upward bias on coarse grids, which is also visible as positive excess loss in \Cref{fig:sens_rcll_grid}.
For the low-censoring task, the mean RCLL increases from approximately 1.50 on the full grid to 1.58 when only 2\% of the event times are retained.
This approximation bias decreases monotonically as the grid becomes denser and is essentially eliminated once at least 50\% of the unique event times are used.
The same qualitative behavior is observed for the interpolated reference model under high censoring (\Cref{tab:rcll_grid_high_cens}), although the approximation error decreases more slowly and remains slightly visible even on the full prediction grid.

Inspection of the absolute RCLL scores in \Cref{tab:rcll_grid,tab:rcll_grid_high_cens} reveals a second systematic effect of prediction grid resolution.
While the interpolated reference converges to the exact likelihood as the grid becomes denser, all fitted models exhibit the opposite trend: coarsening the prediction grid systematically decreases their RCLL scores, making them appear artificially closer to the true model.
This effect is present under both censoring settings, although it is substantially stronger under low censoring, where the event-density component contributes more heavily to the RCLL.
For example, in the low-censoring task, the RSF score drops from $3.073$ (full grid) to $2.5$ (2\% grid), a relative improvement of $18.6\%$; the Oracle parametric model decreases from $1.722$ to $1.583$ ($8.1\%$ improvement); and even the Cox model and Kaplan–Meier estimator show smaller but consistent improvements of $1.3\%$ and $1.5\%$ respectively (\Cref{tab:rcll_grid}).
Under higher censoring, the same tendency remains visible for the more strongly misspecified models, whereas models closer to the truth exhibit considerably smaller and less systematic changes across grid resolutions, with interpolation bias becoming less consistent overall (\Cref{tab:rcll_grid_high_cens}).

\begin{table}[h]
\scriptsize
\centering
\caption{RCLL scores across 100 Monte Carlo replicates for varying time grid sizes $B$ in the low-censoring task with $n_{\mathrm{test}} = 250$.
Values are reported as mean $\pm$ standard deviation.
Grid sizes correspond to proportions of the full set of unique event times ($B = 758$).
`Truth (true $S$, true $f$)' denotes the exact right-censored likelihood.}
\label{tab:rcll_grid}
\renewcommand{\arraystretch}{1.1}
\begin{tabularx}{\textwidth}{>{\raggedright\arraybackslash}l C C C C C C}
\toprule
Model & $B = 16$ (2\%) & $B = 38$ (5\%) & $B = 76$ (10\%) & $B = 152$ (20\%) & $B = 379$ (50\%) & $B = 758$ (100\%) \\
\midrule
Truth (true $S$, true $f$)
& \makecell{1.500 $\pm$ 0.058}
& \makecell{1.501 $\pm$ 0.057}
& \makecell{1.503 $\pm$ 0.055}
& \makecell{1.502 $\pm$ 0.062}
& \makecell{1.501 $\pm$ 0.053}
& \makecell{1.501 $\pm$ 0.060} \\
Truth (true $S$, interp $f$)
& \makecell{1.581 $\pm$ 0.106}
& \makecell{1.531 $\pm$ 0.075}
& \makecell{1.517 $\pm$ 0.064}
& \makecell{1.509 $\pm$ 0.065}
& \makecell{1.501 $\pm$ 0.053}
& \makecell{1.501 $\pm$ 0.059} \\
Oracle
& \makecell{1.583 $\pm$ 0.108}
& \makecell{1.535 $\pm$ 0.076}
& \makecell{1.531 $\pm$ 0.074}
& \makecell{1.533 $\pm$ 0.073}
& \makecell{1.597 $\pm$ 0.062}
& \makecell{1.722 $\pm$ 0.080} \\
LLogis
& \makecell{1.589 $\pm$ 0.113}
& \makecell{1.537 $\pm$ 0.075}
& \makecell{1.532 $\pm$ 0.074}
& \makecell{1.532 $\pm$ 0.067}
& \makecell{1.602 $\pm$ 0.061}
& \makecell{1.728 $\pm$ 0.081} \\
Weibull
& \makecell{1.668 $\pm$ 0.126}
& \makecell{1.591 $\pm$ 0.083}
& \makecell{1.580 $\pm$ 0.084}
& \makecell{1.588 $\pm$ 0.081}
& \makecell{1.650 $\pm$ 0.065}
& \makecell{1.772 $\pm$ 0.081} \\
LogNorm\_scale
& \makecell{1.620 $\pm$ 0.111}
& \makecell{1.579 $\pm$ 0.078}
& \makecell{1.581 $\pm$ 0.077}
& \makecell{1.581 $\pm$ 0.068}
& \makecell{1.639 $\pm$ 0.059}
& \makecell{1.773 $\pm$ 0.079} \\
LogNorm\_int
& \makecell{1.849 $\pm$ 0.204}
& \makecell{1.730 $\pm$ 0.100}
& \makecell{1.710 $\pm$ 0.070}
& \makecell{1.720 $\pm$ 0.068}
& \makecell{1.792 $\pm$ 0.065}
& \makecell{1.913 $\pm$ 0.075} \\
LogNorm
& \makecell{1.904 $\pm$ 0.233}
& \makecell{1.771 $\pm$ 0.105}
& \makecell{1.756 $\pm$ 0.068}
& \makecell{1.768 $\pm$ 0.063}
& \makecell{1.834 $\pm$ 0.062}
& \makecell{1.964 $\pm$ 0.074} \\
Cox\_int
& \makecell{2.006 $\pm$ 0.259}
& \makecell{1.813 $\pm$ 0.142}
& \makecell{1.797 $\pm$ 0.112}
& \makecell{1.809 $\pm$ 0.082}
& \makecell{1.895 $\pm$ 0.058}
& \makecell{2.032 $\pm$ 0.070} \\
KM
& \makecell{2.010 $\pm$ 0.266}
& \makecell{1.812 $\pm$ 0.142}
& \makecell{1.802 $\pm$ 0.119}
& \makecell{1.817 $\pm$ 0.087}
& \makecell{1.901 $\pm$ 0.056}
& \makecell{2.041 $\pm$ 0.069} \\
RSF
& \makecell{2.500 $\pm$ 0.190}
& \makecell{2.794 $\pm$ 0.174}
& \makecell{2.936 $\pm$ 0.153}
& \makecell{3.042 $\pm$ 0.179}
& \makecell{3.080 $\pm$ 0.161}
& \makecell{3.073 $\pm$ 0.148} \\
\bottomrule
\end{tabularx}
\end{table}

The opposite biases observed for the interpolated reference and the fitted models can be understood from the way densities are approximated on discrete prediction grids.
Since $f(t) = -S'(t)$, densities on a discrete time grid are approximated by finite differences, yielding a piecewise-constant density within each interval \citep{Kvamme2021ContinuousNetworks}:

\[
f(t) \approx f_i = - \frac{S(t_{i}) - S(t_{i-1})}{t_{i} - t_{i-1}}, \qquad t \in [t_{i-1},t_i).
\]

When neighboring intervals are merged (e.g. $[t_{i-1},t_i]$ and $[t_i,t_{i+1}]$), the coarse-grid density becomes a weighted average of the corresponding fine-grid densities:

\[
\begin{aligned}
f_{\mathrm{coarse}}
&= -\frac{S(t_{i+1}) - S(t_{i-1})}{t_{i+1} - t_{i-1}} \\
&= \frac{t_i - t_{i-1}}{t_{i+1} - t_{i-1}}
\left(-\frac{S(t_i) - S(t_{i-1})}{t_i - t_{i-1}}\right)
\;+\;
\frac{t_{i+1} - t_i}{t_{i+1} - t_{i-1}}
\left(-\frac{S(t_{i+1}) - S(t_i)}{t_{i+1} - t_i}\right) \\
&= \frac{t_i - t_{i-1}}{t_{i+1} - t_{i-1}}\, f_i
\;+\;
\frac{t_{i+1} - t_i}{t_{i+1} - t_{i-1}}\, f_{i+1} \\
&= w_i\,f_i + w_{i+1}\,f_{i+1},
\end{aligned}
\]

with $w_i,w_{i+1} \ge 0$ and $w_i + w_{i+1} = 1$.
Thus, grid coarsening replaces local densities by averages over wider time intervals.

Consider an observed event occurring within one of the two intervals that are merged by grid coarsening.
Its contribution to the event component of the RCLL is $-\log f_i$ or $-\log f_{i+1}$ on the fine grid, depending on the interval in which the event occurs, whereas on the coarsened grid the contribution becomes $-\log f_{\mathrm{coarse}}$.
Because $f_{\mathrm{coarse}}$ is a convex combination of the fine-grid densities,

\[
\min(f_i,f_{i+1}) \;\le\; f_{\mathrm{coarse}} \;\le\; \max(f_i,f_{i+1}).
\]

Since the function $-\log(\cdot)$ is strictly decreasing, we have:

\[
-\log\!\bigl(\max(f_i,f_{i+1})\bigr) \;\le\; -\log f_{\mathrm{coarse}} \;\le\; -\log\!\bigl(\min(f_i,f_{i+1})\bigr).
\]

Consequently, replacing the fine grid by the coarse grid increases the event contribution whenever the event occurs in the higher-density subinterval, and decreases it whenever it occurs in the lower-density subinterval.
Therefore the direction of the average RCLL bias depends on where the true event-time distribution places its probability mass within each merged interval.

The empirical results are consistent with this mechanism.
For the interpolated reference model, the predicted density closely follows the true event density, so observed events tend to occur in regions where the predicted density is locally high.
Averaging neighboring intervals therefore lowers the density assigned to these typical event times, increasing the RCLL and producing the observed upward bias.
In contrast, misspecified models do not accurately reproduce the local shape of the true event density.
Their predicted densities are therefore often lower than the interval average at the times where events actually occur, so grid coarsening increases the effective density assigned to those events, reducing the RCLL and making these models appear artificially closer to the truth.
\\\\
Overall, these results show that the absolute scale of the RCLL is sensitive to prediction grid resolution, with both the magnitude and direction of the resulting bias depending on the interaction between the predicted density and the true event-time distribution.
Nevertheless, provided all models are evaluated on the same prediction grid and the grid is sufficiently dense, relative model ranking and empirical properness remain largely preserved.
The effect of grid resolution is most pronounced under low censoring, where the event-density component contributes more strongly to the RCLL, and becomes progressively weaker as censoring increases.

Finally, \Cref{tab:rcll_grid,tab:rcll_grid_high_cens} also illustrate an important practical property of the RCLL that is independent of prediction grid resolution.
Across all models, absolute RCLL values are substantially smaller under high censoring than under low censoring, even though the underlying event-time distribution is identical and only the censoring distribution differs between the two simulation settings.
A plausible explanation is that a higher censoring rate shifts more observations from the density contribution, $-\log f(t)$, to the survival contribution, $-\log S(t)$.
In our simulations, censoring predominantly occurs while the survival probability remains relatively high (\Cref{fig:KM_DGP}), so the survival contribution is typically much smaller than the density contribution for the better-performing models.
This is also reflected in the numerical values of the predicted densities and survival probabilities on the considered time scale, where $f(t)$ is generally much smaller than $S(t)$, resulting in substantially larger values of $-\log f(t)$ than $-\log S(t)$.
Consequently, absolute RCLL values are most informative for comparing models within the same dataset or simulation setting.
For benchmarks involving datasets with different censoring rates, comparisons based on within-dataset model rankings are more appropriate than comparisons of the absolute RCLL values themselves.

\section{Discussion and Conclusion}\label{sec:discussion}

Survival analysis plays a central role in fields such as healthcare, finance, and engineering, where model-based decisions can have direct consequences for public well-being.
With the increasing adoption of machine learning methods for survival analysis\cite{Wang2019, ODonnell2025, MLSA2026}—particularly complex and often black-box models such as deep neural networks\cite{Wiegrebe2024}—robust approaches to external validation have become essential.
In this context, predictions are typically probabilistic, describing the full distribution of time-to-event outcomes, and therefore require evaluation tools that go beyond discrimination measures.
Scoring rules provide such a framework by assessing the quality of predicted distributions in a principled way.
In particular, proper scoring rules are designed to reward predictions that are close to the true data-generating process, while strictly proper scoring rules uniquely favor the correct distribution.
Among the wide range of available methods, squared and logarithmic losses remain the most commonly used in survival analysis, reflecting their established role in regression and classification settings.
\\\\
In this work, we formalize marginal properness for survival scoring rules directly in terms of the observable outcomes $(T,\Delta)$ under independent censoring $Y \indep C$, providing a conceptually transparent framework aligned with external validation settings.
Within this framework, we conduct a systematic re-evaluation of commonly used survival losses, distinguishing between exact formulations (SCRPS, NLL, RCLL) and IPCW-based scores (SBS, ISBS, IBLL, ISAS).
Our theoretical results yield a clear characterization of when these scoring rules are proper, highlighting a previously underappreciated dependence on tail behavior.
In particular, we show that the Survival Brier Score (SBS), Integrated Survival Brier Score (ISBS)\cite{Graf1999}, and Right-Censored Log-Loss (RCLL)\cite{Avati2020} are strictly proper under standard regularity conditions, but can become improper when these are violated—most notably in the presence of non-vanishing survival tails, as encountered under finite follow-up or in settings with cure fractions.

Complementing the theoretical analysis, we conduct simulation experiments to examine how these properness properties translate into practice.
The first experiment (Sections \ref{sec:large-n-properness} and \ref{sec:robustness}) introduces a flexible simulation framework for assessing empirical properness and detecting violations across a wide range of survival and censoring distributions.
Using Weibull‑generated data with independent censoring, we evaluate the RCLL, SBS, and ISBS under finite sample sizes, focusing on regimes where residual tail mass arises naturally from limited follow‑up.
The second experiment (Section \ref{sec:sensitivity}) investigates sensitivity to model misspecification.
Here we generate data from a non‑proportional hazards log‑normal parametric model with crossing survival curves, and compare the true survival function against a sequence of increasingly misspecified predictions—ranging from omitted interactions or misspecified scale parameters to semi‑parametric and non‑parametric models.
This design reflects realistic deviations from modeling assumptions and places special emphasis on finite‑sample test sets, where model evaluation is performed in practice.

The results show that the SBS is highly sensitive to tail behavior and can exhibit pronounced improperness, especially at later evaluation times.
This aligns with our theoretical findings, where improperness is driven by the residual tail mass $\varepsilon$ and the choice of evaluation point $\tau^*$—through the interaction between $S_Y(\tau^*)$ and $S_C(\tau^*)$—and translates into poor discrimination between models with different degrees of misspecification.
The ISBS is more stable due to temporal integration, but still exhibits finite-sample deterioration and reduced discriminatory power under misspecification, requiring moderate sample sizes (approximately $n \ge 50-100$) for empirical properness violations to become rare.
In contrast, the RCLL shows consistently robust behavior across settings, maintaining properness even in the presence of finite tail mass while providing better separation between competing models.
Although the absolute scale of the RCLL depends on the prediction grid resolution and the censoring rate, these effects have little influence on relative model ranking when all models are evaluated under the same conditions.
These observations are further supported by an additional benchmark on eight real-world survival datasets (Supplementary Material; \Cref{sec:sup}), where the RCLL consistently provided greater separation between competing learners than other metrics.
The benchmark likewise illustrates that absolute RCLL values vary substantially with the censoring rate, reinforcing the notion that the score is most informative for comparing models within the same dataset or evaluation setting.
For multi-dataset benchmarks, comparisons based on within-dataset rankings rather than absolute RCLL values are therefore more appropriate.

Based on these findings, we advise against using the pointwise Survival Brier Score (SBS) in practice, especially at later evaluation times, where residual tail mass can induce pronounced improperness.
The ISBS can be used reliably when the integration horizon is truncated (e.g., up to the 80th–90th percentile of observed follow-up) to mitigate tail instability, and when test sample sizes are at least moderate ($n \ge 50-100$).
Finally, the RCLL offers a theoretically grounded and practically robust alternative, provided that predicted densities are approximated on a sufficiently dense and consistent evaluation grid across all compared models, and that the density‑positivity condition is satisfied through appropriate interpolation or smoothing.
More broadly, our findings suggest that survival evaluation metrics with similar theoretical properties can exhibit markedly different empirical behavior.
This highlights the need for further work that complements theoretical analyses of properness with systematic empirical studies of discrimination, robustness and numerical stability, ultimately providing clearer guidance on metric selection for survival prediction.
\\\\
A core limitation of this work is the assumption that $Y$ and $C$ are independent.
This restriction is common in survival analysis evaluation due to the fundamental non‑identifiability problem\cite{Tsiatis1975}: when censoring and event times are dependent, accurate estimation of survival distributions becomes exceedingly challenging, if not impossible.
In practice, researchers often assume conditional independence ($Y \indep C \mid X$), but the extent to which either assumption holds in real‑world settings remains largely unverified, and the robustness of scoring rules under violations of independence is underexplored.
If independence (marginal or conditional) is violated, no loss function is generally proper, rendering any evaluation difficult to interpret.
Recent guidelines have begun to formalize this concern: Kragh Andersen et al. \cite{KraghAndersen2021} emphasize that, ideally, datasets should document why each subject was censored, and analysts should justify the assumed censoring mechanism.
Yet in practice, the true nature of censoring is often unknown, leaving the independence assumption unverifiable and its potential violations as a persistent source of bias.

The IPCW losses studied here further require estimation of the censoring survival function $S_C$, typically via the Kaplan–Meier (KM) estimator\cite{Kaplan1958} or a conditional model (e.g., Cox regression\cite{cox1972regression}).
While the KM estimator is consistent under marginal independent censoring, many practical benchmark studies \cite{Jaeger2024, Wissel2023, Kvamme2023, Li2024, Stuber2023ACRC, Bommert2022BenchmarkData, Herrmann2021, Kantidakis2020SurvivalTechniques, Haider2020, Burk2026AData, Wissel2025} assume only conditional independent censoring and nevertheless use the unconditional KM estimator, accepting the resulting approximation as a practical compromise.
However, when $Y$ and $C$ are only conditionally independent given covariates, the KM estimator is generally biased \cite{Gerds2006}, and choosing a suitable conditional estimator is non‑trivial: it introduces a recursive dependence between the evaluation loss and the censoring model, and becomes particularly challenging in high‑dimensional settings where covariates far exceed the number of samples ($p \gg n$)\cite{Kvamme2023}.
Only recently have researchers begun to investigate whether the censoring distribution should be estimated from training, test, or the combined data to reduce evaluation bias, with test‑set estimation showing advantages under covariate‑dependent censoring or small sample sizes \cite{Prince2025}.
Emerging work based on copula models \cite{Czado2022, Deresa2025} offers a promising direction to relax the independence assumption.
Extending our properness analysis to conditional independent censoring and to settings with dependent censoring constitutes an important avenue for future work.

Our simulation studies, while informative, have several design limitations.
The first experiment relies exclusively on Weibull distributions, which have vanishing tails by construction; future empirical work should therefore incorporate a wider range of survival distributions—such as log‑normal, Gompertz, or mixture cure models—and include administrative censoring by default to cover more realistic evaluation scenarios.
The second experiment uses a single non‑proportional hazards data-generating process (a log-normal accelerated failure time model with crossing survival curves), chosen to represent a challenging scenario where misspecification is particularly consequential.
The relative ranking of scoring rules, especially with regard to their sensitivity to model misspecification, may differ under other data‑generating mechanisms, such as proportional hazards, weaker or non‑crossing non‑proportional alternatives, heavy‑tailed event times, or time‑varying coefficients.
Extending the range of data-generating processes would therefore further strengthen the generalizability of our empirical conclusions.

Our results demonstrate that the RCLL exhibits strong theoretical and empirical properties, yet it remains less widely adopted than the ISBS.
One barrier is the need to estimate densities from survival or hazard functions, as many survival models do not provide them directly.
Future work should therefore investigate density estimation techniques, including improved interpolation and smoothing methods, that reduce the sensitivity of the absolute RCLL to prediction-grid resolution while preserving computational efficiency.
This direction complements the findings of Yanagisawa \cite{Yanagisawa2023}, who observed that a relatively coarse discretization (16-32 time points) was sufficient for their discrete approximation of the RCLL to closely reproduce the full logarithmic score.
Our results provide additional insight into this observation by showing that coarse prediction grids primarily affect the absolute RCLL values through interpolation bias in the estimated densities, while having comparatively little influence on the relative ranking of competing models when all models are evaluated on the same prediction grid.

Beyond the single‑event, right‑censored setting, survival analysis frequently involves competing risks, recurrent events, and multi‑state models \cite{Spitoni2018PredictionModels}.
The IPCW-based scoring rules examined here have already been applied to cause‑specific outcomes in competing risks \cite{Lee2018, Kretowska2018, Bender2021AAnalysis, Alberge2025}, making our properness findings directly relevant.
Future work should systematically investigate proper losses tailored for these multi‑event settings, as well as for left‑ and interval‑censored data—areas that remain largely unexplored.
\\\\
Our findings offer reassurance: despite theoretical improperness under certain regularity violations, metrics such as the ISBS and RCLL remain empirically reliable when used with adequate sample sizes, truncated evaluation horizons (for ISBS), and consistent prediction grids for density approximation (for RCLL).
For the RCLL, a sufficiently dense grid improves the accuracy of absolute score values, while relative model rankings remain largely stable across grid resolutions.
Moreover, because the absolute scale of the RCLL depends on the censoring rate, comparisons across datasets with different censoring levels are best based on relative model rankings rather than raw score values.
Thus, these metrics remain appropriate for model validation in future research, provided these practical guidelines are followed.
Overall, prior studies relying on ISBS‑based evaluation remain largely trustworthy, while underscoring the need to carefully account for censoring and tail behavior.
To ensure robust and well‑interpreted results, we recommend that future work continue using these scoring rules alongside complementary measures such as concordance indices and calibration metrics \cite{Zhao2024, Lillelund2025StopModels, Steyerberg2010}.

\section*{Financial disclosure}

JZ received funding from the European Union’s
Horizon 2020 research and innovation program (grant n.\ 101016851), project PANCAIM.
RDB has received support from the Research Council of Norway (RCN) through the FRIPRO PLUMBIN’ (grant n.\ 323985) and the Centre of Excellence Integreat (grant n.\ 332645).

\section*{Conflict of interest}

The authors declare no potential conflict of interests.

\section*{Supporting information}
\phantomsection
\label{sec:sup}
R scripts and data for the simulation studies in Section \ref{sec:experiments} are available online at \url{https://github.com/survival-org/scoring-rules-2024}.
The repository includes all code, results, supplementary tables and figures, and an HTML report summarizing the findings, available at \url{https://survival-org.github.io/scoring-rules-2024/}.
Reproducibility details and instructions on how to run the analyses are also included as part of the report.

\section*{Acknowledgements}
\phantomsection
The authors would like to thank the editor and two referees for comments that improved both the clarity of the writing and the quality of the
research.

\section*{Declaration of generative AI and AI-assisted technologies in the writing process}
\phantomsection
During the preparation of this work, the authors used DeepSeek (V3) and ChatGPT (GPT-5) to assist with language refinement, structural editing, and proof verification of the mathematical derivations.
After using these tools, the authors carefully reviewed and revised all content and take full responsibility for the final version of the manuscript.

\bibliographystyle{unsrtnat}
\bibliography{library}

\newpage

\appendix
\numberwithin{table}{section}
\numberwithin{figure}{section}

\section{Loss Expectation under Independent Censoring}
\label{sec:appendixA}

We begin by providing some definitions and lemmas that are used in the proof of \Cref{prop:exp}.

\allowdisplaybreaks

\begin{Def}
Let $X$ be an absolutely continuous random variable and let $\Delta$ be a discrete random variable. Then,
\begin{enumerate}[i)]
\item The \emph{mixed joint density} of $(X,\Delta)$ is defined by
\begin{equation}
\label{eq:mixed_joint_dens}
f_{X,\Delta}(x,\delta) = f_{X\mid\Delta}(x\mid\delta)P(\Delta = \delta)
\end{equation}
where $f_{X\mid\Delta}(x\mid\delta)$ is the conditional probability density function of $X$ given $\Delta = \delta$.
\item The \emph{mixed joint cumulative distribution function} of $(X,\Delta)$ is given by
\begin{equation}
F_{X,\Delta}(x,\delta) =  \sum_{\delta' \leq \delta} \int_{-\infty}^x f_{X,\Delta}(u, \delta') \ du
\end{equation}
\end{enumerate}
\end{Def}

\begin{Lem}
\label{lem:dens_Y_D}
Let $Y,C$ be jointly absolutely continuous random variables supported on the Reals with joint density function $f_{Y,C}(y,c)$ and let $\Delta = \II(Y \leq C)$, then the mixed joint density of $(Y,\Delta)$ is given by

\begin{equation}
\label{eq:dens_Y_D}
f_{Y,\Delta}(y,\delta) =
\begin{cases}
\int^\infty_y \ f_{Y,C}(y,c) \ dc, & \delta = 1 \\
\int^y_{-\infty} f_{Y,C}(y,c) \ dc, & \delta = 0
\end{cases} \\
\end{equation}

In addition if $Y \indep C$, then

\begin{equation}
\label{eq:dens_Y_D_indep}
f_{Y,\Delta}(y,\delta) =
\begin{cases}
f_Y(y)\,S_C(y), & \delta = 1 \\
f_Y(y)\,F_C(y), & \delta = 0
\end{cases}
\end{equation}

\end{Lem}
\begin{proof}
Proof follows by transformation of random variables via the joint cumulative distribution function.

The joint cumulative distribution function of $(Y,C)$ is defined by,
\[
F_{Y,C}(y,c) = \int^y_{-\infty}\int^c_{-\infty} f_{Y,C}(s,t) \ dt \ ds
\]

By definition of indicator variables, $\Delta = 1$ if and only if $C \geq Y$ and $0$ otherwise and so:

\[
\begin{aligned}
F_{Y,\Delta}(y,\delta)
&=P(Y \leq y, \Delta \leq \delta) \\
&= \begin{cases}
P(Y \leq y), & \delta = 1 \\
P(Y \leq y, C < Y), & \delta = 0 \\
\end{cases}
\end{aligned}
\]

where the first case follows as $\Delta \in \{0,1\}$ and hence $P(Y \leq y, \Delta \leq 1) = P(Y \leq y)$ as $\Delta$ is marginalised out. The second case follows as $\Delta \in \{0,1\}$ and so $P(\Delta \leq 0) = P(\Delta = 0)$, and by definition of indicator variables $\Delta = 0$ if and only if $C < Y$. Now focusing on the second case,

\[
\begin{aligned}
P(Y \leq y, C < Y)
&= \iint_{\{t \leq y,c<t\}} f_{Y,C}(t,c) \ dt \ dc \\
&= \int^y_{-\infty}\int^t_{-\infty} f_{Y,C}(t,c) \ dc \ dt
\end{aligned}
\]

The first line follows by definition of joint probabilities and the second by change of notation. Now,

\begin{equation}
\label{eq:lem_joints_cdf}
F_{Y,\Delta}(y,\delta) =
\begin{cases}
F_Y(y), & \delta = 1 \\
\int^y_{-\infty}\int^t_{-\infty} f_{Y,C}(t,c) \ dc \ dt, & \delta = 0 \\
\end{cases}
\end{equation}

The joint density of $(Y,\Delta)$ (\Cref{eq:dens_Y_D}) follows by differentiating the cumulative distribution function:

\[
\begin{aligned}
f_{Y,\Delta}(y,\delta)
&=
\begin{cases}
\frac{\partial \left(F_{Y,\Delta}(y,1) - F_{Y,\Delta}(y,0)\right)}{\partial y}, & \delta = 1 \\
\frac{\partial F_{Y,\Delta}(y,0)}{\partial y}, & \delta = 0
\end{cases} \\
&=
\begin{cases}
\frac{\partial}{\partial y} \left(\int^y_{-\infty}\int^\infty_{-\infty} f_{Y,C}(t,c) \ dc \ dt - \int^y_{-\infty}\int^t_{-\infty} f_{Y,C}(t,c) \ dc \ dt\right), & \delta = 1 \\
\frac{\partial}{\partial y} \left(\int^y_{-\infty}\int^t_{-\infty} f_{Y,C}(t,c) \ dc \ dt\right), & \delta = 0
\end{cases} \\
&=
\begin{cases}
\frac{\partial}{\partial y} \left(\int^y_{-\infty}\left(\int^\infty_{-\infty} f_{Y,C}(t,c) \ dc - \int^t_{-\infty} f_{Y,C}(t,c) \ dc \right) \ dt\right), & \delta = 1 \\
\int^y_{-\infty} f_{Y,C}(y,c) \ dc, & \delta = 0
\end{cases} \\
&=
\begin{cases}
\frac{\partial}{\partial y} \left( \int^y_{-\infty}\int^\infty_t f_{Y,C}(t,c) \ dc \ dt\right), & \delta = 1 \\
\int^y_{-\infty} f_{Y,C}(y,c) \ dc, & \delta = 0
\end{cases} \\
&=
\begin{cases}
\int^\infty_y f_{Y,C}(y,c) \ dc, & \delta = 1 \\
\int^y_{-\infty} f_{Y,C}(y,c) \ dc, & \delta = 0
\end{cases} \\
\end{aligned}
\]

where the first equality is the definition of the joint mixed propability density function in terms of the cumulative distribution function, the second equality follows by substituting \Cref{eq:lem_joints_cdf}, the third by taking the partial derivative of the second case over $y$ and by linearity of integration in the first case, the fourth by subtracting the inner integrals, and the fifth by taking the partial derivative of the first case over $y$.

Continuing, if $Y,C$ are independent, then $f_{Y,C}(y,c) = f_Y(y)\,f_C(c)$ by definition of independence.
Substituting this result into \Cref{eq:dens_Y_D}, we have:

\[
\begin{aligned}
f_{Y,\Delta}(y,\delta)
&=
\begin{cases}
\int^\infty_y \ f_Y(y)\,f_C(c) \ dc, & \delta = 1 \\
\int^y_{-\infty} f_Y(y)\,f_C(c) \ dc, & \delta = 0
\end{cases} \\
&=
\begin{cases}
f_Y(y) \int^\infty_y \ f_C(c) \ dc, & \delta = 1 \\
f_Y(y) \int^y_{-\infty} f_C(c) \ dc, & \delta = 0
\end{cases} \\
&=
\begin{cases}
f_Y(y)\,S_C(y), & \delta = 1 \\
f_Y(y)\,F_C(y), & \delta = 0
\end{cases} \\
\end{aligned}
\]

where the first equality holds as $Y,C$ independent, the second by properties of integration, and the third by definition of the cumulative distribution and survival functions.
\end{proof}

\begin{Lem}
\label{lem:dens_Y_C}
Let $Y,C$ be jointly absolutely continuous random variables supported on the Reals with joint density function $f_{Y,C}(y,c)$ and let $\Delta = \II(Y \leq C)$, then the mixed joint density of $(C,\Delta)$ is given by

\begin{equation}
f_{C,\Delta}(c,\delta) =
\begin{cases}
\int^c_{-\infty} f_{Y,C}(y,c) \ dy, & \delta = 1\\
\int^\infty_c \ f_{Y,C}(y,c) \ dy, & \delta = 0
\end{cases} \\
\end{equation}

In addition if $Y \indep C$, then

\begin{equation}
f_{C,\Delta}(c,\delta) =
\begin{cases}
f_C(c)\,F_Y(c), & \delta = 1\\
f_C(c)\,S_Y(c), & \delta = 0
\end{cases} \\
\end{equation}
\end{Lem}

\begin{proof}
Proofs follow analogously to \Cref{lem:dens_Y_D}.
\end{proof}

\begin{proof}[Proof of Proposition~\ref{prop:exp}](\textbf{Expectation of Losses under Independent Censoring})

\[
\scriptsize
\begin{aligned}
\underset{(T,\Delta)}{\EE}\Big[\phi(\hat{S}, T, \Delta)\Big]
&= \sum_{\delta \in \{0,1\}} \int_0^\infty \phi(\hat{S}, t, \delta)\, f_{T,\Delta}(t,\delta)\, dt 
&& \text{definition of joint expectation over } (T,\Delta) \\[1em]
&= \int_0^\infty 
\Big[
P(\Delta=1)\,\phi(\hat{S}, t \mid 1)\, f_{T \mid \Delta}(t \mid 1)
+ P(\Delta=0)\,\phi(\hat{S}, t \mid 0)\, f_{T \mid \Delta}(t \mid 0)
\Big] dt 
&& \text{law of total expectation} \\[1em]
&= \int_0^\infty 
\Big[
f_{T,\Delta}(t,1)\, \phi(\hat{S}, t \mid 1)
+ f_{T,\Delta}(t,0)\, \phi(\hat{S}, t \mid 0)
\Big] dt
&& \text{Equation }\ref{eq:mixed_joint_dens} \text{ with } X=T, \Delta \in \{0,1\}\\
&= \int_0^\infty 
\Big[
f_{Y,\Delta}(t,1)\, \phi(\hat{S}, t \mid 1)
+ f_{C,\Delta}(t,0)\, \phi(\hat{S}, t \mid 0)
\Big] dt
&& Y \equiv T\mid\Delta = 1; C \equiv T\mid\Delta = 0 \\[1em]
&= \int^\infty_0 \left[ f_Y(t)\,S_C(t)\,\phi(\hatS, t \mid 1) + f_C(t)\,S_Y(t)\,\phi(\hatS, t \mid 0) \right] dt
&& \text{\Cref{lem:dens_Y_D} and \Cref{lem:dens_Y_C} as } Y \indep C
\end{aligned}
\]

For the last derivation we note that equivalent expressions for the joint density of the observed data 
$(T,\Delta)$ under independent censoring are given in 
Klein \& Moeschberger \citep{Klein2003SurvivalData}, Example 3.10, p.~76.
\end{proof}

\begin{remark}
\label{rmk:cov_dep_exp}
The prediction $\hat S$ may represent either an unconditional survival distribution or a covariate-dependent prediction $\hat S(\cdot\mid X)$.
In the latter case, the expectation is formally taken with respect to the joint distribution of $(X,T,\Delta)$:
\[
\underset{(T,\Delta,X)}{\EE}\Big[\phi(\hat S(\cdot\mid X), T, \Delta)\Big].
\]

rather than only $(T,\Delta)$.
The derivation of the expectation formula in Proposition~\ref{prop:exp} uses only the marginal observed-data density $f_{T,\Delta}$ and the independent
censoring assumption $Y\indep C$; it does not require conditioning on $X$ or the assumption of conditional independent censoring $Y\indep C\mid X$.
Covariate-dependent predictions are therefore accommodated within the marginal properness framework.
\end{remark}

\section{Properness Proofs}
\label{sec:appendixB}

All notation used in the proofs, unless stated otherwise, follows \Cref{sec:notation}.

\subsection{Properness of Convex Combinations of Losses}

Let $\mathcal{Y}$ denote the outcome space and let $\mathcal{P}$ be a class of probability distributions on $\mathcal{Y}$. 
For a loss function $L$ and distributions $P,Q \in \mathcal{P}$, we define the (expected) risk of predicting $Q$ when the true data-generating distribution is $P$ as
\[
\mathcal{R}_L(P,Q) := \underset{Y \sim P}{\EE}\!\left[L(Q,Y)\right].
\]

A loss $L$ is called \emph{proper} on $\mathcal{P}$ if, for every $P \in \mathcal{P}$,
\[
\mathcal{R}_L(P,P) \le \mathcal{R}_L(P,Q) \qquad \text{for all } Q \in \mathcal{P},
\]
and it is called \emph{strictly proper} if equality holds only when $Q = P$ \cite{Gneiting2007}.
Given the above definitions, we can now proceed to state the following proposition:

\begin{Prop}[\textbf{Properness Preservation for Loss Function Mixtures}]
\label{prop:proper-sum}
Let $L_A$ and $L_B$ be two losses and let $\alpha \in [0,1]$ be a constant independent of $L_A,L_B$, the prediction, and the outcome.
Define $L_\alpha$ as the convex sum of two loss functions:

\[
L_\alpha := \alpha \,L_A + (1 - \alpha)\,L_B
\]

Then,
\begin{enumerate}
\item[(1)] If $L_A$ and $L_B$ are proper, then $L_\alpha$ is proper.
\item[(2)] If one of $L_A$ or $L_B$ is strictly proper and the other is proper, and the strictly proper component receives positive weight---i.e. $\alpha>0$ or $1-\alpha>0$, depending on which component is strictly proper---then $L_\alpha$ is strictly proper.
\end{enumerate}
\end{Prop}

\begin{proof}
\textbf{(1) Properness.} 
Fix $P,Q \in \mathcal{P}$.
By independence of $\alpha$ and linearity of expectation,
\begin{equation}
\label{eq:mixture_risk}
\mathcal{R}_{L_\alpha}(P,Q) = \alpha \, \mathcal{R}_{L_A}(P,Q) + (1-\alpha) \, \mathcal{R}_{L_B}(P,Q).    
\end{equation}

Since $L_A$ and $L_B$ are proper,
\[
\mathcal{R}_{L_A}(P,P) \le \mathcal{R}_{L_A}(P,Q), 
\qquad
\mathcal{R}_{L_B}(P,P) \le \mathcal{R}_{L_B}(P,Q).
\]
Multiplying by the nonnegative weights $\alpha$ and $1-\alpha$ and adding yields
\[
\mathcal{R}_{L_\alpha}(P,P) \le \mathcal{R}_{L_\alpha}(P,Q),
\]
so $L_\alpha$ is proper.

\textbf{(2) Strict properness.}
Assume without loss of generality that $L_A$ is strictly proper and $L_B$ is proper, and that $\alpha>0$.
Fix $P \in \mathcal{P}$ and let $Q \ne P$.
Define
\[
\Delta_A := \mathcal{R}_{L_A}(P,Q) - \mathcal{R}_{L_A}(P,P) > 0,
\qquad
\Delta_B := \mathcal{R}_{L_B}(P,Q) - \mathcal{R}_{L_B}(P,P) \ge 0.
\]
Then, using \Cref{eq:mixture_risk}, the risk difference of the mixture is
\[
\mathcal{R}_{L_\alpha}(P,Q) - \mathcal{R}_{L_\alpha}(P,P) = \alpha \, \Delta_A + (1-\alpha) \, \Delta_B \ge \alpha \, \Delta_A > 0.
\]
Thus
\[
\mathcal{R}_{L_\alpha}(P,Q) > \mathcal{R}_{L_\alpha}(P,P)
\qquad \text{for all } Q \ne P,
\]
establishing strict properness for $L_\alpha$.
\end{proof}

\subsection{SCRPS is not proper}
\label{proof:scrps_not_proper}

\begin{proof}
We demonstrate that the SCRPS scoring rule is not marginally proper by constructing a counterexample where the expected loss under a misspecified prediction $\hatS$ is lower than under the correctly specified distribution $S_Y$.

Let $\mathcal{T} \subseteq \NNReals$ and let $\calP$ be a family of absolutely continuous distributions over $\mathcal{T}$ containing at least two elements. 
Let $\xi, \upsilon \in \calP$ denote distributions over $\mathcal{T}$.
Let $Y \sim \xi$ be the true survival time, and let $C \indep Y$ be an independent censoring time with support on $\mathcal{T}$. 
Define the observed time and censoring indicator as $T := \min\{Y, C\}$ and $\Delta := \II(Y \le C)$.
Let $\hat Y \sim \upsilon$ denote a predicted distribution, and let $\hat{S}$ and $\hat{F}$ be its survival and cumulative distribution functions, respectively.

Taking the expectation of $L_{SCRPS}$ \eqref{eq:SCRPS}, we have:

\[
\small
\begin{aligned}
\underset{(T,\Delta)}{\EE}\Big[L_{\text{SCRPS}}(\hat S, T, \Delta)\Big] &= \int^\infty_0 f_Y(t)S_C(t)\left(\int^t_0 \hat{F}^2(\tau) \ d\tau + \int^\infty_t \hat{S}^2(\tau) \ d\tau\right) \, dt \ + \\
&\quad \int^\infty_0 f_C(t)S_Y(t)\left(\int_0^t \hat{F}^2(\tau) \ d\tau\right) dt && \text{\Cref{prop:exp} as } Y \indep C \\
&= \int^\infty_0 f_Y(t)S_Y(t)\left(2\int^t_0 \hat{F}^2(\tau) \ d\tau + \int^\infty_t \hat{S}^2(\tau) \ d\tau\right) \ dt && \text{let both } C \sim \xi \text{ and } Y \sim \xi\\
&= \int^\infty_0 e^{-2t} \left(2\int^t_0 (1 - e^{-\mu\tau})^2\ d\tau + \int^\infty_t e^{-2\mu\tau} \ d\tau\right) \ dt && \text{let } \xi = \Exp(1), \hat{Y} \sim \Exp(\mu) \\
&= \int^\infty_0 e^{-2t} \left(\frac{-6 - e^{-2 t \mu} + 8 e^{-t \mu} + 4 t \mu}{2\mu} \right) \ dt &&\text{integration} \\
&= \frac{2 \mu^3 + \mu + 2}{4 (\mu + 2) \mu(\mu + 1)} &&\text{integration} \\
\end{aligned}
\]

Let $\mu = 1.5$ so that $\hat Y \sim \Exp(1.5)$ and recall $Y \sim \Exp(\mu = 1)$. Then:

\[
\mathbb{E}[L_{\text{SCRPS}}(\hatS, T, \Delta)] \approx 0.1952 
< 
\mathbb{E}[L_{\text{SCRPS}}(S_Y, T, \Delta)] \approx 0.2083
\]

This violates the definition of marginal properness in \Cref{def:surv_proper}, since the expected loss is lower for a misspecified model.
Thus, $L_{\text{SCRPS}}$ is not a proper scoring rule.
\end{proof}

\subsection{NLL is not proper}
\label{proof:nll_not_proper}

\begin{proof}
We demonstrate that the NLL scoring rule is not marginally proper by constructing a counterexample where the expected loss under a misspecified prediction $\hatS$ is lower than under the correctly specified distribution $S_Y$.

Let $\mathcal{T} \subseteq \NNReals$ and let $\calP$ be a family of absolutely continuous distributions over $\mathcal{T}$ containing at least two elements. 
Let $\xi, \upsilon \in \calP$ denote distributions over $\mathcal{T}$.
Let $Y \sim \xi$ be the true survival time, and let $C \indep Y$ be an independent censoring time with support on $\mathcal{T}$. 
Define the observed time and censoring indicator as $T := \min\{Y, C\}$ and $\Delta := \II(Y \le C)$.
Let $\hat Y \sim \upsilon$ denote a predicted distribution, and let $\hat{S}$ be its survival distribution function, with positive density $\hat{f}$.

Taking the expectation of $L_{NLL}$ \eqref{eq:NLL}, we have:
\[
\footnotesize
\begin{aligned}
\underset{(T,\Delta)}{\EE}\Big[L_{NLL}(\hat S, T, \Delta)\Big] 
&= \int_0^\infty \left[ f_Y(t) S_C(t) \, (-\log \hat f(t)) + f_C(t) S_Y(t) \, (-\log \hat f(t)) \right] dt
&& \text{\Cref{prop:exp} as } Y \indep C \\
&= -2 \int^\infty_0 f_Y(t)S_Y(t) \log(\hat f(t)) \ dt 
&& \text{let both } C \sim \xi \text{ and } Y \sim \xi \\
&= - 2 \int_0^\infty e^{-2t} \log(\mu e^{-\mu t}) \, dt && \text{let } Y \sim \Exp(1),\ \hat Y \sim \Exp(\mu), \mu > 0 \\
&= - 2 \log \mu \int_0^\infty e^{-2t} dt + 2 \mu \int_0^\infty t e^{-2t} dt \\
&= - 2 \log \mu \cdot \frac{1}{2} + 2 \mu \cdot \frac{1}{2^2}
&& \int_0^\infty t^n e^{-\alpha t} dt = \frac{\Gamma(n+1)}{\alpha^{n+1}} \overset{n \in \mathbb{N}}{=} \frac{n!}{\alpha^{n+1}}\\
&= \frac{\mu}{2} - \log \mu.
\end{aligned}
\]

Let $\mu = 2$ so that $\hat Y \sim \Exp(2)$, while the true distribution is $Y \sim \Exp(\mu = 1)$.
For all subsequent calculations, we use the natural logarithm, as is standard in the scoring rules literature.
Then:
\[
\mathbb{E}[L_{NLL}(\hatS, T, \Delta)] 
= \frac{2}{2} - \log 2
\approx 0.3068
< 
\mathbb{E}[L_{NLL}(S_Y, T, \Delta)] 
= \frac{1}{2} - \log 1 
= 0.5.
\]

Since the expected loss is lower for a misspecified prediction $\hatS$, this violates the definition of marginal properness in \Cref{def:surv_proper}.
Thus, $L_{NLL}$ is not a proper scoring rule.
\end{proof}

\subsection{Properness of RCLL}
\label{proof:rcll-proof}

\begin{proof}
Let $Y$ and $C$ be independent, continuous random variables on $[0,\infty)$ with absolutely continuous survival functions $S_Y,S_C$ and densities $f_Y=-S_Y'$, $f_C=-S_C'$ a.e..
Let $\hat Y$ be a predicted survival time with survival function $\hat S$ and density $\hat f=-\hat S'$ a.e.
Assume $f_Y,f_C,\hat f>0$ and $S_Y,S_C,\hat S>0$ for all finite $t\ge0$, and that all functions are integrable.
Under these conditions, the quantities $\log S(t)$ and $\log f(t)$ are well-defined for finite $t$, and all integrals involving them over $[0,\infty)$ exist.
\\\\
Taking the expectation of $L_{RCLL}$ \eqref{eq:RCLL} using \Cref{prop:exp} as $Y \indep C$, we have:
\[
\mathcal{R}(\hat S) := \underset{(T,\Delta)}{\EE}\Big[L_{RCLL}(\hat S, T, \Delta)\Big] = 
-\int_0^\infty \Big[f_Y(t) S_C(t) \log \hat f(t) + f_C(t) S_Y(t) \log \hat S(t)\Big] \, dt
\]

To study RCLL's properness, we write the difference of expected losses, given that $\mathcal{R}(\hat S)$ is the expected RCLL loss under the predicted survival function and $\mathcal{R}(S_Y)$ the expected loss under the true survival function:

\begin{equation}
\label{eq:rcll_exp_diff}
\footnotesize
\begin{aligned}
\mathcal{R}(\hat S) - \mathcal{R}(S_Y)
&= -\int_0^\infty \Big[f_Y(t) S_C(t) \log \hat f(t) + f_C(t) S_Y(t) \log \hat S(t) \Big] \, dt \\
&\quad + \int_0^\infty \Big[f_Y(t) S_C(t) \log f_Y(t) + f_C(t) S_Y(t) \log S_Y(t) \Big] \, dt \\
&= \int_0^\infty \Big[
f_Y(t) S_C(t) \log \frac{f_Y(t)}{\hat f(t)}
+ f_C(t) S_Y(t) \log \frac{S_Y(t)}{\hat S(t)}
\Big] \, dt 
&& \text{combine integrals, use log properties} \\
&= \int_0^\infty \Big[
f_Y(t) S_C(t) \log \frac{f_Y(t) S_C(t)}{\hat f(t) S_C(t)}
+ f_C(t) S_Y(t) \log \frac{f_C(t) S_Y(t)}{f_C(t) \hat S(t)} \Big] \, dt
&& \text{multiply and divide by}\ S_C,f_C>0 \\
&= \int_0^\infty \Big[w_1(t) \log \frac{w_1(t)}{\hat w_1(t)} + w_0(t) \log \frac{w_0(t)}{\hat w_0(t)} \Big] dt
\end{aligned}
\end{equation}

where we introduced the functions $w_1(t), \hat w_1(t), w_0(t), \hat w_0(t)$, defined as:
\begin{equation}
\label{eq:w_true}
w_1(t) := f_Y(t) S_C(t), \qquad w_0(t) := f_C(t) S_Y(t).
\end{equation}
\begin{equation}
\label{eq:w_pred}
\hat w_1(t) := \hat f(t) S_C(t), \qquad \hat w_0(t) := f_C(t) \hat S(t).
\end{equation}

The functions $w_1(t)$ and $w_0(t)$ represent the \emph{mixed joint event-time densities} for the observed data under the true model, corresponding to the likelihood contributions of uncensored and censored observations, respectively (\Cref{lem:dens_Y_D}, \Cref{lem:dens_Y_C}).
Specifically, $w_1(t) = f_Y(t)S_C(t)$ gives the joint probability density that the event occurs at time $t$ and is observed (i.e., $Y=t,\, Y<C$, censoring hasn't happened yet), while $w_0(t) = f_C(t)S_Y(t)$ corresponds to the joint probability density that censoring occurs at time $t$ before the event (i.e., $C=t,\, C<Y$).
Analogously, $\hat w_1(t)$ and $\hat w_0(t)$ denote the corresponding quantities implied by the predicted survival model $\hat S$.
These functions thus decompose the expected loss integral into two additive components—one corresponding to observed events and one to censored cases—each reflecting a separate likelihood component of the right-censored data under the assumption of independence between $Y$ and $C$.
Note that our initial positivity assumptions for the survival and density functions imply $w_1(t), \hat w_1(t), w_0(t), \hat w_0(t)> 0$ for all finite $t \ge 0$ (and also integrable), meaning that both event and censoring processes have nonzero probability across the time domain—that is, some events are observed and some are censored.

We now define the following constant quantities representing the total probabilities of observing an event or a censoring occurrence:

\begin{equation}
\label{eq:Z_true}
Z_1 := \int_0^\infty w_1(t) \, dt = P(Y < C), \qquad Z_0 := \int_0^\infty w_0(t) \, dt = P(C < Y).
\end{equation}

Due to the positivity assumptions, we have that $Z_1, Z_0 > 0$.
In general, these do not sum up to one.
Under the assumptions of independence of $Y$ and $C$, continuity (so ties are not possible, i.e. $P(Y=C)=0$), and allowing for residual mass at infinity in the event and censoring distributions, we obtain:
\begin{equation}
\label{eq:Z_total}
Z := Z_1 + Z_0 = \int_0^\infty \Big[w_1(t)+w_0(t)\Big]\,dt = 1 - \varepsilon.
\end{equation}

The remaining mass $\varepsilon \ge 0$ corresponds to the probability that neither the event nor censoring occurs in finite time, i.e.
\begin{equation}
\label{eq:varepsilon_rcll}
\varepsilon := P(Y = \infty, C = \infty).
\end{equation}

Under independence of $Y$ and $C$ this factorizes as:
\begin{equation}
\varepsilon := P(Y = \infty)\,P(C = \infty) = \varepsilon_Y\,\varepsilon_C.
\end{equation}

where $\varepsilon_Y := \lim_{t\to\infty} S_Y(t) \ge 0$ and $\varepsilon_C := \lim_{t\to\infty} S_C(t) \ge 0$ denote the residual mass of the event time and censoring time distributions, respectively.
We observe that if at least one of the two distributions is proper, i.e. $\varepsilon_Y=0$ or $\varepsilon_C=0$ (regularity condition holds), we recover $Z = 1$ and hence $Z_0+Z_1=1$.

Next, we define analogous quantities for each observation type under the predicted model:

\begin{equation}
\label{eq:Z_pred}
\hat Z_1 := \int_0^\infty \hat w_1(t) \, dt, \qquad \hat Z_0 := \int_0^\infty \hat w_0(t) \, dt.
\end{equation}

Under the same assumptions of continuity and positivity, these satisfy
\begin{equation}
\label{eq:Z_total_pred}
\hat Z := \hat Z_1 + \hat Z_0 = 1 - \hat\varepsilon,
\end{equation}
where
\begin{equation}
\hat\varepsilon := P(\hat Y = \infty,C = \infty)
= \varepsilon_{\hat Y} \, \varepsilon_C.
\end{equation}

Here, $\varepsilon_{\hat Y} := \lim_{t \to \infty} \hat S(t)$ denotes the residual mass of the predicted survival distribution, while $\varepsilon_C$ is defined as above.

From these, we define normalized joint densities over the space of observed data $(t,\delta) \in [0,\infty]\times\{0,1\}$:

\begin{equation}
\label{eq:p_densities}
p(t,\delta) :=
\begin{cases}
\dfrac{w_1(t)}{Z}, & \delta = 1\\[6pt]
\dfrac{w_0(t)}{Z}, & \delta = 0
\end{cases}
\qquad
q(t,\delta) :=
\begin{cases}
\dfrac{\hat w_1(t)}{\hat Z}, & \delta = 1\\[6pt]
\dfrac{\hat w_0(t)}{\hat Z}, & \delta = 0
\end{cases}
\end{equation}

By construction, these functions are nonnegative and integrate to one:
\begin{equation}
\label{eq:norm_densities_rcll}
\sum_{\delta \in \{0,1\}} \int_0^\infty p(t,\delta)\,dt = 1,
\qquad
\sum_{\delta \in \{0,1\}} \int_0^\infty q(t,\delta)\,dt = 1.    
\end{equation}

so both $p$ and $q$ define valid probability densities on the joint space of observed times and observation types.

Now we substitute the above expressions in the expected RCLL loss difference:

\begin{equation}
\footnotesize
\begin{aligned}
\label{eq:rcll_exp_diff_KL}
\mathcal{R}(\hat S) - \mathcal{R}(S_Y)
&= \int_0^\infty \Big[w_1(t) \log \frac{w_1(t)}{\hat w_1(t)} + w_0(t) \log \frac{w_0(t)}{\hat w_0(t)} \Big] dt
&& \text{\Cref{eq:rcll_exp_diff}} \\[4pt]
&= \int_0^\infty \Big[Z\,p(t,1) \log \frac{Z\,p(t,1)}{\hat Z\,q(t,1)} + Z\,p(t,0) \log \frac{Z\,p(t,0)}{\hat Z\,q(t,0)} \Big] dt
&& \text{\Cref{eq:p_densities}} \\[4pt]
&= \int_0^\infty \Big[ Z\,p(t,1)\Big(\log\frac{p(t,1)}{q(t,1)} + \log\frac{Z}{\hat Z}\Big)
+ Z\,p(t,0)\Big(\log\frac{p(t,0)}{q(t,0)} + \log\frac{Z}{\hat Z}\Big)\Big]\,dt 
&& \text{split each log} \\[6pt]
&= Z\int_0^\infty p(t,1)\log\frac{p(t,1)}{q(t,1)}\,dt
+ Z\log\frac{Z}{\hat Z}\int_0^\infty p(t,1)\,dt \\
&\quad + Z\int_0^\infty p(t,0)\log\frac{p(t,0)}{q(t,0)}\,dt
+ Z\log\frac{Z}{\hat Z}\int_0^\infty p(t,0)\,dt
&& \text{distribute and separate integrals} \\[6pt]
&= Z \sum_{\delta \in \{0,1\}} \int_0^\infty p(t,\delta)\log\frac{p(t,\delta)}{q(t,\delta)}\,dt
+ Z\log\frac{Z}{\hat Z} \sum_{\delta \in \{0,1\}} \int_0^\infty p(t,\delta)\,dt
&& \text{collect terms using joint notation} \\[6pt]
&= Z \sum_{\delta \in \{0,1\}} \int_0^\infty p(t,\delta)\log\frac{p(t,\delta)}{q(t,\delta)}\,dt
+ Z\log\frac{Z}{\hat Z}
&& \text{\Cref{eq:norm_densities_rcll}} \\[6pt]
&= Z\,D_{\mathrm{KL}}(p\|q) + Z\log\frac{Z}{\hat Z}
&& \text{identify the KL term}
\end{aligned}
\end{equation}

The first term is a (scaled) Kullback–Leibler (KL) divergence between the normalized joint densities $p$ and $q$ over observed data $(t,\delta)$.
It is non-negative and equals zero if and only if $p(t,\delta) = q(t,\delta)$ almost everywhere.
The second term, reflects differences in the residual (tail) mass at infinity between the true and predicted models, as $Z = 1 - \varepsilon$ \eqref{eq:Z_total} and $\hat Z = 1 - \hat\varepsilon$ \eqref{eq:Z_total_pred}.

We distinguish between two cases, based on the value of $\varepsilon$.

\medskip\noindent
\textbf{Case 1.} If $\varepsilon = \varepsilon_Y \varepsilon_C = 0$ (i.e., at least one tail vanishes, e.g. $\lim_{t\to\infty} S_Y(t)=0$ or $\lim_{t\to\infty} S_C(t)=0$), then $Z = 1$ and the risk difference reduces to

\begin{equation}
\mathcal{R}(\hat S) - \mathcal{R}(S_Y) = D_{\mathrm{KL}}(p\|q) + \log\frac{1}{\hat Z} \;\ge\; 0.
\end{equation}

as $\hat Z \le 1$.
Equality holds if and only if both terms equal zero, i.e.
\begin{equation}
p(t,\delta) = q(t,\delta) \quad \text{ and } \quad \hat Z = 1.
\end{equation}

Substituting back the definitions of $p,q$ \eqref{eq:p_densities} and $\hat Z$ \eqref{eq:Z_total_pred}, this implies

\begin{equation}
\begin{aligned}
& w_1(t) = \hat w_1(t), \quad w_0(t) = \hat w_0(t) \quad \text{for all } t \ge 0, \\
& \hat \varepsilon = 0.
\end{aligned}
\end{equation}

Using \eqref{eq:w_true} and \eqref{eq:w_pred}, this yields
\begin{equation}
f_Y(t) = \hat f(t), \quad S_Y(t) = \hat S(t) \quad \text{for all finite } t \ge 0,
\end{equation}
i.e., the predicted model exactly matches the true survival distribution.
Therefore, when the regularity assumptions hold, the RCLL scoring rule is \emph{strictly proper} under independent censoring, following the marginal properness definition \ref{def:surv_proper}.

\medskip\noindent\textbf{Case 2.} If $\varepsilon = P(Y=\infty, C=\infty) > 0$ (both tails nonzero), then $Z = 1-\varepsilon < 1$.

Since the first term in \eqref{eq:rcll_exp_diff_KL} is always $\ge 0$ (KL divergence and $Z$ positive), the sign of the whole expression depends on the second term, which can in general be negative, zero, or positive.
In particular, if $\hat Z > Z$ (equivalently $\hat\varepsilon < \varepsilon \Rightarrow \varepsilon_{\hat Y} < \varepsilon_Y$), then $\log(Z/\hat Z) < 0$ and the second term becomes negative.
This suggests that, if the predicted model is such that $p$ and $q$ are identical (or very close), then the KL term becomes sufficiently small, and the overall difference could potentially become negative.

We now construct an explicit counterexample where the true survival and censoring distributions have residual tail mass (so $\varepsilon>0$) and show numerically that $\mathcal{R}(\hat S) < \mathcal{R}(S_Y)$ for a misspecified $\hat S$.

Let $Y$ and $C$ be independent with identical mixture distributions:
\begin{align*}
S_Y(t) = 0.2 + 0.8 e^{-t}, \qquad f_Y(t) = 0.8 e^{-t},\\
S_C(t) = 0.2 + 0.8 e^{-t}, \qquad f_C(t) = 0.8 e^{-t}.
\end{align*}

Then
\[
\varepsilon = P(Y=\infty)P(C=\infty) = \varepsilon_Y \cdot \varepsilon_C =0.2 \cdot 0.2 = 0.04 > 0,
\qquad Z = 1 - \varepsilon = 0.96.
\]

Consider a misspecified model with smaller residual mass $\varepsilon_{\hat Y} = 0.1$:
\[
\hat S(t) = 0.1 + 0.9 e^{-t}, \qquad \hat f(t) = 0.9 e^{-t}.
\]

then
\[
\hat \varepsilon = \varepsilon_{\hat Y} \cdot \varepsilon_C = 0.1 \cdot 0.2 = 0.02 < \varepsilon,
\qquad \hat Z = 1 - \hat \varepsilon = 0.98 > Z
\]

The risk difference \eqref{eq:rcll_exp_diff} is:
\[
I := \mathcal{R}(\hat S)-\mathcal{R}(S_Y) = \int_0^\infty\!\Big[ f_Y(t)S_C(t)\log\frac{f_Y(t)}{\hat f(t)} + f_C(t) S_Y(t) \log\frac{S_Y(t)}{\hat S(t)}\Big] dt.
\]

Substituting the expressions and setting $x = e^{-t} \in (0,1]$ and therefore $dt = -dx/x$, yields
\[
I = \int_0^1 0.8(0.2+0.8x)\left[\log\frac{8}{9}+\log\frac{0.2+0.8x}{0.1+0.9x}\right] dx.
\]

Numerical evaluation of this integral gives
\[
I \approx -0.0117 < 0.
\]

Thus $\mathcal{R}(\hat S)-\mathcal{R}(S_Y) < 0$, i.e. the misspecified model yields a lower expected loss than the true model.
This confirms that the RCLL scoring rule is \emph{not proper} following \Cref{def:surv_proper} under independent censoring when $\varepsilon>0$.

\medskip\noindent
\textbf{Summary.} Under the independence assumption between $Y$ and $C$, the scoring rule $L_{RCLL}$ is strictly proper (following Definition~\ref{def:surv_proper} of marginal properness) if and only if at least one of the true survival functions $S_Y$ or $S_C$ has no residual mass at infinity, i.e., $\varepsilon = P(Y=\infty,C=\infty)=0$.
If $\varepsilon>0$ (both $Y$ and $C$ have positive probability of being infinite), then the scoring rule is improper, because a predicted model with a smaller residual mass $\varepsilon_{\hat Y} < \varepsilon_Y$ can yield a strictly lower expected loss than the true model, as demonstrated by the counterexample above.
\end{proof}

\subsection{Properness of SBS}
\label{proof:sbs-proof}

\begin{proof}
Let $\tau^*\in\mathcal T \subseteq \NNReals$ be a fixed positive evaluation time with $S_C(\tau^*)>0$.
We assume that functions $f_C,f_Y,\hat f,S_C,S_Y,\hat S$ are continuous on $[0,\infty)$ and non-negative.

Taking $L_{SBS}$ \eqref{eq:SBS} and assuming independence of $Y$ and $C$, the expected Survival Brier Score is:

\begin{equation}
\label{eq:exp_sbs}
\footnotesize
\begin{aligned}
\underset{(T,\Delta)}{\EE}\Big[L_{SBS}(\hatS, T, \Delta)\Big]
&= \int^\infty_0 \Big [  f_Y(t)S_C(t) L_{SBS}(\hatS,t,\delta=1) + f_C(t)S_Y(t) L_{SBS}(\hatS,t,\delta=0) \Big ] \, dt && \text{\Cref{prop:exp}} \\
&= \int^\infty_0 \Big [ f_Y(t)S_C(t) \left( \frac{\hatS^2(\tau^*) \II(t \leq \tau^*)}{S_C(t)} + \frac{\hatF^2(\tau^*) \II(t > \tau^*)}{S_C(\tau^*)}\right) + \\ 
&f_C(t)S_Y(t) \left(\frac{\hatF^2(\tau^*) \II(t > \tau^*)}{S_C(\tau^*)}\right) \Big ] \ dt && \text{def. SBS \eqref{eq:SBS}} \\
&= \hatS^2(\tau^*) \int^{\tau^*}_0 f_Y(t)  \ dt + \frac{\hatF^2(\tau^*)}{S_C(\tau^*)} \int^\infty_{\tau^*} f_Y(t)S_C(t) \ dt + \frac{\hatF^2(\tau^*)}{S_C(\tau^*)} \int^\infty_{\tau^*} f_C(t)S_Y(t)  \ dt 
&& \text{split integrals} \\
&= \hatS^2(\tau^*)F_Y(\tau^*) + \frac{\hatF^2(\tau^*)}{S_C(\tau^*)} \int^\infty_{\tau^*} \Big [ f_Y(t)S_C(t) +  f_C(t)S_Y(t) \Big ] \ dt 
&& \text{def. $F_Y$}
\\
\end{aligned}
\end{equation}

Note the following relationship (product rule, assuming $S_Y,S_C$ are differentiable):

\begin{equation}
\label{eq:product_rule}
\frac{d}{dt} \Big( S_Y(t)S_C(t) \Big)= \frac{dS_Y(t)}{dt} S_C(t) + S_Y(t) \frac{dS_C(t)}{dt} = -\Big( f_Y(t)S_C(t) + S_Y(t)f_C(t) \Big)    
\end{equation}

Substituting the above expression into \eqref{eq:exp_sbs} yields:

\[
\footnotesize
\begin{aligned}
\underset{(T,\Delta)}{\EE}\Big[L_{SBS}(\hat S, T, \Delta)\Big]
&= \hat S^2(\tau^*)F_Y(\tau^*) - \frac{\hat F^2(\tau^*)}{S_C(\tau^*)} \int^\infty_{\tau^*} \frac{d}{dt} \Big(S_Y(t)S_C(t)\Big) \ dt \\
&= \hat S^2(\tau^*)F_Y(\tau^*) - \frac{\hat F^2(\tau^*)}{S_C(\tau^*)} \Big[S_Y(t)S_C(t)\Big]^\infty_{t=\tau^*} \\
&= \hat S^2(\tau^*)F_Y(\tau^*) - \frac{\hat F^2(\tau^*)}{S_C(\tau^*)} \Big(\varepsilon_Y \varepsilon_C - S_Y(\tau^*)S_C(\tau^*)\Big)
&& \text{let } \varepsilon_Y := \lim_{t\to\infty} S_Y(t) \ge 0, \, \varepsilon_C := \lim_{t\to\infty} S_C(t) \ge 0 \\
&= \hat S^2(\tau^*)F_Y(\tau^*)
+ \hat F^2(\tau^*)\Big(S_Y(\tau^*) - \frac{\varepsilon}{S_C(\tau^*)}\Big)\;
&& \text{let } \varepsilon := \varepsilon_Y \varepsilon_C \ge 0
\end{aligned}
\]

So the expectation of $L_{SBS}$ at time $\tau^*$ depends only on the predicted value $\hatS(\tau^*)$.
For notational brevity set $x := \hat S(\tau^*), a := S_Y(\tau^*) \in [0,1], b := S_C(\tau^*)>0$.
Then:
\begin{equation}
\label{eq:sbs-expectation}
\mathcal{G}(x) := \EE[L_{SBS}] = x^2(1-a) + (1-x)^2\Big(a-\frac{\varepsilon}{b}\Big)
\end{equation}

The first derivative of the expected score is:

\[
\begin{aligned}
\mathcal{G}'(x) 
&= 2x(1-a) - 2(1-x)\Big(a-\frac{\varepsilon}{b}\Big) \\
&= 2x\Big(1-a+a-\frac{\varepsilon}{b}\Big)-2\Big(a - \frac{\varepsilon}{b} \Big) \\
&= 2x\Big(1-\frac{\varepsilon}{b}\Big)-2\Big(a-\frac{\varepsilon}{b}\Big)
\end{aligned}
\]

and setting $\mathcal{G}'(x)=0$ to get the stationary point yields (assuming $\varepsilon\neq b$):

\begin{equation}
\label{eq:SBS_minimum}
\boxed{\,x^\star = \frac{a - \varepsilon/b}{1 - \varepsilon/b}\, }
\end{equation}

We note that $\lim\limits_{t\to\infty} S_C(t) \le S_C(\tau^*) \Rightarrow \varepsilon_C \le b$ for finite $\tau^*$ (similarly $\varepsilon_Y \le a$) and since $0 \le \varepsilon_Y \le 1$, it follows that

\begin{equation}
0 \le \varepsilon_C \le b \underset{\varepsilon_Y \ge 0}{\Rightarrow}
0 \le \varepsilon_Y \cdot \varepsilon_C \le \varepsilon_Y \cdot b \underset{\varepsilon_Y \le 1}{\le} 1 \cdot b \Rightarrow 
0 \le \varepsilon \le b \Rightarrow 
0 \le 1 - \frac{\varepsilon}{b} \le 1
\end{equation}

with equality $\varepsilon = b$ corresponding to $\varepsilon_Y=1, \varepsilon_C=b$, which is a degenerate case with zero events, flat censoring survival after $\tau^*$ and thus an uninformative survival process.
Therefore for realistic survival data, $\varepsilon < b$, implying that the denominator in \eqref{eq:SBS_minimum} satisfies $1 - \varepsilon/b > 0$.
Moreover, since the stationary point must represent a valid survival probability, the numerator of \eqref{eq:SBS_minimum} must be non-negative, i.e. $a-\varepsilon/b \ge 0 \Rightarrow ab \ge \varepsilon \Rightarrow ab \ge \varepsilon_Y \varepsilon_C$, which is true by definition (time point $\tau^*$ is finite).

The second derivative of $\mathcal{G}(x)$ is
\begin{equation}
\label{eq:SBS_2deriv}
\mathcal{G}''(x) = 2\Big(1-\frac{\varepsilon}{b}\Big) > 0.
\end{equation}

Hence the expected loss $\mathcal{G}$ is strictly convex and has a unique global minimum.
We distinguish between two cases, based on the value of $\varepsilon$:

\medskip\noindent\textbf{Case 1.} If $\varepsilon = \varepsilon_Y \varepsilon_C = 0$ (at least one tail vanishes — e.g., $\lim\limits_{t\to\infty} S_Y(t)=0$ or $\lim\limits_{t\to\infty} S_C(t)=0$), then \eqref{eq:SBS_minimum} and \eqref{eq:SBS_2deriv} become:

\begin{equation}
\label{eq:sbs-proper-solution}
x^\star = a, \quad \mathcal{G}''(x) = 2 > 0.    
\end{equation}

Thus the unique global minimum is found at $\hat S(\tau^*) = S_Y(\tau^*)$, i.e., $L_{SBS}$ is \emph{strictly proper} at $\tau^*$.

\medskip\noindent\textbf{Case 2.} If $0 < \varepsilon < b$ (both tails nonzero):
\begin{equation}
\label{eq:sbs-improper-min}
x^\star = \frac{a - \varepsilon/b}{1 - \varepsilon/b} \neq a, \quad \mathcal{G}''(x) > 0.    
\end{equation}
The expected loss remains convex with a unique minimum, but the minimum is shifted:
\begin{equation}
\label{eq:SBS_shift}
x^\star - a = -\frac{\varepsilon(1-a)}{b-\varepsilon} \le 0    
\end{equation}
since $a\in[0,1] \Rightarrow 1-a \ge 0$ and $b-\varepsilon > 0, \varepsilon > 0$.
Therefore $x^\star \le a$ and the predicted optimum \textit{underestimates} the true survival probability.
In this case, $L_{SBS}$ is \emph{improper}. 
Note that when $a = S_Y(\tau^*) = 1$, from \eqref{eq:SBS_minimum} we obtain $x^\star = 1$, so the minimizer coincides with the truth.
This is another degenerate observational case where $L_{SBS}$ is proper: no events occur on $[0,\tau^*]$ (everyone survives past $\tau^*$), so the score is minimized at the trivial prediction $\hat S(\tau^*)=1$.

\medskip\noindent
\textbf{Summary.} Under the independence assumption between $Y$ and $C$, and following Definition~\ref{def:surv_proper} of marginal properness, $L_{SBS}$ is strictly proper at each fixed $\tau^*$ with $S_C(\tau^*)>0$, whenever at least one of the survival functions $S_Y$ or $S_C$ is proper, i.e., has a vanishing tail ($\varepsilon=0$).
If both tails persist ($\varepsilon>0$), the expected loss remains strictly convex but attains its minimum at a misspecified prediction $x^\star < S_Y(\tau^*)$, which underestimates the true survival probability; consequently, $L_{SBS}$ is not proper in this case. 
The only exception is the degenerate scenario $S_Y(\tau^*)=1$, where no events occur before $\tau^*$ and the trivial prediction $\hat S(\tau^*)=1$ is optimal, rendering $L_{SBS}$ strictly proper again.
\end{proof}

\begin{remark}[First-order approximation of the improperness bias]
\label{rem:sbs_bias_approx}

We consider a first-order expansion of the deviation derived in \eqref{eq:SBS_shift} (replacing back $a = S_Y(\tau^*)$ and $b = S_C(\tau^*)$):

\begin{equation}
\label{eq:sbs_shift_2}
x^\star - S_Y(\tau^*) = -\frac{\varepsilon\bigl(1 - S_Y(\tau^*)\bigr)}{S_C(\tau^*) - \varepsilon}.
\end{equation}

For a fixed $\tau^*$ and small residual tail mass $\varepsilon > 0$, we expand the denominator as

\begin{equation}
\frac{1}{S_C(\tau^*) - \varepsilon} 
= \frac{1}{S_C(\tau^*)} \cdot \frac{1}{1 - \varepsilon/S_C(\tau^*)}
= \frac{1}{S_C(\tau^*)} \left(1 + \frac{\varepsilon}{S_C(\tau^*)} + O(\varepsilon^2)\right).
\end{equation}

Substituting this expansion in \eqref{eq:sbs_shift_2} yields

\begin{equation}
x^\star - S_Y(\tau^*) 
= -\varepsilon\bigl(1 - S_Y(\tau^*)\bigr) \cdot \left[ \frac{1}{S_C(\tau^*)} + \frac{\varepsilon}{S_C(\tau^*)^2} + O(\varepsilon^2) \right].
\end{equation}

Retaining only the first-order term in $\varepsilon$, we obtain the approximate bias

\begin{equation}
\boxed{
x^\star - S_Y(\tau^*) \;\approx\; -\; \varepsilon \;\cdot\; \frac{1 - S_Y(\tau^*)}{S_C(\tau^*)}
}
\end{equation}

which is strictly negative for $\varepsilon > 0$ and vanishes as $\varepsilon \to 0$.
\end{remark}

\subsection{Properness of ISBS}
\label{proof:isbs-proof}

\begin{proof}

Let $\tau^*\in\mathcal T \subseteq \NNReals$ be a fixed finite time point with $S_C(\tau)>0$ for all $\tau \in [0,\tau^*]$.
We assume that functions $f_C,f_Y,\hat f,S_C,S_Y,\hat S$ are continuous on $[0,\infty)$ and non-negative.

Recall that the ISBS \eqref{eq:ISBS} is the time integral of the pointwise in-$\tau$ SBS \eqref{eq:SBS}:
\[
L_{ISBS}(\hat S,t,\delta,\tau^*)
:= \int^{\tau^*}_0 L_{SBS}(\hat S, t, \delta, \tau) \ d\tau
\]

Under the independence assumption ($Y\indep C$), we calculate the expectation of the ISBS using Proposition~\ref{prop:exp}, and observe that the result expands to a double integral over $(t,\tau)$:

\begin{equation}
\underset{(T,\Delta)}{\EE}\Big[L_{ISBS}(\hatS, T, \Delta)\Big] 
= \underset{(T,\Delta)}{\EE}\!\left[\int_0^{\tau^*} L_{SBS}(\hat S, t, \delta, \tau)\,d\tau\right]
= \int_{0}^{\infty}\int_{0}^{\tau^*} \mathcal{K}(\tau,t)\,d\tau\,dt
\end{equation}

where $\mathcal{K}(\tau,t)$ denotes the nonnegative integrand obtained by substituting $\phi=L_{\mathrm{SBS}}(\hat S,t,\delta, \tau)$ into the right-hand side of Equation \eqref{eq:expprop} (same as the integrand in the second line of \eqref{eq:exp_sbs}, using $\tau$ instead of $\tau^*$).
Because $\mathcal{K}(\tau,t)\ge 0$ for all $\tau\in[0,\tau^*]$, $t\ge0$ (each SBS term is a sum of squared terms divided by $S_C(\cdot)>0$), Tonelli's theorem permits swapping the order of integration. Swapping yields

\begin{equation}
\label{eq:E_ISBS_swap}
\begin{aligned}
\underset{(T,\Delta)}{\EE}\Big[L_{ISBS}(\hat S, T, \Delta)\Big] 
&= \int_{0}^{\tau^*}\!\left(\int_{0}^{\infty} \mathcal{K}(\tau,t)\,dt\right)\,d\tau \\
&= \int_0^{\tau^*} \underset{(T,\Delta)}{\EE}\Big[L_{SBS}(\hat S, T, \Delta)\Big]\,d\tau.
\end{aligned}
\end{equation}

From \eqref{eq:E_ISBS_swap}, the expected ISBS is the time integral of expected SBS losses evaluated pointwise in $\tau$.
Therefore, the same derivation as in the proof in \Cref{proof:sbs-proof}—determining $\mathcal{G}(x)=\EE[L_{SBS}]$, analyzing its curvature $\mathcal{G}''(x)$, and identifying the stationary point $x^\star$ for $x=\hat S(\tau)$—applies directly for each $\tau\in[0,\tau^*]$.
Briefly, for fixed $\tau$, we set
\[
x := \hat S(\tau),\qquad a := S_Y(\tau)\in[0,1],\qquad b := S_C(\tau)>0,
\]
and denote the tail constants
\[
\varepsilon_Y := \lim_{t\to\infty} S_Y(t),\qquad
\varepsilon_C := \lim_{t\to\infty} S_C(t),\qquad
\varepsilon := \varepsilon_Y\varepsilon_C \ge 0.
\]

Then the pointwise-in-$\tau$ expected SBS equals \eqref{eq:sbs-expectation}:

\[
\begin{aligned}
\underset{(T,\Delta)}{\EE}\Big[L_{SBS}(\hat S, T, \Delta)\Big]
&= \hat S^2(\tau)(1-S_Y(\tau)) + (1-\hat S(\tau))^2\Big(S_Y(\tau)-\frac{\varepsilon}{S_C(\tau)}\Big) \\
&= x^2(1-a) + (1-x)^2\Big(a-\frac{\varepsilon}{b}\Big) := \mathcal{G}_\tau(x,a,b)
\end{aligned}
\]

The expected value of the ISBS is therefore (with $x,a,b$ all functions of $\tau$):
\begin{equation}
\label{eq:ISBS_exp_withG}
\underset{(T,\Delta)}{\EE}[L_{ISBS}(\hat S, T, \Delta)] = \int_0^{\tau^*} \mathcal{G}_\tau (x, a, b)\,d\tau.
\end{equation}

Since the integrand $\mathcal{G}\tau$ depends on $\hat S$ only through its value at $\tau$, minimizing the integrated expectation over all $\hat S(\cdot)$ reduces to minimizing each pointwise expected loss $\mathcal{G}\tau(x,a,b)$ for fixed $\tau$ independently.

We distinguish between the following two cases:

\medskip\noindent\textbf{Case 1: } $\varepsilon=0$ (i.e. at least one tail vanishes: $\varepsilon_Y=0$ or $\varepsilon_C=0$).

\begin{equation}
\label{eq:ISBS_proper_exp}
\underset{(T,\Delta)}{\EE}\Big[L_{ISBS}(\hat S, T, \Delta)\Big] = \int_0^{\tau^*} \Big[\hat S^2(\tau)(1-S_Y(\tau)) + (1-\hat S(\tau))^2 S_Y(\tau)\Big]\,d\tau.
\end{equation}

From \eqref{eq:sbs-proper-solution}, for every fixed $\tau$, the pointwise expected SBS is uniquely minimized at $x = S_Y(\tau)$.
Integrating over $\tau$ therefore yields that the expected ISBS is uniquely minimized when $\hat S(\tau) = S_Y(\tau)$ for all $\tau \in [0, \tau^*]$.
Therefore, under independent censoring, $L_{ISBS}$ is strictly proper whenever $\varepsilon = 0$.
This result aligns with the known propriety of IBS in the absence of censoring \citep{Gneiting2007}, since \eqref{eq:ISBS_proper_exp} coincides with the IBS expectation in that setting.

\medskip\noindent\textbf{Case 2: } $\varepsilon>0$ (i.e. both tails persist).

From the SBS analysis \eqref{eq:SBS_minimum}–\eqref{eq:SBS_shift}, the pointwise expected loss $\mathcal{G}_\tau(x,a,b)$ remains strictly convex with a unique minimizer \eqref{eq:sbs-improper-min}:

\begin{equation}
\label{eq:ISBS_improper_min}
x^\star = \frac{a - \varepsilon/b}{1 - \varepsilon/b} \le a    
\end{equation}

implying that the expected SBS at each $\tau$ is minimized by a prediction $x^\star = \hat S(\tau)$ which is smaller or equal to the true survival $a = S_Y(\tau)$.
Since survival functions are monotonically decreasing, equality across all $\tau \in [0,\tau^*]$ would imply the degenerate case $S_Y(\tau) \equiv 1$, corresponding to no observed events up to $\tau^*$, where \eqref{eq:ISBS_improper_min} yields $x^\star = 1$ and thus a trivially proper minimizer.
Otherwise, if strict inequality holds for any $\tau$, the integrand in \eqref{eq:ISBS_exp_withG} attains a smaller value for this misspecified $\hat S(\tau)$ than for the true survival $S_Y(\tau)$.
Consequently, the overall expected ISBS is strictly lower for the misspecified survival curve, demonstrating that the loss is improper whenever $\varepsilon > 0$.

\medskip\noindent
\textbf{Summary.} Under $Y\indep C$, and following Definition~\ref{def:surv_proper} of marginal properness, ISBS is \emph{strictly proper} on $[0,\tau^*]$ whenever at least one of the survival functions $S_Y$ or $S_C$ is proper (i.e.\ has a vanishing tail, $\varepsilon=0$).
If both tails persist ($\varepsilon>0$), the expected ISBS can be reduced by a misspecified prediction $\hat S(\tau)<S_Y(\tau)$, leading to systematic underestimation of survival and rendering $L_{ISBS}$ improper.
The only exception is the degenerate scenario $S_Y(\tau)=1$, where no events occur before $\tau^*$, and the trivial prediction $\hat S(\tau)=1$ remains optimal.
\end{proof}

\section{Experiments}
\label{sec:appendixC}

\begin{Prop}[\textbf{Expected Survival at the Sample Maximum}]
\label{prop:s_tmax}
Let $T_1, \dots, T_n$ be i.i.d. samples from a continuous distribution with survival function $S(t) = P(T > t)$, and let $M_n = \max\{T_1, \dots, T_n\}$. Then
\begin{equation}
\mathbb{E}[S(M_n)] = \frac{1}{n+1}.
\end{equation}
Consequently, as $n \to \infty$, $\mathbb{E}[S(M_n)] \to 0$.
\end{Prop}

\begin{proof}
Let $F(t) = 1 - S(t)$ denote the cumulative distribution function of $T$.  
By the probability integral transform, the variables
\[
U_i = F(T_i), \quad i=1,\dots,n,
\]
are independent and uniformly distributed on $(0,1)$.

The sample maximum satisfies
\[
F(M_n) = \max\{F(T_1), \dots, F(T_n)\} = \max\{U_1, \dots, U_n\}.
\]
Let $U_{(n)} = \max(U_1, \dots, U_n)$ denote the largest order statistic.
It is a standard result for uniform order statistics \cite{David2004OrderStatistics}
that
\[
U_{(n)} \sim \text{Beta}(n,1).
\]
Therefore,
\[
1 - U_{(n)} \sim \text{Beta}(1,n).
\]

For a Beta$(a,b)$ random variable, the expectation equals $a/(a+b)$. Hence,
\[
\mathbb{E}[1 - U_{(n)}] = \frac{1}{1+n}.
\]

Finally, since
\[
S(M_n) = 1 - F(M_n) = 1 - U_{(n)},
\]
we obtain
\[
\mathbb{E}[S(M_n)] = \frac{1}{n+1}.
\]

The limit statement follows immediately, as $1/(n+1) \to 0$ when $n \to \infty$.
\end{proof}

\begin{table}[ht]
\centering
\scriptsize
\caption{Empirical violations of properness for the SBS and ISBS using estimated censoring distribution $\hat{G}(t)$ via Kaplan--Meier across $K = 10{,}000$ simulations.
SBS was evaluated at the 10th percentile ($q_{0.1}$), median ($q_{0.5}$), and 90th percentile ($q_{0.9}$) of observed times.
ISBS was computed using 50 equidistant time points between the 5th and 80th percentiles of observed times.
Columns report the number of violations (\#Viol.), violation rate (Rate), and 
the mean score difference (true minus predicted) among simulations flagged as violations ($\bar{D}_{L}^{\text{Viol}}$).}
\label{tab:sbs_violations_est}
\renewcommand{\arraystretch}{1.2}
\begin{tabularx}{\textwidth}{c|CCC|CCC|CCC|CCC}
\toprule
$n$ & \multicolumn{3}{c|}{\textbf{SBS}($q_{0.1}$)} & \multicolumn{3}{c|}{\textbf{SBS}($q_{0.5}$)} & \multicolumn{3}{c|}{\textbf{SBS}($q_{0.9}$)} & \multicolumn{3}{c}{\textbf{ISBS}} \\
 & \#Viol. & Rate & $\bar{D}_{\text{SBS}}^{\text{Viol}}$ & \#Viol. & Rate & $\bar{D}_{\text{SBS}}^{\text{Viol}}$ & \#Viol. & Rate & $\bar{D}_{\text{SBS}}^{\text{Viol}}$ & \#Viol. & Rate & $\bar{D}_{\text{ISBS}}^{\text{Viol}}$ \\
\midrule
10    & 4325 & 0.432  & 0.00528 & 695  & 0.0695 & 0.00632 & 780 & 0.078 & 0.00962 & 320 & 0.032 & 0.00145 \\
25    & 442  & 0.0442 & 0.00115 & 374  & 0.0374 & 0.00355 & 414 & 0.0414 & 0.00401 & 59  & 0.0059 & 0.000360 \\
50    & 378  & 0.0378 & 0.000695 & 261 & 0.0261 & 0.00171 & 320 & 0.032 & 0.00249 & 21 & 0.0021 & 0.000176 \\
100   & 186  & 0.0186 & 0.000361 & 200 & 0.020  & 0.000895 & 243 & 0.0243 & 0.00134 & 0  & 0      & --       \\
250   & 53   & 0.0053 & 0.000182 & 120 & 0.012  & 0.000407 & 153 & 0.0153 & 0.000578 & 0  & 0      & --       \\
500   & 15   & 0.0015 & 0.000130 & 51  & 0.0051 & 0.000276 & 89  & 0.0089 & 0.000327 & 0  & 0      & --       \\
1000  & 0    & 0      & --       & 19  & 0.0019 & 0.000139 & 36  & 0.0036 & 0.000226 & 0  & 0      & --       \\
2500  & 0    & 0      & --       & 0   & 0      & --       & 8   & 0.0008 & 0.000145 & 0  & 0      & --       \\
5000  & 0    & 0      & --       & 0   & 0      & --       & 0   & 0      & --       & 0  & 0      & --       \\
10000 & 0    & 0      & --       & 0   & 0      & --       & 0   & 0      & --       & 0  & 0      & --       \\
\bottomrule
\end{tabularx}
\end{table}

\begin{table}[ht]
\scriptsize
\centering
\caption{RCLL scores across 100 Monte Carlo replicates for varying time grid sizes $B$ in the high-censoring task with $n_{\mathrm{test}} = 250$.
Values are reported as mean $\pm$ standard deviation.
Grid sizes correspond to proportions of the full set of unique event times ($B = 363$).
`Truth (true $S$, true $f$)' denotes the exact right-censored likelihood.}
\label{tab:rcll_grid_high_cens}
\renewcommand{\arraystretch}{1.1}
\begin{tabularx}{\textwidth}{>{\raggedright\arraybackslash}l C C C C C C C}
\toprule
Model & $B = 8$ (2\%) & $B = 19$ (5\%) & $B = 37$ (10\%) & $B = 73$ (20\%) & $B = 182$ (50\%) & $B = 363$ (100\%) \\
\midrule
Truth (true $S$, true $f$)
& \makecell{0.681 $\pm$ 0.055}
& \makecell{0.688 $\pm$ 0.046}
& \makecell{0.680 $\pm$ 0.052}
& \makecell{0.679 $\pm$ 0.047}
& \makecell{0.682 $\pm$ 0.049}
& \makecell{0.679 $\pm$ 0.049} \\
Truth (true $S$, interp $f$)
& \makecell{0.733 $\pm$ 0.081}
& \makecell{0.724 $\pm$ 0.064}
& \makecell{0.694 $\pm$ 0.061}
& \makecell{0.687 $\pm$ 0.057}
& \makecell{0.684 $\pm$ 0.050}
& \makecell{0.682 $\pm$ 0.051} \\
Oracle
& \makecell{0.737 $\pm$ 0.078}
& \makecell{0.724 $\pm$ 0.060}
& \makecell{0.698 $\pm$ 0.063}
& \makecell{0.698 $\pm$ 0.061}
& \makecell{0.709 $\pm$ 0.053}
& \makecell{0.739 $\pm$ 0.063} \\
LLogis
& \makecell{0.748 $\pm$ 0.086}
& \makecell{0.734 $\pm$ 0.065}
& \makecell{0.704 $\pm$ 0.065}
& \makecell{0.704 $\pm$ 0.062}
& \makecell{0.714 $\pm$ 0.054}
& \makecell{0.744 $\pm$ 0.063} \\
Weibull
& \makecell{0.767 $\pm$ 0.085}
& \makecell{0.755 $\pm$ 0.070}
& \makecell{0.729 $\pm$ 0.071}
& \makecell{0.730 $\pm$ 0.071}
& \makecell{0.736 $\pm$ 0.058}
& \makecell{0.765 $\pm$ 0.063} \\
LogNorm\_scale
& \makecell{0.785 $\pm$ 0.084}
& \makecell{0.760 $\pm$ 0.065}
& \makecell{0.733 $\pm$ 0.061}
& \makecell{0.733 $\pm$ 0.059}
& \makecell{0.743 $\pm$ 0.053}
& \makecell{0.776 $\pm$ 0.064} \\
LogNorm\_int
& \makecell{0.812 $\pm$ 0.089}
& \makecell{0.795 $\pm$ 0.068}
& \makecell{0.776 $\pm$ 0.064}
& \makecell{0.771 $\pm$ 0.060}
& \makecell{0.791 $\pm$ 0.062}
& \makecell{0.812 $\pm$ 0.069} \\
LogNorm
& \makecell{0.846 $\pm$ 0.091}
& \makecell{0.822 $\pm$ 0.066}
& \makecell{0.805 $\pm$ 0.063}
& \makecell{0.799 $\pm$ 0.060}
& \makecell{0.818 $\pm$ 0.062}
& \makecell{0.841 $\pm$ 0.069} \\
Cox\_int
& \makecell{0.927 $\pm$ 0.154}
& \makecell{0.880 $\pm$ 0.100}
& \makecell{0.849 $\pm$ 0.074}
& \makecell{0.857 $\pm$ 0.077}
& \makecell{0.898 $\pm$ 0.077}
& \makecell{0.937 $\pm$ 0.068} \\
KM
& \makecell{0.941 $\pm$ 0.147}
& \makecell{0.909 $\pm$ 0.091}
& \makecell{0.884 $\pm$ 0.070}
& \makecell{0.891 $\pm$ 0.074}
& \makecell{0.932 $\pm$ 0.074}
& \makecell{0.978 $\pm$ 0.071} \\
RSF
& \makecell{1.246 $\pm$ 0.226}
& \makecell{1.324 $\pm$ 0.129}
& \makecell{1.363 $\pm$ 0.151}
& \makecell{1.414 $\pm$ 0.143}
& \makecell{1.464 $\pm$ 0.138}
& \makecell{1.480 $\pm$ 0.139} \\
\bottomrule
\end{tabularx}
\end{table}

\begin{figure}[h]
\centering
\includegraphics[width=1.0\textwidth]{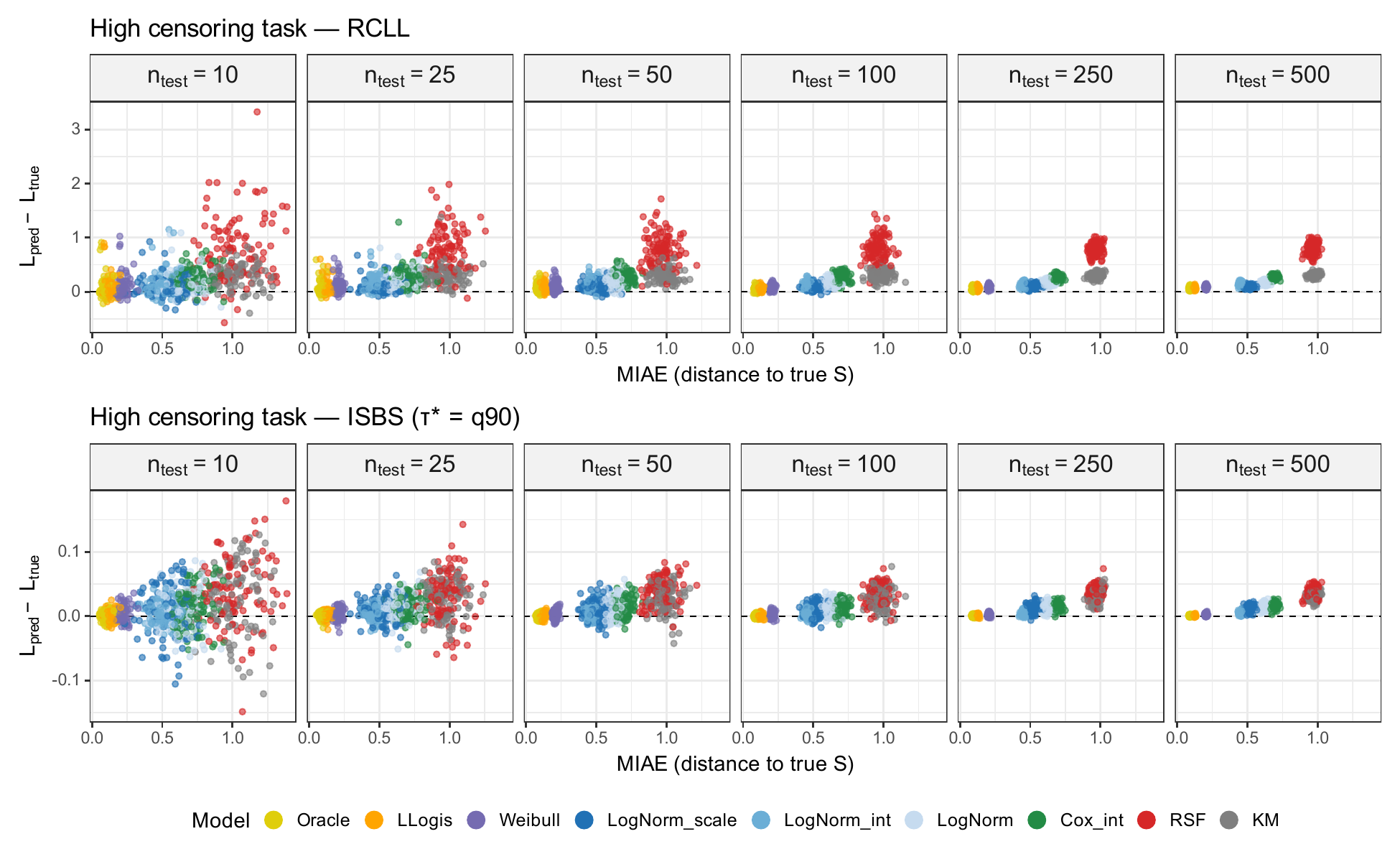}
\caption{Effect of test set size on sensitivity to model misspecification under high censoring for RCLL (top) and ISBS (bottom).
Each panel corresponds to a test set size $n_{\mathrm{test}} \in \{10, 25, 50, 100, 250, 500\}$.
Points represent 100 Monte Carlo replicates and show the excess loss $\Delta L = L_{\mathrm{pred}}(\hat S) - L_{\mathrm{true}}(S_Y)$ versus the mean integrated absolute error (MIAE) between $\hat S$ and $S_Y$.
Models are ordered along the horizontal axis by increasing discrepancy from the true survival function (larger MIAE).
ISBS is integrated up to the 90th percentile of observed times.}
\label{fig:sens_vary_n_high_cens}
\end{figure}

\end{document}